\input ./style/arxiv-general.cfg
\documentclass[MSNbibl,number,seceqn,citesort,dvips]{arxbj}
\makeatletter
   \@ifpackageloaded{graphicx}{}{\usepackage{graphicx}}
\makeatother
\usepackage{upgreek}

\aid{0}
\volume{21}
\issue{4}
\pubyear{2015}
\firstpage{2024}
\lastpage{2072}
\doi{10.3150/14-BEJ634} 
\docsubty{FLA}

\makeatletter
\newcommand{\ds}{\displaystyle}
\newproclaim{definition}{Definition}[section]
\newtheorem{lemma}{Lemma}[section]
\newtheorem{proposition}{Proposition}[section]
\newtheorem{theorem}{Theorem}[section]
\newremark{step}{Step}
\newremark{remark}{Remark}[section]
\newtheorem{corollary}{Corollary}[section]
\makeatother

\begin{document}
\begin{frontmatter}

\title{Local asymptotic mixed normality property
for nonsynchronously observed diffusion~processes}
\runtitle{LAMN for nonsynchronously observed diffusion}

\begin{aug}
\author[A]{\inits{T.}\fnms{Teppei}~\snm{Ogihara}\corref{}\ead[label=e1]{ogihara@ism.ac.jp}}
\address[A]{The Institute of Statistical Mathematics, 10-3 Midori-cho, Tachikawa, Tokyo 190-8562, Japan. \printead{e1}}
\end{aug}

\received{\smonth{10} \syear{2013}}
\revised{\smonth{2} \syear{2014}}

%
\begin{abstract}
We prove the local asymptotic mixed normality (LAMN) property for a
family of probability measures
defined by parametrized diffusion processes with nonsynchronous observations.
We assume that observation times of processes are independent of processes
and we will study asymptotics when the maximum length of observation
intervals goes to zero in probability.
We also prove that the quasi-maximum likelihood estimator and the
Bayes-type estimator proposed
in Ogihara and Yoshida (\textit{Stochastic Process. Appl.}
\textbf{124} (2014) 2954--3008)
are asymptotically efficient.
\end{abstract}

%
\end{frontmatter}

\section{Introduction}

Given a probability space $(\Omega, \mathcal{F}, P)$ with a filtration
$\mathbf{F}=(\mathcal{F}_t)_{t\in[0,T]}$,
we consider a two-dimensional $\mathbf{F}$-adapted process $Y=\{Y_t\}
_{0\leq t\leq T}=\{(Y^1_t,Y^2_t)\}_{0\leq t\leq T}$
satisfying the following stochastic differential equation:
%
\begin{equation}
\label{SDE} \mathrm{d}Y_t=\mu(t,Y_t,
\sigma_{\ast})\,\mathrm{d}t+b(t,Y_t,\sigma_{\ast})
\,\mathrm{d}W_t,\qquad t \in[0,T],
\end{equation}
where\vspace*{1pt} $\{W_t\}_{0\leq t\leq T}$ is a two-dimensional standard $\mathbf{F}$-Wiener process,
$b=(b^{ij})_{1\leq i,j\leq2}\dvtx\break  [0,T]\times\mathbb{R}^2 \times\Lambda
\to\mathbb{R}^2\otimes\mathbb{R}^2$ is a Borel function,
$\mu=(\mu^1,\mu^2)$ is a $\mathbb{R}^2$-valued function, $\sigma_{\ast}
\in\Lambda$, and $\Lambda$ is a bounded open subset of $\mathbb{R}^d$.

We will consider the problem of estimating the unknown true value
$\sigma_{\ast}$ of the parameter
by nonsynchronous observations $\{Y^1_{S^{n,i}}\}_{i=0}^{\ell_{1,n}}$
and $\{Y^2_{T^{n,j}}\}_{j=0}^{\ell_{2,n}}$,
where $\{S^{n,i}\}_{i=0}^{\ell_{1,n}}$ and $\{T^{n,j}\}_{j=0}^{\ell
_{2,n}}$ are observation times of $Y^1$ and $Y^2$, respectively.

The problem of nonsynchronous observations appears when we study
statistical inference for high-frequency financial data.
Hayashi and Yoshida \cite{hay-yos01} pointed out that simple `\textit{synchronization}' methods such as linear interpolation or
`previous-tick' interpolation do not work well
for covariation estimation. They constructed a consistent estimator of
the quadratic covariation of processes.
On the other hand, Malliavin and Mancino \cite{mal-man} proposed an
estimator based on a~Fourier analytic method,
and Ogihara and Yoshida \cite{ogi-yos} constructed a quasi-maximum
likelihood estimator and a Bayes-type estimator
for a statistical model of nonsynchronously observed diffusion processes.
There are also several studies about covariation estimation under
nonsynchronous observations and market microstructure noise.
See Barndorff-Nielsen \textit{et al.} \cite{bar-et-al}, Christensen, Kinnebrock
and Podolskij \cite{chr-et-al}, A\"{\i}t-Sahalia, Fan and
Xiu \cite{ait-et-al}, Bibinger \textit{et al.} \cite{bib-et-al}, for example.

In this work, we will study the local asymptotic mixed normality (LAMN)
property of a statistical model of nonsynchronously observed diffusion
processes.
The definition of the LAMN property is as follows (Jeganathan \cite{jeg}).

%
\begin{definition}\label{LAMN-def}
Let $P_{\sigma,n}$ be a probability measure on some measurable space
$(\mathcal{X}_n,\mathcal{A}_n)$ for each $\sigma\in\Lambda$ and $n\in
\mathbb{N}$.
Then the family $\{P_{\sigma,n}\}_{\sigma,n}$ satisfies the local
asymptotic mixed normality (LAMN) property at $\sigma=\sigma_{\ast}$
if there exist a sequence $\{b_n\}_{n\in\mathbb{N}}$ of positive numbers,
$d\times d$ symmetric random matrices $\Gamma_n,\Gamma$ and
$d$-dimensional random vectors $\mathcal{N}_n,\mathcal{N}$
such that $\Gamma$ is positive definite a.s., $P_{\sigma_{\ast
},n}[\Gamma_n  \mbox{ is  positive  definite}]=1$  $(n\in\mathbb{N})$,
$b_n\to\infty$,
\[
\log\frac{\mathrm{d}P_{\sigma_{\ast}+b_n^{-1/2}u,n}}{\mathrm{d}P_{\sigma_{\ast},n}}- \biggl(u^{\star}\sqrt{\Gamma_n}
\mathcal{N}_n-\frac{1}{2}u^{\star}\Gamma _n
u \biggr)\to0
\]
in $P_{\sigma_{\ast},n}$-probability as $n\to\infty$ for any
$u\in\mathbb{R}^d$, where $\star$ represents transpose.
Moreover, $\mathcal{N}$ follows the $d$-dimensional standard normal
distribution, $\mathcal{N}$ is independent of $\Gamma$ and
$\mathcal{L}(\mathcal{N}_n,\Gamma_n|P_{\sigma_{\ast},n})\to\mathcal
{L}(\mathcal{N},\Gamma)$ as $n\to\infty$.
\end{definition}

The LAMN property is significantly related to asymptotic efficiency of
estimators.
Let $E_{\sigma}$ denote expectation with respect to $P_{\sigma,n}$.
Jeganathan \cite{jeg} proved the minimax theorem:
%
\begin{equation}
\label{minimax-ineq}
\lim_{\alpha\to\infty}\liminf_{n\to\infty}\sup
_{|u|\leq\alpha}E_{\sigma
_{\ast}+b_n^{-1/2}u}\bigl[l\bigl(\bigl|b_n^{1/2}
\bigl(V_n-\sigma_{\ast}-b_n^{-1/2}u
\bigr)\bigr|\bigr)\bigr] \geq E\bigl[l\bigl(\bigl|\Gamma^{-1/2}\mathcal{N}\bigr|\bigr)
\bigr]
\end{equation}
for any estimators $\{V_n\}_n$ and any function $l\dvtx [0,\infty
)\to[0,\infty)$ which is nondecreasing and \mbox{$l(0)=0$},
when the family $\{P_{\sigma,n}\}_{\sigma,n}$ has the LAMN
property at $\sigma=\sigma_{\ast}$.
This inequality gives lower bounds of risk functions of estimation
errors.
In particular, this inequality gives a lower bound of asymptotic
variance of estimators if $l(x)=x^2$.
When estimators $\{V_n\}_n$ attain the lower bound of (\ref{minimax-ineq}),  they are called \textit{asymptotically efficient}.

In a statistical model with independent identically distributed random
variables, the maximum likelihood estimator and the Bayes estimator
have minimal asymptotic variance under certain regularity conditions.
See Chapter I of Ibragimov and Has'minskii \cite{ibr-has03} for the details.
The LAMN property is proved for a statistical model of
one-dimensional
diffusion
process with synchronous, equispaced observations in Dohnal  \cite{doh}, and then the
results are extended to a multi-dimensional diffusion in  Gobet \cite{gob},
by using a Malliavin calculus approach.
On the other hand, Gobet \cite{gob02} proved the  LAN property (that means the LAMN
property with a deterministic $\Gamma$)
for ergodic diffusion process when the end time $T$ of observations
goes to infinity.

The aim of this paper is to show the LAMN property for nonsynchronously
observed diffusion processes,
and consequently have the minimax theorem (\ref{minimax-ineq}).
We also prove that the quasi-maximum likelihood estimator and the
Bayes-type estimator  proposed in Ogihara and Yoshida \cite{ogi-yos,ogi-yos02} are asymptotically efficient.
Ogihara and Yoshida \cite{ogi-yos} constructed an estimator
of quadratic covariation of the processes based on the quasi-maximum
likelihood estimator
and verified that the variance of estimation error of the estimator is much
smaller than that of the Hayashi--Yoshida estimator in a simple example.

When the observations occur in synchronous manner, the log-likelihood
ratio\break  $\log(\mathrm{d}P_{\sigma_{\ast}+b_n^{-1/2}u,n}/\mathrm{d}P_{\sigma_{\ast},n})$\vspace*{1pt}
is decomposed into differences of logarithms of transition density functions.
A Malliavin calculus approach enables us to apply limit theorems to
these differences, and consequently to obtain the LAMN property, as
seen in Gobet \cite{gob}.
However, when the sampling scheme is a nonsynchronous one, the
log-likelihood ratio does not have such a simple form and we cannot
apply the Malliavin calculus approach in Gobet \cite{gob} directly.
Instead, we will define stochastic processes that `\textit{connect}' the
process $Y$ and an Euler--Maruyama approximation process (Section~\ref{preliminaries}),
and we prove asymptotic equivalence of the log-likelihood ratio of $Y$
and that of Euler--Maruyama approximation.
Since the log-likelihood ratio of Euler--Maruyama approximation is
asymptotically equivalent to the quasi-log-likelihood ratio in Ogihara
and Yoshida \cite{ogi-yos}
and the quasi-log-likelihood ratio has a LAMN-type property, we obtain
the LAMN property of the model.

This paper is organized as follows. Section~\ref{main-results-section}
presents assumptions and main theorems.
Section~\ref{preliminaries} contains some preliminary results.
In Section~\ref{fundamental-section}, we introduce fundamental lemmas,
some notation and the result in Ogihara and Yoshida \cite{ogi-yos} with
respect to a LAMN-type property of the quasi-log-likelihood ratio.
Section~\ref{malliavin-section} gives some results in Malliavin calculus,
and Section~\ref{tightness-section} is devoted to prove tightness of
some log-likelihood ratios, which is used in the proof of the LAMN property.
We complete the proof of the main theorem in Section~\ref{LAMN-proof-section}.

\section{Main results}\label{main-results-section}

We begin with some general conventions.
For a real number $x$, $[x]$ denotes the maximum integer which is not
greater than $x$.
Let us denote by $|K|$ the length of interval $K$.
For a matrix $M$, $\Vert M \Vert$ represents the operator norm
of $M$
and $M^{\star}$ represents transpose of $M$.
Let $\mathcal{E}_l$ be the unit matrix of size $l$ and $\delta_{i,j}$
be Kronecker's delta function.
We denote $|x|^2=\sum_{i_1,\ldots, i_k}|x_{i_1,\ldots, i_k}|^2$ for
$x=\{x_{i_1,\ldots, i_k}\}_{i_1,\ldots, i_k}$.
For a vector $x=(x_1,\ldots, x_k)$,
we denote $\partial_{x}^l=(\frac{\partial^l}{\partial x_{i_1}\cdots\,
\partial x_{i_l}})_{i_1,\ldots, i_l=1}^k$.
We use the symbol $C$ for a generic positive constant varying from line
to line.
We denote by $\to^{s\mbox{-}\mathcal{L}}$ stable convergence of a
random sequence,
which is stronger than convergence in distribution and weaker
than convergence in probability.
See Aldous and Eagleson \cite{ald-eag} or Jacod \cite{jac} for the
definition and fundamental properties of stable convergence.

Let us start with some definitions and assumptions. The end time $T>0$
of observations is assumed to be a fixed constant.
We assume that the parameter space $\Lambda$ satisfies Sobolev's
inequality, that is, for any $p>d$, there exists $C>0$ such that
\[
\sup_{x\in\Lambda} \bigl|u(x)\bigr|\leq C\sum_{k=0,1}
\bigl\Vert\partial ^k_xu(x)\bigr\Vert_p\qquad
\bigl(u\in C^1(\Lambda)\bigr).
\]
It is the case if $\Lambda$ has a Lipschitz boundary (see Adams \cite{ada}, Adams and Fournier \cite{ada-fou}).

Let $\{\ell_{1,n}\}_{n\in\mathbb{N}}$ and $\{\ell_{2,n}\}_{n\in\mathbb
{N}}$ be sequences of positive integer-valued random variables,
the observation times $\Pi_n=((S^{n,i})_{i=0}^{\ell
_{1,n}},(T^{n,j})_{j=0}^{\ell_{2,n}})$ satisfy $S^{n,0}=T^{n,0}=0$,
$S^{n,\ell_{1,n}}=T^{n,\ell_{2,n}}=T$
and random times $\{S^{n,i}\}_i,\{T^{n,j}\}_j$ be monotone increasing
with respect to $i,j$.
Moreover, we assume that $\sigma(\{\Pi_n\}_n)$ is independent of $\{
(Y_t,W_t)\}_{0\leq t\leq T}$.
We assume that $\Pi_n$ and $Y_0$ do not depend on $\sigma_{\ast}$.

Let $b^k=(b^{k1},b^{k2})$ for $k=1,2$, where $\{b^{ij}\}_{i,j}$ are
elements of the diffusion coefficient~$b$.
Let $I^i=[S^{n,i-1},S^{n,i}), J^j=[T^{n,j-1},T^{n,j})$, $r_n=\max_{i,j}(|I^i|\vee|J^j|)$,
$\mathcal{E}^1(t)=\{\delta_{i,i'}1_{\{I^i\cap[0,t)\neq\varnothing\} }\}
_{i,i'=1}^{\ell_{1,n}}$,
$\mathcal{E}^2(t)=\{\delta_{j,j'}1_{\{J^j\cap[0,t)\neq\varnothing\} }\}
_{j,j'=1}^{\ell_{2,n}}$ for\vspace*{1pt} $t\in(0,T]$
and $G$ be an $\ell_{1,n}\times\ell_{2,n}$ matrix with the elements
$G_{ij}=|I^i\cap J^j||I^i|^{-1/2}|J^j|^{-1/2}$.
Moreover, let
\begin{eqnarray*}
\mathcal{U} &=& \bigl\{\bar{u}=\bigl(\bigl(s^i\bigr)_{i=0}^{L^1},
\bigl(t^j\bigr)_{j=0}^{L^2}\bigr);
L^1,L^2\in\mathbb{N},\\
&&\hspace*{5pt} 0=s^0<s^1<
\cdots<s^{L^1}=T, 0=t^0<t^1<\cdots
<t^{L^2}=T\bigr\},
\end{eqnarray*}
and we denote $\mathsf{X}_{\bar{u}}=((\mathsf{X}^1_{s^i})_{i=0}^{L^1},(\mathsf{X}^2_{t^j})_{j=0}^{L^2})$
and $\mathsf{X}_{\bar{v}}=((\mathsf{X}^1_{v^i})_{i=0}^{L},
(\mathsf{X}^2_{v^j})_{j=0}^{L})$ for a two-dimensional stochastic process
$\mathsf{X}=\{(\mathsf{X}^1_t,\mathsf{X}^2_t)\}_{0\leq t\leq T}$,
$\bar{u}=((s^i)_{i=0}^{L_1},(t^j)_{j=0}^{L_2})\in\mathcal{U}$ and $\bar
{v}=(v^i)_{i=0}^L$ satisfying $0=v^0<\cdots<v^L=T$.
Let $Y^{(\sigma)}=\{Y^{(\sigma)}_t\}_{0\leq t\leq T}$ denote the
two-dimensional diffusion process satisfying $(\ref{SDE})$ with a
parameter $\sigma$ and $Y^{(\sigma)}_0=Y_0$.
Let
$P_{\sigma,n}$ be the distribution\vspace*{1pt} of $(\Pi_n,Y^{(\sigma)}_{\Pi_n})$.

Our purpose is to obtain the LAMN property of probability measures
$\{P_{\sigma,n}\}_{\sigma\in\Lambda,n\in\mathbb{N}}$ of nonsynchronous
observations $(\Pi_n,Y^{(\sigma)}_{\Pi_n})$.
For this purpose, we will introduce several assumptions.
First, we consider conditions for the process $Y$.
\begin{enumerate}[{[$A1$]}]
\item[{[$A1$]}] \mbox{}
\begin{enumerate}[3.]
\item[1.] For $0\leq i+j\leq3$ and $0\leq k \leq4$, the derivatives
$\partial_t^i\partial_x^j\partial_{\sigma}^kb$ and $\partial_t^i\partial
_x^j\partial_{\sigma}^k\mu$
exist and are continuous with respect to $(t,x,\sigma)$.
Moreover, $\partial_x\mu,\partial_xb$ are bounded uniformly in
$[0,T]\times\mathbb{R}^2\times\Lambda$.
\item[2.] A matrix $(rb(t_1,x_1,\sigma)+(1-r)b(t_2,x_2,\sigma
))(rb(t_1,x_1,\sigma)+(1-r)b(t_2,x_2,\sigma))^{\star}$ is positive
definite for any $r\in[0,1]$, $t_1,t_2\in[0,T]$, $x_1,x_2\in\mathbb
{R}^2$ and $\sigma\in\Lambda$.
\item[3.] $E[|Y_0|^2]<\infty$.
\end{enumerate}
\end{enumerate}
Condition $[A1]$ is similar conditions to that for the LAMN property of
the statistical model with synchronous, equispaced observations in
Gobet \cite{gob}.
We do not need further conditions for the process $Y$.
If the diffusion coefficient $b$ is symmetric and positive definite, we
have $[A1]$ 2.

Second, we give assumptions of observation times.
Let $\{b_n\}_{n\in\mathbb{N}}$ be a sequence of positive numbers such
that $b_n\geq1$ and $b_n\to\infty$ as $n\to\infty$.

\begin{enumerate}[{[$A2$]}]
\item[{[$A2$]}] There\vspace*{1pt} exist positive constants $\{\delta_j\}_{j=1}^3$
such that $(5\delta_1+4\delta_3)\vee(3\delta_1+2\delta
_2+2\delta_3) \vee(3\delta_1/2+3\delta_2)<1/2$
and the following conditions hold true:
\begin{enumerate}[2.]
\item[1.] $r_n=\mathrm{O}_p(b_n^{-1+\delta_1})$.
\item[2.]
%
\begin{eqnarray}
\label{A2-eq1}
\hspace*{-36pt}\lim_{n\to\infty}b_n^2\sup
_{j_1,j_2\in\mathbb{N},|j_1-j_2|\geq
b_n^{\delta_2}}P \biggl[\ell_{1,n}\geq j_1\vee
j_2 \mbox{ and }  \frac
{|S^{n,j_2}-S^{n,j_1}|}{|j_2-j_1|}\leq b_n^{-1-\delta_3}
\biggr]=0,\qquad
\\
\label{A2-eq2}
\hspace*{-36pt}\lim_{n\to\infty}b_n^2\sup
_{j_1,j_2\in\mathbb{N},|j_1-j_2|\geq
b_n^{\delta_2}}P \biggl[\ell_{2,n}\geq j_1\vee
j_2 \mbox{ and }  \frac
{|T^{n,j_2}-T^{n,j_1}|}{|j_2-j_1|}\leq b_n^{-1-\delta_3}
\biggr]=0.\qquad
\end{eqnarray}
%
\end{enumerate}
\end{enumerate}
Condition $[A2]$ 2. controls the probability that too many observations
occur in some local interval.
For example, if we set $S^{n,i}=iT/n^2$ for $0\leq i\leq n$,
$S^{n,i}=(i+1-n)T/n$ for $n+1\leq i\leq2n-1$
and $T^{n,j}=jT/n$ for $0\leq j\leq n$, then we can easily see that
$[A2]$ 2. is not satisfied for $b_n\equiv n$.
In this setting, extremely many observations of $Y^1$ occur in the
interval $[0,T/n]$ compared to other intervals.
Condition $[A2]$ is a condition to exclude observations with such
extremely different frequency.
 This condition is necessary to obtain asymptotic equivalence
between the true log-likelihood ratios and the quasi-log-likelihood
ratios defined later
(Lemmas \ref{coeff-move-lemma1}, \ref{integration-lemma} and \ref{coeff-move-lemma2}),
and to obtain convergence results of the quasi-log-likelihood ratios
(Theorem~\ref{ogi-yos-thm}).

We need one more condition for observation times.
\begin{enumerate}[{[$A3$]}]
\item[{[$A3$]}] There exist $\sigma(\{\Pi_n\}_n)$-measurable\vspace*{1pt}
left-continuous processes $a_0(t)$ and $c_0(t)$ such that
$\int^T_0a_0(t)\,\mathrm{d}t\vee\int^T_0c_0(t)\,\mathrm{d}t<\infty$ almost\vspace*{1pt} surely,
$b_n^{-1}\operatorname{tr}(\mathcal{E}^1(t))\to^p \int^t_0a_0(s)\,\mathrm{d}s$ and
$b_n^{-1}\operatorname{tr}(\mathcal{E}^2(t))\to^p\int^t_0c_0(s)\,\mathrm{d}s$
as $n\to\infty$ for $t\in(0,T]$. Moreover, at least one of the
following conditions holds true:
\begin{enumerate}[2.]
\item[1.] There exist $\eta\in(0,1)$ and a $\sigma(\{\Pi_n\}
_n)$-measurable process $a(z,t)$ such that $a$ is continuous with
respect to $z$,
left-continuous with respect to $t$, $\int^T_0a(z,t)\,\mathrm{d}t<\infty$ a.s. and
\[
b_n^{-1}\operatorname{tr}\bigl(\mathcal{E}^1(t) \bigl(
\mathcal{E}_{\ell_{1,n}}-z^2GG^{\star
}\bigr)^{-1}
\bigr)\to^p \int^t_0a(z,s)\,
\mathrm{d}s
\]
as $n\to\infty$ for $t\in(0,T]$ and $z\in\mathbb{C},|z|<\eta$.
\item[2.] There exist $\eta\in(0,1)$ and a $\sigma(\{\Pi_n\}
_n)$-measurable process $c(z,t)$ such that $c$ is continuous with
respect to $z$,
left-continuous with respect to $t$, $\int^T_0c(z,t)\,\mathrm{d}t<\infty$ a.s. and
\[
b_n^{-1}\operatorname{tr}\bigl(\mathcal{E}^2(t) \bigl(
\mathcal{E}_{\ell_{2,n}}-z^2G^{\star
}G\bigr)^{-1}
\bigr)\to^p \int^t_0c(z,s)\,
\mathrm{d}s
\]
as $n\to\infty$ for $t\in(0,T]$ and $z\in\mathbb{C},|z|<\eta$.
\end{enumerate}
\end{enumerate}
In particular, $[A3]$ implies tightness of $\{b_n^{-1}(\ell_{1,n}+\ell
_{2,n})\}_n$.

Lemma~4 in Ogihara and Yoshida \cite{ogi-yos} shows that \textit{both} 1.
and 2. in $[A3]$ hold true if $r_n\to^p0$ and $[A3]$ holds true,
that is, the first statement of $[A3]$ and \textit{either} 1. or 2. in
$[A3]$ hold true.
Moreover, $a$ and $c$ are analytic with respect to $z$ and
$a(z,t)-a(0,t)=c(z,t)-c(0,t)$ for any $z\in\mathbb{C},|z|<\eta$ and
$t\in[0,T]$ almost surely,
assuming that $r_n\to^p0$ and $[A3]$ (Lemmas 3~and~4 and
Proposition~2 in \cite{ogi-yos}).
We will give tractable sufficient conditions of $[A2]$ and $[A3]$ in
Lemmas \ref{B1toA3} and \ref{B1toA2}.

The intuitive meaning of $[A3]$ is as follows. If $\mu\equiv0$ and
$b(t,x,\sigma)$ does not depend on $(t,x)$,
then $Y$ is a Wiener process and we obtain
\[
\log(\mathrm{d}P_{\sigma_{\ast}+b_n^{-1/2}u,n}/\mathrm{d}P_{\sigma_{\ast},n})=H_n\bigl(\sigma
_{\ast}+b_n^{-1/2}u\bigr)\circ(
\Pi,Y_{\Pi})-H_n(\sigma_{\ast})\circ(\Pi
,Y_{\Pi}),
\]
where $H_n(\sigma)$ is defined in (\ref{Hn-def}).
Roughly speaking, $H_n(\sigma)\circ(\Pi,Y_{\Pi})$ is asymptotically
equivalent~to
\begin{eqnarray*}
&& E\bigl[H_n(\sigma)\circ(\Pi,Y_{\Pi})|\Pi\bigr] \\
&&\quad  =-
\frac{|b^1|^2(\sigma_{\ast})}{2|b^1|^2(\sigma)}\operatorname{tr}\bigl(\bigl(\mathcal {E}_{\ell_{1,n}}-
\rho^2GG^{\star}\bigr)^{-1}\bigr) -\frac{|b^2|^2(\sigma_{\ast})}{2|b^2|^2(\sigma)}\operatorname{tr}\bigl(\bigl(\mathcal {E}_{\ell_{2,n}}-\rho^2G^{\star}G
\bigr)^{-1}\bigr)
\\
&&\qquad {}+\frac{b^1\cdot b^2(\sigma_{\ast})}{|b^1||b^2|(\sigma)}\operatorname{tr}\bigl(\rho \bigl(\mathcal{E}_{\ell_{1,n}}-
\rho^2GG^{\star}\bigr)^{-1}GG^{\star}\bigr) -
\frac{1}{2}\log\det S(\sigma),
\end{eqnarray*}
where $\rho=\rho(\sigma)=b^1\cdot b^2|b^1|^{-1}|b^2|^{-1}(\sigma)$.
Therefore, it is natural to assume conditions about asymptotic
behaviors of
$\operatorname{tr}((\mathcal{E}_{\ell_{1,n}}-\rho^2GG^{\star})^{-1})$ and $\operatorname{tr}((\mathcal{E}_{\ell_{2,n}}-\rho^2G^{\star}G)^{-1})$
in this special case of $\mu$ and $b$.
Since the diffusion coefficient of the diffusion process $Y$
in general is locally approximated by a constant
and asymptotic contribution of drift coefficient $\mu$ is negligible,
$[A3]$ is suitable for specifying asymptotic behaviors of
log-likelihood ratios in general cases.

Let $B^k_t=B^k_t(\sigma)=|b^k(t,Y_t,\sigma_{\ast})|/|b^k(t,Y_t,\sigma
)|$ for $k=1,2$, $\rho_t=\rho_t(\sigma)=b^1\cdot
b^2|b^1|^{-1}\* |b^2|^{-1} (t,Y_t,\sigma)$ and
%
\begin{eqnarray}
\Gamma&=& \int^T_0 \biggl\{
\partial_za\bigl(\rho_t(\sigma_{\ast}),t\bigr)
\frac
{(\partial_{\sigma}\rho_t(\sigma_{\ast}))^2}{\rho_t(\sigma_{\ast})}1_{\{
\rho_t(\sigma_{\ast})\neq0\} } \nonumber\\
&&\label{gamma-def}\hspace*{20pt}{}+2a\bigl(\rho_t(
\sigma_{\ast}),t\bigr) \bigl(\partial_{\sigma}B^1_t(
\sigma_{\ast
})\bigr)^2+2c\bigl(\rho_t(
\sigma_{\ast}),t\bigr) \bigl(\partial_{\sigma}B^2_t(
\sigma_{\ast
})\bigr)^2
\\
\nonumber
&&\hspace*{20pt}{}-\bigl(a\bigl(\rho_t(\sigma_{\ast}),t\bigr)-a(0,t)\bigr)
\biggl(\frac{\partial_{\sigma}\rho
_t(\sigma_{\ast})}{\rho_t(\sigma_{\ast})}1_{\{\rho_t(\sigma_{\ast})\neq
0\} }-\partial_{\sigma}B^1_t(
\sigma_{\ast})-\partial_{\sigma
}B^2_t(
\sigma_{\ast}) \biggr)^2 \biggr\}\,\mathrm{d}t.
\end{eqnarray}

We also assume the following condition.
\begin{enumerate}[{[$H$]}]
\item[{[$H$]}] The $d\times d$ random matrix $\Gamma$ is positive
definite almost surely.
\end{enumerate}

We can now formulate our main theorem.

%
\begin{theorem}\label{main}
Assume $[A1]$--$[A3]$ and $[H]$. Then the family $\{P_{\sigma,n}\}
_{\sigma,n}$ defined by nonsynchronous observations
$(\Pi_n,Y_{\Pi_n})$ has the LAMN property at $\sigma=\sigma_{\ast}$,
where $\mathcal{N}$ in Definition~\ref{LAMN-def} is a random variable
on an extension of $(\Omega,\mathcal{F},P)$, $\mathcal{N}$ is
independent of $\mathcal{F}$ and $\Gamma$ in Definition~\ref{LAMN-def}
is defined by (\ref{gamma-def}).
Moreover, $\mathcal{N}_n$ and $\Gamma_n$ can be taken so that
$(\mathcal{N}_n,\Gamma_n)\circ(\Pi_n,Y_{\Pi_n}) \to^{s\mbox{-}\mathcal{L}}(\mathcal{N},\Gamma)$.
\end{theorem}

Conditions $[A2]$, $[A3]$ and $[H]$ are  often not easy to
check for practical settings.
 So we see some easily tractable sufficient conditions for
these conditions.
\begin{enumerate}[{[$B1$]}]
\item[{[$B1$]}]
There\vspace*{1pt} exists exponential $\alpha$-mixing simple point process $\{\bar
{N}_t\}_{t\geq0}=\{(\bar{N}^1_t,\bar{N}^2_t)\}_{t\geq0}$ such that
$\bar{N}_0=0$,
$S^{n,i}=\inf\{t\geq0; \bar{N}^1_{b_nt}\geq i\}\wedge T$, $T^{n,j}=\inf
\{t\geq0; \bar{N}^2_{b_nt}\geq j\}\wedge T$
and the distribution of $(\bar{N}^i_{t+t_k}-\bar
{N}^i_{t+t_{k-1}})_{k=1}^M$ does not depend on $t\geq0$
for $M\in\mathbb{N}$, $0\leq t_0<t_1<\cdots< t_M$ and $i=1,2$.
Moreover,
$E[|\bar{N}_1|^q]<\infty$ and $\limsup_{u\rightarrow\infty}\max_{i=1,2}u^qP[\bar{N}^i_u=0]<\infty$
for any $q>0$.
\item[{[$H'$]}] There exists a constant $\varepsilon>0$ such that $|bb^{\star
}(t,x,\sigma_1)-bb^{\star}(t,x,\sigma_2)|\geq\varepsilon|\sigma_1-\sigma_2|$
for any $t\in[0,T]$, $x\in\mathbb{R}^2$ and $\sigma_1,\sigma_2\in
\Lambda$.
\end{enumerate}
For example, we can easily see that condition $[B1]$ is satisfied if
the processes $\{\bar{N}^1\}_{t\geq0}$ and $\{\bar{N}^2\}_{t\geq0}$
are two independent homogeneous Poisson processes.

The following lemma is proved in Section~6, Proposition~4 and Remark~2
in Ogihara and Yoshida \cite{ogi-yos}.
(We also use some localization techniques.)

%
\begin{lemma}\label{B1toA3}
\textup{1}. Condition $[B1]$ implies $[A3]$.
\begin{enumerate}[2.]
\item[2.] Assume $[A1]$, $[B1]$ and $[H']$. Then $[H]$ holds true.
\end{enumerate}
\end{lemma}

Let\vspace*{1pt} $\mathbf{N}^1_t=\sum_{i=1}^{\ell_{1,n}}1_{\{S^{n,i}\leq t\} }$ and
$\mathbf{N}^2_t=\sum_{j=1}^{\ell_{2,n}}1_{\{T^{n,j}\leq t\} }$.
Then we also have the following. The proof is left in  the \hyperref[app]{Appendix}.

%
\begin{lemma}\label{B1toA2}
Let $q>0$. Assume that there exists $n_0\in\mathbb{N}$ such that
\[
\sup_{n\geq n_0}\max_{1\leq i\leq2}\sup
_{0\leq t\leq
T-b_n^{-1}}E\bigl[\bigl(\mathbf{N}^i_{t+b_n^{-1}}-
\mathbf{N}^i_t\bigr)^q\bigr]<\infty.
\]
Then (\ref{A2-eq1}) and (\ref{A2-eq2}) hold true for any $\delta_2>3/q$
and $\delta_3>3/q$. In particular, $[B1]$ implies $[A2]$.
\end{lemma}

\begin{remark}
Conditions $[B1]$ and $[H']$ are the simplest sufficient conditions of
$[A2]$, $[A3]$ and~$[H]$. More detailed discussion about sufficient
conditions of $[A3]$ and $[H]$ can be found
in Sections~4 and~6 in Ogihara and Yoshida \cite{ogi-yos} and Section~4
in Uchida and Yoshida \cite{uch-yos02}.
\end{remark}

By Theorem~\ref{main}, we obtain the minimax theorem (\ref
{minimax-ineq}) under the conditions in Theorem~\ref{main}.
In the rest of this section, we will prove that the quasi-maximum
likelihood estimator and the Bayes-type estimator defined in Ogihara
and Yoshida \cite{ogi-yos}
attain the lower bound in (\ref{minimax-ineq}) under certain conditions.
So these estimators are asymptotically efficient in this sense.
For these purposes, we use the scheme of
Yoshida \cite{yos05} which leads to convergence of moments of
estimators.

We will make the assumptions for asymptotic efficiency of estimators.
We denote $\omega_{\alpha}(g)=\sup_{t\neq s}|g(t)-g(s)|/|t-s|^{\alpha}$
for $\alpha\in(0,1/2)$ and an $\alpha$-H\"{o}lder continuous
function $g\dvtx [0,T]\to\mathbb{R}$.
Let $\mathbf{K}(\bar{u})=\{[s^{i-1},s^i)\}_{i=1}^{L^1}\cup\{
[t^{j-1},t^j)\}_{j=1}^{L^2}$ and
$\{\theta(p,l;\bar{u})\}_{1\leq l\leq L^1+L^2,p\in\mathbb{Z}_+}$ be
defined by
$\theta(0,l;\bar{u})=[s^{l-1},s^l)$  $(1\leq l \leq L^1)$, $\theta
(0,l;\bar{u})=[t^{l-L^1-1},t^{l-L^1})$  $(L^1<l\leq L^1+L^2)$ and
\begin{eqnarray*}
\theta(p,l;\bar{u}) &=& \bigcup\bigl\{K_{2p};K_1,\ldots,
K_{2p}\in\mathbf{K}(\bar {u}),\\
&&\hspace*{17pt} K_1\cap\theta(0,l;\bar{u})
\neq\varnothing,K_j\cap K_{j-1}\neq \varnothing \ (2\leq j
\leq2p)\bigr\}
\end{eqnarray*}
for $p\in\mathbb{N}$, $\bar{u}=((s^i)_{i=0}^{L^1},(t^j)_{j=0}^{L^2})\in
\mathcal{U}$ and $1\leq l\leq L^1+L^2$.
That is, the interval $\theta(p,l;\bar{u})$ is the union of intervals
which are reached by $2p$ \textit{transfers} from $\theta(0,l;\bar{u})$.
Let $\theta_{p,l}=\theta(p,l;\Pi)$.

Let $q>2$, $\delta\in(0,1)$, $\delta'\geq1$ and $\eta\in(0,1)$.
\begin{enumerate}[{[$C4\mbox{-}q,\delta'$]}]
\item[{[$C1$]}] \mbox{}
\begin{enumerate}[5.]
\item[1.] The functions $b$ and $\mu$ have continuous derivatives $\partial
_t^i\partial_x^j\partial_{\sigma}^kb,\partial_t^{i'}\partial
_x^{j'}\partial_{\sigma}^{k'}\mu$
and satisfy
\begin{eqnarray*}
 \sup_{t\in[0,T],\sigma\in\Lambda}\bigl|\partial_t^i
\partial_x^j\partial _{\sigma}^kb(t,x,
\sigma)\bigr| &\leq  & C\bigl(1+|x|\bigr)^C \quad\mbox{and}\\
 \sup_{t\in[0,T],\sigma\in\Lambda}\bigl|
\partial _t^{i'}\partial_x^{j'}
\partial_{\sigma}^{k'}\mu(t,x,\sigma)\bigr| &\leq &  C\bigl(1+|x|\bigr)^C
\end{eqnarray*}
for $0\leq i+j\leq3, 0\leq k\leq4$, $0\leq i'+j'+k'\leq1$ and $x\in
\mathbb{R}^2$.
\item[2.] The derivatives $\partial_x\mu$ and $\partial_xb$ are bounded
uniformly in $[0,T]\times\mathbb{R}^2\times\Lambda$.

\item[3.] $\inf_{t,x,\sigma} \det bb^{\star}(t,x,\sigma)>0$.
\item[4.] $\sup_{\sigma}\sup_{0\leq t\leq T}E[|Y^{(\sigma)}_t|^q]<\infty$
for any $q>0$.
\item[5.] The function $\partial_{\sigma}^kb$ can be continuously extended
to $[0,T]\times\mathbb{R}^2\times\operatorname{clos}(\Lambda)$ for $0\leq
k\leq4$, where $\operatorname{clos}(\Lambda)$ represents the closure of $\Lambda$.
\end{enumerate}

\item[{[$C2\mbox{-}q,\delta$]}] $E[r_n^q]=\mathrm{O}(b_n^{-\delta q})$.
\item[{[$C3\mbox{-}q,\eta$]}] There exist $n_0\in\mathbb{N}$, $\alpha
\in(0,1/2-1/q)$
and $\sigma(\{\Pi_n\}_n)$-measurable\vspace*{1pt} left-continuous processes $\{
a_p(t)\}_{p\in\mathbb{Z}_+}$ and $\{c_p(t)\}_{p\in\mathbb{Z}_+}$
such that $\int^T_0(a_p\vee c_p)(t)\,\mathrm{d}t \in L^q(\Omega)$ for $p\in\mathbb
{Z}_+$, $E[(\ell_{1,n}+\ell_{2,n})^q]<\infty$ for $n\in\mathbb{N}$ and
\begin{eqnarray*}
&& E \Biggl[ \Biggl(b_n^{\eta} \Biggl|b_n^{-1}
\sum_{i=1}^{\ell
_{1,n}}g\bigl(S^{n,i-1}\bigr)
\bigl(\bigl(GG^{\star}\bigr)^p\bigr)_{ii}-\int
^T_0g(t)a_p(t)\,\mathrm{d}t \Biggr|
\Biggr)^q \Biggr]
\\
&&\qquad{}\vee E \Biggl[ \Biggl(b_n^{\eta} \Biggl|b_n^{-1}
\sum_{j=1}^{\ell
_{2,n}}g\bigl(T^{n,j-1}\bigr)
\bigl(\bigl(G^{\star}G\bigr)^p\bigr)_{jj}-\int
^T_0g(t)c_p(t)\,\mathrm{d}t \Biggr|
\Biggr)^q \Biggr]
\\
&& \quad\leq C(p+1)^C \Bigl(\sup_t\bigl|g(t)\bigr|^q+
\omega_{\alpha}(g)^q \Bigr)
\end{eqnarray*}
for $n\geq n_0$, $p\in\mathbb{Z}_+$ and any $\alpha$-H\"{o}lder continuous function $g$ on $[0,T]$.
\item[{[$C4\mbox{-}q,\delta'$]}]
\begin{eqnarray*}
&& \lim_{n\to\infty} \Biggl\{E \Biggl[\bigl(b_n^{-{q}/{2}}
\vee r_n^q\bigr)\sum_{p=0}^{\infty}
\frac{(\sum_{l=1}^{\ell_{1,n}+\ell_{2,n}}|\theta
_{2p+2,l}|)^q}{(p+1)^{q\delta'}} \Biggr]\\
&&\hspace*{5pt} \qquad {}\vee E \Biggl[ \Biggl(b_n^{-1}
\sum_{p_1,p_2=0}^{\infty}\frac{\sum_{l_1,l_2=1}^{\ell_{1,n}+\ell
_{2,n}}|\theta_{2p_1+3,l_1}\cap\theta_{2p_2+3,l_2}|}{(p_1+1)^{\delta
'}(p_2+1)^{\delta'}}
\Biggr)^{{q}/{2}} \Biggr] \Biggr\}\\
&& \quad =0.
\end{eqnarray*}
\end{enumerate}

Condition $[C3\mbox{-}q,\eta]$ is a stronger condition than
$[A3]$ and is required to obtain
moment convergence of estimation errors.
For any $q>2$ and $\eta\in(0,1)$, we can prove that $[C3\mbox{-}q,\eta]$ implies~$[A3]$. See Section~3.1 in Ogihara and Yoshida \cite{ogi-yos}  for the details.

Condition $[C4\mbox{-}q,\delta']$ is a technical condition to obtain
the asymptotic properties of  quasi-likelihood ratios and its
derivatives.
This condition together with Lemma~13 in \cite{ogi-yos} enable us to
apply martingale limit theorems to the quasi-likelihood
ratios, and hence
it is essential to obtain asymptotic properties of
quasi-likelihood ratios.
See Propositions 3 and 10 in \cite{ogi-yos} and their proofs for the details.

Let $B^k_t(\sigma_1;\sigma_2)=|b^k(t,Y^{(\sigma_2)}_t,\sigma
_2)|/|b^k(t,Y^{(\sigma_2)}_t,\sigma_1)|$ for $k=1,2$,
$\rho_t(\sigma_1;\sigma_2)=b^1\cdot b^2|b^1|^{-1}\* |b^2|^{-1}(t,Y^{(\sigma
_2)}_t,\sigma_1)$, and
\begin{eqnarray*}
\mathcal{Y}(\sigma_1;\sigma_2)&=&\int
^T_0 \biggl\{-\frac{(B^1_t(\sigma
_1;\sigma_2))^2}{2}a\bigl(
\rho_t(\sigma_1;\sigma_2),t\bigr)-
\frac{(B^2_t(\sigma
_1;\sigma_2))^2}{2}c\bigl(\rho_t(\sigma_1;
\sigma_2),t\bigr)
\\
&&\hspace*{20pt}{}+B^1_tB^2_t(
\sigma_1;\sigma_2) \bigl(a\bigl(\rho_t(
\sigma_1;\sigma _2),t\bigr)-a_0(t)\bigr)
\frac{\rho_{t}(\sigma_2;\sigma_2)}{\rho_t(\sigma_1;\sigma
_2)}1_{\{\rho_t(\sigma_1;\sigma_2)\neq0\} }\\
&&\hspace*{20pt}{}+\frac{a_0(t)}{2}+\frac{c_0(t)}{2}
+a_0(t)\log B^1_t(\sigma_1;
\sigma_2)
\\
&&\hspace*{20pt}{}+c_0(t)\log B^2_t(\sigma_1;
\sigma_2)+\int^{\rho_t(\sigma_1;\sigma
_2)}_{\rho_{t}(\sigma_2;\sigma_2)}
\frac{a(\rho,t)-a_0(t)}{\rho}\,\mathrm{d}\rho \biggr\}\,\mathrm{d}t,
\nonumber
\end{eqnarray*}
where $\{a(z,t)\}$ and $\{c(z,t)\}$ are in $[A3]$.

\begin{enumerate}[{[$C5$]}]
\item[{[$C5$]}] There exist a family $\{\tilde{c}_q\}_{q>0}$ of
positive constants and an open set $\Lambda'$ satisfying $\sigma_{\ast
}\in\Lambda'\subset\Lambda$ such that
$\sup_{\sigma_2\in\Lambda'} P[\inf_{\sigma_1\in\Lambda\setminus\{
\sigma_2\}}(-\mathcal{Y}(\sigma_1;\sigma_2)/|\sigma_1-\sigma_2|^2)\leq
r^{-1}]\leq \tilde{c}_q/r^q$ for $r>0$ and $q>0$.
\end{enumerate}
We see that $[C5]$ implies $[H]$, by using the relations $\mathcal
{Y}(\sigma_{\ast};\sigma_{\ast})=\partial_{\sigma_1}\mathcal{Y}(\sigma
_1;\sigma_{\ast})|_{\sigma_1=\sigma_{\ast}}=0$,
$\Gamma=-\partial^2_{\sigma_1}\mathcal{Y}(\sigma_1;\sigma_{\ast
})|_{\sigma_1=\sigma_{\ast}}$,
and hence $\inf_{\sigma\neq\sigma_{\ast}}(-\mathcal{Y}(\sigma;\sigma
_{\ast})/|\sigma-\sigma_{\ast}|^2)\leq\inf_{u\neq0}u^{\star}\Gamma
u/(2|u|^2)$.
Condition $[C5]$ and $[H]$ are conditions about
identifiability of statistical models.
We only need $[H]$ to have Theorem~\ref{main}. However, we need $[C5]$
to obtain asymptotic efficiency of estimators.

Ogihara and Yoshida \cite{ogi-yos} proposed a quasi-log-likelihood
function $H_n$ defined by
%
\begin{equation}
\label{Hn-def}
H_n(\sigma) \circ(\Pi,Y_{\Pi})=-
\tfrac{1}{2}Z^{\star}S^{-1}(\sigma )Z-\tfrac{1}{2}\log
\det S(\sigma),
\end{equation}
where
%
\begin{equation}
\label{Z-def}
Z=\bigl(\bigl(\bigl(Y^1_{S^{n,i}}-Y^1_{S^{n,i-1}}
\bigr)/\sqrt{\bigl|I^i\bigr|}\bigr)_i^{\star
},\bigl(
\bigl(Y^2_{T^{n,j}}-Y^2_{T^{n,j-1}}\bigr)/
\sqrt{\bigl|J^j\bigr|}\bigr)_j^{\star}\bigr)^{\star},
\end{equation}
$\bar{b}^1_{(i)}=\bar{b}^1_{(i)}(\sigma)=b^1(S^{n,i-1},
Y^1_{S^{n,i-1}}, Y^2_{T^{n,j'}},\sigma)$ for $j'=\max\{j;T^{n,j}\leq
S^{n,i-1}\}$,
$\bar{b}^2_{(j)}=\bar{b}^2_{(j)}(\sigma)=b^2(T^{n,j-1}, Y^1_{S^{n,i'}},
Y^2_{T^{n,j-1}},\sigma)$ for $i'=\max\{i;S^{n,i}\leq T^{n,j-1}\}$ and
%
\begin{equation}
\label{S-def}
S(\sigma)=\pmatrix{ \operatorname{diag}
\bigl(\bigl\{\bigl|\bar{b}^1_{(i)}\bigr|^2\bigr
\}_i\bigr) & \bigl\{\bar{b}^1_{(i)}\cdot\bar
{b}^2_{(j)}G_{ij}\bigr\}_{i,j}
\cr\noalign{\vspace*{3pt}}
\bigl\{\bar{b}^1_{(i)}\cdot\bar{b}^2_{(j)}G_{ij}
\bigr\}_{j,i} & \operatorname{diag}\bigl(\bigl\{ \bigl|\bar{b}^2_{(j)}\bigr|^2
\bigr\}_j\bigr)}.
\end{equation}

An intuitive meaning of $H_n$ is as follows. If $\mu\equiv0$,
$b(t,x,\sigma)$ does not depend on $x$ and $\Pi$ is deterministic,
then $Z$ follows a zero-mean normal distribution.
Moreover, the covariance matrix of $Z$ is approximated as
\begin{eqnarray*}
E \biggl[\frac{Y^1_{S^{n,i}}-Y^1_{S^{n,i-1}}}{\sqrt{|I^i|}}\frac
{Y^1_{S^{n,i'}}-Y^1_{S^{n,i'-1}}}{\sqrt{|I^{i'}|}} \biggr] 
&\sim & \bigl|
\bar{b}^1_{(i)}\bigr|^2(\sigma_{\ast})
\delta_{i,i'},\\
E \biggl[\frac{Y^2_{T^{n,j}}-Y^2_{T^{n,j-1}}}{\sqrt{|J^j|}}\frac
{Y^2_{T^{n,j'}}-Y^2_{T^{n,j'-1}}}{\sqrt{|J^{j'}|}} \biggr]
&\sim& \bigl|\bar{b}^2_{(j)}\bigr|^2(
\sigma_{\ast})\delta_{j,j'},
\\
E \biggl[\frac{Y^1_{S^{n,i}}-Y^1_{S^{n,i-1}}}{\sqrt{|I^i|}}\frac
{Y^2_{T^{n,j}}-Y^2_{T^{n,j-1}}}{\sqrt{|J^j|}} \biggr] &\sim & \bar{b}^1_{(i)}(
\sigma_{\ast})\cdot\bar{b}^2_{(j)}(
\sigma_{\ast})G_{ij}.
\end{eqnarray*}
Hence, $S(\sigma)$ is approximation of the covariance matrix of $Z$.
Therefore, we can say $H_n(\sigma)$ is an approximate log-likelihood function.
These arguments are valid only for this special case of $\mu$, $b$ and
$\Pi$.
However, Ogihara and Yoshida \cite{ogi-yos} define $H_n$ as above for
general cases of $\mu$, $b$ and $\Pi$
and studied the quasi-maximum likelihood estimator and the Bayes-type
estimator constructed by $H_n$.

Let $\pi\dvtx \Lambda\to(0,\infty)$ be a bounded continuous function.
The quasi-maximum likelihood estimator $\hat{\sigma}_n$ and the
Bayes-type estimator $\tilde{\sigma}_n$ for the prior density $\pi$ are
defined by
$\hat{\sigma}_n=\operatorname{argmax}_{\sigma\in\operatorname{clos}(\Lambda)}H_n(\sigma)$ and
\[
\tilde{\sigma}_n= \biggl(\int_{\Lambda}\exp
\bigl(H_n(\sigma)\bigr)\pi(\sigma)\,\mathrm{d}\sigma \biggr)^{-1}\int
_{\Lambda}\sigma\exp\bigl(H_n(\sigma)\bigr)\pi(
\sigma)\,\mathrm{d}\sigma.
\]

Let $\sigma^n_u=\sigma_{\ast}+b_n^{-1/2}u$ for $u\in\mathbb{R}^d$.

%
\begin{theorem}\label{ogi-yos-eff}
Let $\delta\in(0,1/2)$. Assume that $0<\inf_{\sigma}\pi(\sigma)\leq
\sup_{\sigma}\pi(\sigma)<\infty$ and that for any $q>0$, there exist
$\delta'\geq1$ and $q'\in\mathbb{N}$ satisfying $2q'>q$ such that
$[C1]$, $[C2\mbox{-}(2q'),\delta]$, $[C3\mbox{-}(2q'),\delta]$,
$[C4\mbox{-}(2q'),\delta']$, $[C5]$ hold.
Then
\begin{eqnarray}
\lim_{\alpha\to\infty}\liminf_{n\to\infty}\sup
_{|u|\leq\alpha}E_{\sigma
^n_u}\bigl[l\bigl(\bigl|b_n^{1/2}
\bigl(\hat{\sigma}_n-\sigma^n_u\bigr)\bigr|
\bigr)\bigr] &=& E\bigl[l\bigl(\bigl|\Gamma ^{-1/2}\mathcal{N}\bigr|\bigr)\bigr],
\nonumber
\\
\lim_{\alpha\to\infty}\liminf_{n\to\infty}\sup
_{|u|\leq\alpha}E_{\sigma
^n_u}\bigl[l\bigl(\bigl|b_n^{1/2}
\bigl(\tilde{\sigma}_n-\sigma^n_u\bigr)\bigr|
\bigr)\bigr] &=& E\bigl[l\bigl(\bigl|\Gamma ^{-1/2}\mathcal{N}\bigr|\bigr)\bigr]
\nonumber
\end{eqnarray}
for any continuous function $l\dvtx [0,\infty)\to[0,\infty)$ that is
nondecreasing, $l(0)=0$ and of at most polynomial growth.
\end{theorem}

%
\begin{remark}
Theorems \ref{ogi-yos-eff}~and~\ref{main} and the minimax theorem by
Jeganathan \cite{jeg} imply
that estimators $\hat{\sigma}_n$ and $\tilde{\sigma}_n$ are
asymptotically efficient under $[A1]$, $[A2]$ and the conditions in
Theorem~\ref{ogi-yos-eff}.
\end{remark}

\begin{pf*}{Outline of the proof of Theorem~\ref{ogi-yos-eff}}
Let $\mathbf{G}_u=b_n^{1/2}(\hat{\sigma}_n\circ(\Pi, Y_{\Pi}^{(\sigma
^n_u)}) - \sigma^n_u)$. Then Theorem~2 in Ogihara and Yoshida \cite{ogi-yos} yields
%
\begin{equation}
\label{riskFunction-conv}
\lim_{n\to\infty}E\bigl[l\bigl(|\mathbf{G}_0|\bigr)
\bigr]=E\bigl[l\bigl(\bigl|\Gamma^{-1/2}\mathcal{N}\bigr|\bigr)\bigr].
\end{equation}
Moreover, for any $\varepsilon,\delta>0$, there exists $n_1\in\mathbb{N}$
such that $\sup_{|u|\leq\alpha}P[|\mathbf{G}_u-\mathbf{G}_0|>\delta]<\varepsilon
$ for $n\geq n_1$,
by a similar argument to the proof of 1. of Theorem~2 in
Ogihara and Yoshida \cite{ogi-yos} and relations $E[\sup_t|Y^{(\sigma
^n_u)}_t-Y^{(\sigma_{\ast})}_t|^q]\leq C_qb_n^{-q/2}|u|^q$ for any
$q\geq2$.

Furthermore, we obtain $\sup_{|u|\leq\alpha}E[|\mathbf{G}_u|^q]<\infty$
for any $\alpha>0$, $q>0$ and sufficiently large $n$,
by a similar argument to the proof of Proposition~5 in Ogihara and
Yoshida \cite{ogi-yos}.

Then for any $\varepsilon>0$, there exist $M'$, $n'$ and $\delta$ such that
\begin{eqnarray*}
\sup_{|u|\leq\alpha}\bigl|E\bigl[l\bigl(|\mathbf{G}_u|\bigr)\bigr]-E
\bigl[l\bigl(|\mathbf{G}_0|\bigr)\bigr]\bigr| &\leq& \sup_{|u|\leq\alpha}\bigl|E
\bigl[l\bigl(|\mathbf{G}_u|\bigr)-l\bigl(|\mathbf{G}_0|\bigr),
|\mathbf{G}_u|\vee|\mathbf{G}_0|\leq M'\bigr]\bigr| + \varepsilon
\\
&\leq& \sup_{|x|\leq M'} l(x) \sup_{|u|\leq\alpha}P\bigl[|
\mathbf{G}_u-\mathbf{G}_0|\geq\delta\bigr]+2\varepsilon<3
\varepsilon
\end{eqnarray*}
for $n\geq n'$, by continuity of $l$.

Hence, we obtain
\[
\lim_{\alpha\to\infty}\liminf_{n\to\infty}\sup
_{|u|\leq\alpha
}E\bigl[l\bigl(|\mathbf{G}_u|\bigr)\bigr]=\lim
_{n\to\infty}E\bigl[l\bigl(|\mathbf{G}_0|\bigr)\bigr]=E\bigl[l
\bigl(\bigl|\Gamma ^{-1/2}\mathcal{N}\bigr|\bigr)\bigr]
\]
by (\ref{riskFunction-conv}). We can similarly obtain the result for
the Bayes-type estimator $\tilde{\sigma}_n$.
\end{pf*}

The following corollary is obtained by the argument in Section~6 in
Ogihara and Yoshida \cite{ogi-yos}.

%
\begin{corollary}
Assume that $0<\inf_{\sigma}\pi(\sigma)\leq\sup_{\sigma}\pi
(\sigma)<\infty$
and that $[C1],[B1]$ and $[H']$ hold. Then the results in Theorem~\ref{ogi-yos-eff} hold true.
\end{corollary}

\section{Preliminary results}\label{preliminaries}

In the rest of this paper, we will prove Theorem~\ref{main}.
For this purpose, we will prove asymptotic equivalence
between the log-likelihood ratio $\log(\mathrm{d}P_{\sigma^n_u}/\mathrm{d}P_{\sigma_{\ast
}})(Y_{\Pi})$ of the processes $Y^{(\sigma)}$
and the quasi-log-likelihood ratio $H_n(\sigma^n_u)-H_n(\sigma_{\ast})$.
Then we obtain Theorem~\ref{main} since $H_n(\sigma^n_u)-H_n(\sigma_{\ast})$ has a LAMN-type property.

This section is devoted to some auxiliary results.
We use Malliavin calculus techniques and prove estimates for transition
density functions and their derivatives in Section~\ref{malliavin-section}.
Section~\ref{tightness-section} is devoted to prove some tightness
results of log-likelihood ratios.
These results play essential roles in the proof of Theorem~\ref{main}
in Section~\ref{LAMN-proof-section}.

\subsection{Some fundamental results}\label{fundamental-section}

In this subsection, we define Euler--Maruyama-type processes and
related notation.
We also introduce a LAMN-type property of $H_n(\sigma^n_u)-H_n(\sigma_{\ast})$.

First, we prepare several fundamental lemmas. The first one is about
localization.
To obtain Theorem~\ref{main}, it is sufficient to consider the
following stronger condition $[A1']$ instead of $[A1]$.
\begin{enumerate}[{[$A1'$]}]
\item[{[$A1'$]}] Condition $[A1]$ holds true, $|Y_0|\leq M$ a.s. for some
$M>0$, and $b$, $\mu$ and their derivatives
are bounded on $[0,T]\times\mathbb{R}^2 \times\Lambda$.
Moreover, there exist positive constants $\eta_{\min}$ and $\eta_{\max
}$ such that
\begin{eqnarray*}
\eta_{\min}\mathcal{E}_2 &\leq & \bigl(rb(t_1,x_1,
\sigma )+(1-r)b(t_2,x_2,\sigma)\bigr)
\bigl(rb(t_1,x_1,\sigma)+(1-r)b(t_2,x_2,
\sigma )\bigr)^{\star}\\
&\leq & \eta_{\max}\mathcal{E}_2
\end{eqnarray*}
for any $r\in[0,1]$, $t_1,t_2\in[0,T]$, $x_1,x_2\in\mathbb{R}^2$ and
$\sigma\in\Lambda$.
\item[{[$L$]}]  There exists a $d$-dimensional standard normal
random variable $\mathcal{N}$ on an extension of $(\Omega,\mathcal{F},P)$
such that $\mathcal{N}$ is independent of $\mathcal{F}$,
$-b_n^{-1}\partial_{\sigma}^2H_n(\sigma_{\ast})\circ(\Pi,Y_{\Pi})\to^p
\Gamma$,
\[
\log\frac{\mathrm{d}P_{\sigma_{\ast}+b_n^{-1/2}u,n}}{\mathrm{d}P_{\sigma_{\ast},n}}- \biggl(u^{\star}b_n^{-1/2}
\partial_{\sigma}H_n(\sigma_{\ast})+
\frac
{1}{2}u^{\star}b_n^{-1}
\partial_{\sigma}^2H_n(\sigma_{\ast})u
\biggr) \to0
\]
in $P_{\sigma_{\ast},n}$ probability, and
\[
b_n^{-1/2}\partial_{\sigma}H_n(
\sigma_{\ast}) \circ(\Pi,Y_{\Pi}) \to ^{s\mbox{-}\mathcal{L}}
\Gamma^{1/2}\mathcal{N}
\]
for $\Gamma$ defined in (\ref{gamma-def}).
\end{enumerate}

Let $\mathcal{H}=\{\omega\in\Omega; -\partial_{\sigma
}^2H_n(\sigma_{\ast})(\omega)  \mbox{ is   positive   definite}\}$,
%
\begin{equation}
\label{Gamman-def}
\Gamma_n=-b_n^{-1}
\partial_{\sigma}^2H_n(\sigma_{\ast})1_{\mathcal
{H}}+
\mathcal{E}_d1_{\mathcal{H}^c}, \qquad \mathcal{N}_n=\bigl(-
\partial_{\sigma}^2H_n(\sigma_{\ast
})
\bigr)^{-1/2}\partial_{\sigma}H_n(
\sigma_{\ast})1_{\mathcal{H} }.
\end{equation}
%

%
\begin{lemma}\label{A1toA1}
Assume that $[L]$ holds true under $[A1']$, $[A2]$ and $[A3]$.
Then Theorem~\ref{main} holds true with $\Gamma_n$ and $\mathcal{N}_n$
in (\ref{Gamman-def}).
\end{lemma}

\begin{pf}
Similar to the proof of Lemma~4.1. in Gobet \cite{gob} and we omit the details.
\end{pf}

The second lemma is Lemma~11 in Ogihara and Yoshida \cite{ogi-yos}.

%
\begin{lemma}\label{cond-conv}
Let $\{\mathsf{X}_n\}_{n\in\mathbb{N}}$ be a sequence of integrable random
variables on some probability space $(\Omega',\mathcal{F}',P')$ and
$\{\mathcal{G}_n\}_{n\in\mathbb{N}}$ be sub $\sigma$-fields of $\mathcal{F}'$.
Assume $E'[\mathsf{X}_n|\mathcal{G}_n]\to^p 0$ as $n\to\infty$. Then $\mathsf{X}_n\to^p0$ as $n\to\infty$.
\end{lemma}

Moreover, the following lemma is proved similarly to Lemma~\ref{cond-conv}.

%
\begin{lemma}\label{cond-conv-to-conv}
Let $\Theta$ be a set,
$\{\mathsf{X}_{n,\lambda}\}_{n\in\mathbb{N},\lambda\in\Theta}$ be a family
of integrable random variables on some probability space $(\Omega
',\mathcal{F}',P')$ and
$\{\mathcal{G}_n\}_{n\in\mathbb{N}}$ be sub $\sigma$-fields of $\mathcal{F}'$.
Assume that for any $\varepsilon>0$, there exists $M>0$ such that
$\sup_{n,\lambda}P'[E'[|\mathsf{X}_{n,\lambda}||\mathcal{G}_n]>M]<\varepsilon$.
Then for any $\varepsilon>0$, there exists $M>0$ such that
$\sup_{n,\lambda}P'[|\mathsf{X}_{n,\lambda}|>M]<\varepsilon$.
\end{lemma}

Let $0\leq s<t\leq T$, $r\in[0,1]$ and $\sigma\in\Lambda$. Under
$[A1]$, a stochastic differential equation
\[
\left\{
\begin{array}{l} \mathcal{X}^{r,\sigma}_v=
\bigl\{(1-r)\mu\bigl(v+s,\mathcal{X}^{r,\sigma
}_v,\sigma\bigr)+r
\mu(s,z_0,\sigma)\bigr\}\,\mathrm{d}v
\\[3pt]
\phantom{\mathcal{X}^{r,\sigma}_v=}{}+\bigl\{(1-r)b\bigl(v+s,\mathcal{X}^{r,\sigma}_v,\sigma\bigr)
+rb(s,z_0,\sigma)\bigr\} \,\mathrm{d}W_{v+s},\qquad  v\in[0,t-s],
\\[3pt]
\mathcal{X}^{r,\sigma}_0 = z_0
\end{array}
\right.
\]
has a unique strong solution $\{\mathcal{X}^{r,\sigma}_v\}_{0\leq v\leq
t-s}$. Let $p(z_1;z_0,r,s,t,\sigma)$ be
the probability density function of $\mathcal{X}^{r,\sigma}_{t-s}$.

The following lemma is classical estimate. See Theorem~1 in Aronson
\cite{aro} or Proposition~5.1 in Gobet \cite{gob}.

\begin{lemma}\label{density-est}
Assume $[A1']$. Then there exist positive constants $\mu_1<\mu_2$ and
$C>1$ such that
\[
\frac{1}{C}\frac{\mu_2}{2\uppi(t-s)}\exp \biggl(-\frac{\mu
_2|z_1-z_0|^2}{2(t-s)} \biggr)\leq
p(z_1;z_0,r,s,t,\sigma)\leq C\frac{\mu
_1}{2\uppi(t-s)}\exp
\biggl(-\frac{\mu_1|z_1-z_0|^2}{2(t-s)} \biggr),
\]
for $0\leq s< t\leq T$, $r\in[0,1]$, $z_0,z_1\in\mathbb{R}^2$ and
$\sigma\in\Lambda$.
\end{lemma}

We will define some further notation.
Let $n_u$ be the minimum positive integer satisfying $\{\sigma^n_{vu}\}
_{n\geq n_u,0\leq v\leq1}\subset\Lambda$ for $u\in\mathbb{R}^d$. For
$\bar{u}=((s^i)_i,(t^j)_j)\in\mathcal{U}$,
let $\check{u}=\{\check{u}^{k}(\bar{u})\}_{k=0}^{L_0(\bar{u})}$ be a
strictly increasing sequence of the elements of $\bar{u}$ such that
$\check{u}$  is equal to $\bar{u}$ as a set.
Let $\Delta\check{u}^k=\check{u}^k-\check{u}^{k-1}$, $k_1(i)=k_1(i;\bar
{u})$ be $k$ satisfying $s^i=\check{u}^k$ and $k_2(j)=k_2(j;\bar{u})$
be $k$ satisfying $t^j=\check{u}^k$,
\begin{eqnarray*}
i(k) &=& i(k;\bar{u}) = \max\bigl\{i; \mbox{ there  exists } j  \mbox{ such  that }
s^i\leq t^j\leq\check{u}^{k-1}\bigr\},
\\
j(k) &=& j(k;\bar{u}) = \max\bigl\{j; \mbox{ there  exists }  i  \mbox{ such  that }
t^j\leq s^i\leq\check{u}^{k-1}\bigr\}.
\end{eqnarray*}
We define random times $\check{U}^k=\check{u}^k(\Pi)$ and $\check{U}=\{
\check{U}^k\}_k$.

For $\bar{u}=((s^i)_{i=0}^{L^1},(t^j)_{j=0}^{L^2})\in\mathcal{U}$ and
$z=((x_k)_{k=0}^{L_0(\bar{u})},(y_k)_{k=0}^{L_0(\bar{u})})\in\mathbb
{R}^{2L_0(\bar{u})+2}$,
we denote $\bar{z}=((x_{k_1(i)})_{i=1}^{L^1},(y_{k_2(j)})_{j=1}^{L^2})$
and $\hat{z}=((x_k)_{k\notin\{k_1(i);0\leq i\leq L^1\}},(y_k)_{k\notin\{k_2(j);0\leq j\leq L^2\}})$.

Now, let us define stochastic processes that \textit{connect} the process
$Y^{(\sigma)}$ and an Euler--Maruyama process.
Let $u\in\mathbb{R}^d$, $\bar{u}\in\mathcal{U}$, $r\in[0,1]$
and $n\geq n_u$. Under $[A1']$, there exists a unique two-dimensional
stochastic process
$Y^{r,u}=\{Y^{r,u}_t\}_{0\leq t\leq T}=\{(Y^{r,u,1}_t(\bar
{u}),\break Y^{r,u,2}_t(\bar{u}))\}_{0\leq t\leq T}$ satisfying
\begin{eqnarray*}
Y^{r,u}_t&=&Y_0+\sum
_{k=1}^{L_0}\int^{t\wedge\check{u}^{k}}_{t\wedge
\check{u}^{k-1}}
\bigl\{(1-r)\mu\bigl(s,Y^{r,u}_s,\sigma^n_u
\bigr)+r\mu\bigl(\check {u}^{k-1},Y^{r,u}_{\check{u}^{k-1}},
\sigma^n_u\bigr)\bigr\}\,\mathrm{d}s
\\
&&{}+\sum_{k=1}^{L_0}\int
^{t\wedge\check{u}^{k}}_{t\wedge\check
{u}^{k-1}}\bigl\{(1-r)b\bigl(s,Y^{r,u}_s,
\sigma^n_u\bigr) +rb\bigl(\check {u}^{k-1},Y^{r,u}_{\check{u}^{k-1}},
\sigma^n_u\bigr)\bigr\}\,\mathrm{d}W_s,\qquad  t\in[0,T].
\end{eqnarray*}
Then we have $Y^{0,0}\equiv Y$.

Moreover, we define $\check{p}^r_{k,u}(z_0,z_1)=p(z_1;z_0,r,\check
{u}^{k-1},\check{u}^k,\sigma^n_u)$, $\check
{p}^{r,(1)}_{k,u}(z_0,z_1)=\partial_{\sigma}p(z_1;z_0,r,\allowbreak \check
{u}^{k-1},\check{u}^k,\sigma^n_u)$,
\[
\mathbb{P}^r_u(z,\bar{u})=\prod_{k=1}^{L_0}
\check {p}^r_{k,u}\bigl((x_{k-1},y_{k-1}),(x_k,y_k)
\bigr), \qquad \bar{\mathbb {P}}^r_u(z_0,\bar{z},
\bar{u})=\int\mathbb{P}^r_u(z;\bar{u})\,\mathrm{d}\hat{z}
\]
and $P^r_u=\mathbb{P}^r_u(z,\bar{u})P_{Y_0}(\mathrm{d}z_0)\,\mathrm{d}\bar{z}\,\mathrm{d}\hat{z}$
for $z=((x_k)_{k=0}^{L_0}, (y_k)_{k=0}^{L_0})\in\mathbb{R}^{2L_0+2}$.

Then synchronous observations $Y^{r,u}_{\check{u}}$ follow the
distribution $P^r_u$.
Moreover, we have
\begin{eqnarray*}
P_{\sigma^n_u,n} &=& P_{(\Pi,Y^{0,u}_{\Pi})}=P_{Y^{0,u}_{\bar
{u}}}(\mathrm{d}z_0\,\mathrm{d}
\bar{z}|\Pi=\bar{u})P_{\Pi}(\mathrm{d}\bar{u}) =P_{Y^{0,u}_{\bar{u}}}(
\mathrm{d}z_0\,\mathrm{d}\bar{z})P_{\Pi}(\mathrm{d}\bar{u})\\
&=& \bar{\mathbb{P}}^0_u(z_0,
\bar{z},\bar{u})P_{Y_0}(\mathrm{d}z_0)\,\mathrm{d}\bar{z}P_{\Pi
}(\mathrm{d}\bar{u}).
\end{eqnarray*}
Therefore, we obtain
%
\begin{equation}
\label{ratio-equ1}
\log\frac{\mathrm{d}P_{\sigma^n_u,n}}{\mathrm{d}P_{\sigma_{\ast},n}}(z_0,\bar{z},\bar{u}) =\log
\frac{\bar{\mathbb{P}}^0_u}{\bar{\mathbb{P}}^0_0}(z_0,\bar{z},\bar{u}).
\end{equation}
So it is sufficient to investigate the asymptotic behavior of $\log
(\bar{\mathbb{P}}^0_u/\bar{\mathbb{P}}^0_0)$.

For each function with respect to $(z,\bar{u})$ or $(z_0,\bar{z},\bar
{u})$, we often omit the variable $\bar{u}$.

The following theorem gives a LAMN-type property of $H_n$ (Proposition~3 and Proposition~10 in Ogihara and Yoshida \cite{ogi-yos}).

%
\begin{theorem}\label{ogi-yos-thm}
Assume $[A1'],[A2]$ and $[A3]$. Then there exists a random variable
$\mathcal{N}$ on an extension of $(\Omega,\mathcal{F},P)$ such that
$\mathcal{N}$ is independent of $\mathcal{F}$, $-b_n^{-1}\partial
_{\sigma}^2H_n(\sigma_{\ast})\circ(\Pi,Y_{\Pi})\to^p \Gamma$,
\[
\bigl\{H_n\bigl(\sigma^n_u
\bigr)-H_n(\sigma_{\ast})-\bigl(u^{\star}b_n^{-1/2}
\partial _{\sigma}H_n(\sigma_{\ast})+u^{\star}b_n^{-1}
\partial_{\sigma
}^2H_n(\sigma_{\ast})u/2
\bigr) \bigr\}\circ(\Pi,Y_{\Pi})\to^p 0,
\]
$b_n^{-1/2}\partial_{\sigma}H_n(\sigma_{\ast})\circ(\Pi,Y_{\Pi})\to
^{s\mbox{-}\mathcal{L}} \Gamma^{1/2}\mathcal{N}$
as $n\to\infty$, where $\Gamma$ is defined by (\ref{gamma-def}).
\end{theorem}

\begin{remark}
Though we need an assumption ``$\partial^k_{\sigma}b \ (0\leq k\leq4)$
can be extended to a continuous function on $[0,T]\times\mathbb
{R}^2\times\bar{\Lambda}$''
to apply the results in Ogihara and Yoshida \cite{ogi-yos},
the assumption can be removed by considering a relatively compact open
subset of $\Lambda$ containing $\sigma_{\ast}$.
\end{remark}

By virtue of  Lemma~\ref{A1toA1}, Theorem~\ref{ogi-yos-thm}
and (\ref{ratio-equ1}), to obtain Theorem~\ref{main}, it is sufficient
to show asymptotic equivalence
of $\log(\bar{\mathbb{P}}^0_u/\bar{\mathbb{P}}^0_0)(Y_{\Pi})$ and
$(H_n(\sigma^n_u)-H_n(\sigma_{\ast}))\circ(\Pi, Y_{\Pi})$
under $[A1']$, $[A2]$ and $[A3]$.
We will prove it in the rest of this paper.

\subsection{Malliavin calculus techniques and estimates for transition
densities}\label{malliavin-section}

We will prepare results of estimates for transition density functions
used later. To this end, we introduce some techniques from Malliavin calculus.
We refer the reader to Chapter II in Nualart \cite{nua} and Gobet \cite{gob} for detailed expositions of this subsection.

We\vspace*{1pt} fix $u\in\mathbb{R}^d$, $\bar{u}\in\mathcal{U}$, $1\leq k\leq
L_0(\bar{u})$ and $n\geq n_u$ here.
For $0\leq r\leq1$ and $x\in\mathbb{R}^2$, consider a~unique
two-dimensional process $\{\mathbf{Y}^{r,u,k,x}_t\}_{t\in[0,\Delta\check
{u}^k]}=\{(\mathbf{Y}^{r,u,k,x,1}_t,\mathbf{Y}^{r,u,k,x,2}_t)\}_{t\in
[0,\Delta\check{u}^k]}$ satisfying
\begin{eqnarray*}
\mathbf{Y}^{r,u,k,x,i}_t &=& x+\int^t_0
\bigl\{(1-r)\mu^{(0),r,i}_s+r\mu^i\bigl(\check
{u}^{k-1},x,\sigma^n_u\bigr) \bigr\}\,
\mathrm{d}s \\
&&{}+\sum_{j=1}^2\int
^t_0 \bigl\{ (1-r)b^{(0),r,i,j}_s+rb^{ij}
\bigl(\check{u}^{k-1},x,\sigma^n_u\bigr) \bigr
\}\,\mathrm{d}W^j_s
\end{eqnarray*}
for $t\in[0,\Delta\check{u}^k]$,
where $\mu^{(q),r,i}_{t,p_1,\ldots, p_q}=(\partial^q_x\mu^i(t+\check
{u}^{k-1},\mathbf{Y}^{r,u,k,x}_t,\sigma^n_u))_{p_1,\ldots, p_q}$,
$b^{(q),r,i,j}_{t,p_1,\ldots, p_q}= (\partial^q_xb^{ij}(t+\check
{u}^{k-1},\mathbf{Y}^{r,u,k,x}_t,\sigma^n_u))_{p_1,\ldots, p_q}$
for $q\in\mathbb{Z}_+$. We simply denote $\mathbf{Y}^{r}_t=\mathbf{Y}^{r,u,k,x}_t$ and $\mathbf{Y}^{r,i}_t=\mathbf{Y}^{r,u,k,x,i}_t$.

Under $[A1']$, Theorem~39 in Chapter V of Protter \cite{pro} ensures
that $\partial_r\mathbf{Y}^r_t=(\partial_r\mathbf{Y}^{r,1}_t,\partial_r\mathbf{Y}^{r,2}_t)$ exists for any $t\in[0,\Delta\check{u}^k]$ a.s. and satisfies
\begin{eqnarray*}
\partial_r\mathbf{Y}^{r,i}_t&=&\int
^t_0 \biggl[\sum_p(1-r)
\mu ^{(1),r,i}_{s,p}\partial_r\mathbf{Y}^{r,p}_s+
\mu^i\bigl(\check{u}^{k-1},x,\sigma ^n_u
\bigr)-\mu^{(0),r,i}_s \biggr]\,\mathrm{d}s
\\
&&{}+\sum_{j}\int^t_0
\biggl[\sum_p(1-r)b^{(1),r,i,j}_{s,p}
\partial_r\mathbf{Y}^{r,p}_s+b^{ij}
\bigl(\check{u}^{k-1},x,\sigma^n_u
\bigr)-b^{(0),r,i,j}_s \biggr]\,\mathrm{d}W^j_s.
\end{eqnarray*}

Define an isonormal Gaussian process $W$
by $W(\xi)=\int^{\Delta\check{u}^k}_0(\mathrm{d}\xi_t/\mathrm{d}t)\cdot \mathrm{d}W_{t+\check
{u}^{k-1}}$ for an $\mathbb{R}^2$-valued absolutely continuous function\vspace*{1pt}
$\xi=\{\xi_t\}_{0\leq t\leq\Delta\check{u}^k}$
satisfying $\int^{\Delta\check{u}^k}_0|\mathrm{d}\xi_t/\mathrm{d}t|^2\,\mathrm{d}t<\infty$.
We also consider the Malliavin derivative operator $D$ and the
divergence operator $\delta$.

Let $\{\mathsf{V}^{r}_t\}_{t\in[0,\Delta\check{u}^k]}=\{\mathsf{V}^{r,i,j}_t\}_{t\in[0,\Delta\check{u}^k],i,j}$ be a stochastic
process satisfying
\[
\mathsf{V}^{r,i,j}_t=\delta_{ij}+ \sum
_p\int^t_0(1-r)\mu
^{(1),r,i}_{s,p}\mathsf{V}^{r,p,j}_s \,
\mathrm{d}s +\sum_{p,q}\int^t_0(1-r)b^{(1),r,i,q}_{s,p}{
\sf V}^{r,p,j}_s \,\mathrm{d}W^q_s,
\]
then the argument in Sections~2.2 and 2.3 of Nualart \cite{nua} yields
\[
D^j_t\mathbf{Y}^{r,i}_{\Delta\check{u}^k}=\sum
_{p,q}\mathsf{V}^{r,i,q}_{\Delta\check{u}^k}\bigl(
\bigl(\mathsf{V}^r_t\bigr)^{-1}
\bigr)_{qp} \bigl[(1-r)b^{(0),r,p,j}_t+rb^{pj}
\bigl(\check{u}^{k-1},x,\sigma^n_u\bigr) \bigr]
\]
for $t\in[0,\Delta\check{u}^k]$, where $((\mathsf{V}^r_t)^{-1})_{qp}$
represents the element of $(\mathsf{V}^r_t)^{-1}$.
Moreover, we obtain
\[
\sup_{t\in[0,\Delta\check{u}^k]}\bigl(E\bigl[\bigl|D_t\mathbf{Y}^r_{\Delta\check
{u}^k}\bigr|^M
\bigr]\vee E\bigl[\bigl|D_t\partial_r\mathbf{Y}^r_{\Delta\check
{u}^k}\bigr|^M
\bigr]\bigr)<\infty
\]
for $M>0$.

\begin{lemma}\label{p-est}
Let $u\in\mathbb{R}^d$ and $q\geq1$. Assume $[A1']$. Then there exists
a positive constant $C_q$ such that
\[
\sup_k\frac{1}{(\Delta\check{u}^k)^{q/2}}\sup_{0\leq r\leq
1,z_{k-1}}\int \biggl|
\frac{\partial_r\check{p}^r_{k,u}}{\check
{p}^r_{k,u}} \biggr|^q\check{p}^r_{k,u}(z_{k-1},z_k)\,\mathrm{d}z_k
\leq C_q
\]
for any $n\geq n_u$ and $\bar{u}\in\mathcal{U}$.
\end{lemma}

\begin{pf}
Let
$\mathcal{B}^r_t=\{\mathcal{B}^{r,i}_{t,j}\}
_{i,j}=((1-r)b^{(0),r}_t+rb(\check{u}^{k-1},x,\sigma
^n_u))^{-1}V^r_t(V^r_{\Delta\check{u}^k})^{-1}$ for $t\in[0,\Delta
\check{u}^k]$,
where $b^{(0),r}_t=\{b^{(0),r,i,j}_t\}_{ij}$.
Then by a similar argument to the proof of Proposition~4.1. in Gobet
\cite{gob}, we obtain
%
\begin{equation}
\label{density-rep}
\frac{\partial_r\check{p}^r_{k,u}}{\check{p}^r_{k,u}}(z_{k-1},z_k)=
\frac{1}{\Delta\check{u}^k}E\bigl[\delta\bigl(\bigl(\mathcal{B}^r
\bigr)^{\star}\partial_r\mathbf{Y}^{r,u,k,z_{k-1}}_{\Delta\check{u}^k}
\bigr)|\mathbf{Y}^{r,u,k,z_{k-1}}_{\Delta\check{u}^k}=z_k\bigr].
\end{equation}

Moreover, by Proposition~1.3.3. in Nualart \cite{nua}, we obtain
%
\begin{equation}
\label{p-est-eq1}
\delta\bigl(\bigl(\mathcal{B}^r\bigr)^{\star}
\partial_r\mathbf{Y}^r_{\Delta\check{u}^k}\bigr) =\sum
_{i=1}^2 \biggl\{\partial_r
\mathbf{Y}^{r,i}_{\Delta\check{u}^k}\delta \bigl(\mathcal{B}^{r,i}
\bigr) -\int^{\Delta\check{u}^k}_0 D_t\bigl(
\partial_r\mathbf{Y}^{r,i}_{\Delta\check
{u}^k}\bigr) \cdot
\mathcal{B}^{r,i}_t\,\mathrm{d}t \biggr\}.
\end{equation}
Furthermore, we have
%
\begin{eqnarray}
&& E \Biggl[ \Biggl|\sum_{i=1}^2
\int^{\Delta\check{u}^k}_0 D_t\bigl(
\partial_r\mathbf{Y}^{r,i}_{\Delta\check{u}^k}\bigr) \cdot
\mathcal{B}^{r,i}_t\,\mathrm{d}t \Biggr|^q \Biggr]
\nonumber\\
&&\label{p-est-eq2} \quad \leq   C\bigl(\Delta\check{u}^k\bigr)^{q-1}\sum
_{i=1}^2\int^{\Delta\check{u}^k}_0
E\bigl[\bigl|D_t\partial_r\mathbf{Y}^{r,i}_{\Delta\check{u}^k}\bigr|^{2q}
\bigr]^{{1}/{2}}E\bigl[\bigl|\mathcal{B}^{r,i}_t\bigr|^{2q}
\bigr]^{{1}/{2}}\,\mathrm{d}t
\\
&& \quad = \mathrm{O}\bigl(\bigl(\Delta\check {u}^k
\bigr)^{3q/2}\bigr),
\nonumber
\end{eqnarray}
and
%
\begin{equation}
\label{p-est-eq3}
E \Biggl[ \Biggl|\sum_{i=1}^2
\partial_r\mathbf{Y}^{r,i}_{\Delta\check
{u}^k}\delta\bigl(
\mathcal{B}^{r,i}\bigr) \Biggr|^q \Biggr] =\mathrm{O}\bigl(\bigl(\Delta
\check{u}^k\bigr)^q\bigr)\times\sum
_{i=1}^2E\bigl[\bigl|\delta\bigl(\mathcal
{B}^{r,i}\bigr)\bigr|^{2q}\bigr]^{{1}/{2}},
\end{equation}
where we use the fact  that any moments of $(\mathsf{V}^r_{\Delta
\check{u}^k})^{-1}$ are bounded  (see Section~2.3.1. in Nualart~\cite{nua}).

By Propositions 1.3.8. and 1.5.7. in Nualart \cite{nua} and
Clark--Ocone representation formula (Corollary A.2. in Nualart and
Pardoux \cite{nua-par}), we have
%
\begin{eqnarray}
 E\bigl[\bigl|\delta\bigl(\mathcal{B}^{r,i}
\bigr)\bigr|^{2q}\bigr]&=&E \biggl[ \biggl|\int^{\Delta\check
{u}^k}_0
E\bigl[D_t\delta\bigl(\mathcal{B}^{r,i}\bigr)|
\mathcal{F}_{t+\check
{u}^{k-1}}\bigr]\cdot \mathrm{d}W_{t+\check{u}^{k-1}} \biggr|^{2q}
\biggr]
\nonumber
\\
\label{p-est-eq4} &\leq & \bigl(\Delta\check{u}^k\bigr)^{q-1}\int
^{\Delta\check{u}^k}_0 E\bigl[\bigl|D_t\delta\bigl(
\mathcal{B}^{r,i}\bigr)\bigr|^{2q}\bigr]\,\mathrm{d}t
\\
\nonumber
& =& \bigl(
\Delta\check{u}^k\bigr)^{q-1}\int^{\Delta\check{u}^k}_0
E\bigl[\bigl|\mathcal {B}^{r,i}_t+\delta\bigl(D_t
\mathcal{B}^{r,i}\bigr)\bigr|^{2q}\bigr]\,\mathrm{d}t
\\
&=& \mathrm{O}\bigl(\bigl(\Delta
\check{u}^k\bigr)^q\bigr).\nonumber
\end{eqnarray}

By (\ref{p-est-eq1})--(\ref{p-est-eq4}), we obtain $E[|\delta((\mathcal
{B}^r)^{\star}\partial_r\mathbf{Y}^r_{\Delta\check{u}^k})|^q]=\mathrm{O}((\Delta
\check{u}^k)^{3q/2})$.
Therefore, we have
\[
\int \biggl|\frac{\partial_r\check{p}^r_{k,u}}{\check
{p}^r_{k,u}} \biggr|^q\check{p}^r_{k,u}(z_{k-1},z_k)\,\mathrm{d}z_k
\leq C_q\bigl(\Delta\check{u}^k\bigr)^{q/2}
\]
by (\ref{density-rep}).
\end{pf}

The following lemma is proved similarly.

%
\begin{lemma}\label{density-est2}
Let $u\in\mathbb{R}^d$ and $q\geq1$. Assume $[A1']$. Then
\begin{eqnarray*}
&&\frac{\partial_{\sigma}p}{p}\bigl(z_k;z_{k-1},r,\check{u}^{k-1},
\check {u}^k,\sigma^n_u\bigr) =
\frac{1}{\Delta\check{u}^k}E\bigl[\delta\bigl(\bigl(\mathcal{B}^r
\bigr)^{\star}\partial _{\sigma}\mathbf{Y}^r_{\Delta\check{u}^k}
\bigr)|\mathbf{Y}^r_{\Delta\check
{u}^k}=z_k\bigr],
\\
&&\frac{\partial^2_{(r,\sigma)}p}{p}\bigl(z_k;z_{k-1},r,\check{u}^{k-1},
\check {u}^k,\sigma^n_u\bigr) \\
&& \quad =E \Biggl[
\frac{1}{\Delta\check{u}^k}\delta\bigl(\bigl(\mathcal{B}^r
\bigr)^{\star
}\partial^2_{(r,\sigma)}\mathbf{Y}^r_{\Delta\check{u}^k}
\bigr)
\\
&&\qquad\hspace*{12pt}{}+\frac{1}{(\Delta\check{u}^k)^2}\sum_{i=1}^2\delta
\bigl(\mathcal {B}^{r,i}\delta\bigl(\bigl(\mathcal{B}^r
\bigr)^{\star}\partial_{(r,\sigma)}\mathbf{Y}^r_{\Delta\check{u}^k}
\partial_{(r,\sigma)}\mathbf{Y}^{r,i}_{\Delta
\check{u}^k}\bigr)\bigr) \Big|
\mathbf{Y}^r_{\Delta\check{u}^k}=z_k \Biggr],
\end{eqnarray*}
and there exists a constant $C_q>0$ such that
\begin{eqnarray*}
&& \sup_{k,0\leq r'\leq1,0\leq v\leq1, z_{k-1}}\bigl(\Delta\check {u}^k
\bigr)^{-{lq}/{2}}\int\bigl|\partial^j_{\sigma}
\partial^l_r\log p\bigr|^q\bigl(z_k;z_{k-1},r',
\check{u}^{k-1},\check{u}^k,\sigma^n_{vu}
\bigr)\\
&&\qquad \hspace*{111pt}{}\times \check{p}^{r'}_{k,vu}(z_{k-1},z_k)\,\mathrm{d}z_k
\leq C_q
\end{eqnarray*}
for $r\in[0,1]$, $\bar{u}\in\mathcal{U}$, $n\geq n_u$, $1\leq j+l\leq2$.
\end{lemma}

\subsection{Tightness results of some log-likelihood ratios}\label{tightness-section}

In this section,\hspace*{1pt} we will prove some tightness results, which are
necessary later.
First, we prove tightness of $\{\sup_{0\leq r\leq1}|\log(\bar{\mathbb
{P}}^r_u/\bar{\mathbb{P}}^0_0)|(Y_{\Pi})\}_n$.\vspace*{1pt}
To this end, we prove results about the log-likelihood ratio $\log
(\mathbb{P}^r_u/\mathbb{P}^0_0)(Y_{\Pi})$.
Then we prove a key proposition (Proposition~\ref{synch-to-nonsynch})
which enables us to  obtain tightness of a density ratio in a
nonsynchronous scheme
from properties of a density ratio in a synchronous scheme.

\begin{enumerate}[{[$A3'$]}]
\item[{[$A3'$]}] The sequence $\{b_n^{-1}(\ell_{1,n}+\ell_{2,n})\}_n$ is tight.
\end{enumerate}
Since $\operatorname{tr}(\mathcal{E}^1(T))+\operatorname{tr}(\mathcal{E}^2(T))=\ell
_{1,n}+\ell_{2,n}$, $[A3]$ implies $[A3']$.

We prepare some results for the log-likelihood ratio $\log(\mathbb
{P}^r_u/\mathbb{P}^0_0)$.

\begin{lemma}\label{synchro-tight-lemma1}
Let $u\in\mathbb{R}^d$. Assume $[A1']$ and $[A3']$. Then for any
$\varepsilon>0$, there exists $M>0$ such that
\[
\sup_{n\geq n_u} \biggl\{P \biggl[ \biggl|\log\frac{\mathbb{P}^0_u}{\mathbb
{P}^0_0} \biggr|
\bigl(Y^{0,u}_{\check{U}}\bigr)>M \biggr]\vee P \biggl[ \biggl|\log
\frac
{\mathbb{P}^0_u}{\mathbb{P}^0_0} \biggr|\bigl(Y^{0,0}_{\check{U}}\bigr)>M \biggr] \biggr
\}<\varepsilon.
\]
\end{lemma}

\begin{pf}
By Lemmas \ref{density-est2} and \ref{density-est}, we obtain
\begin{eqnarray*}
&& E \biggl[ \biggl|\log\frac{\mathbb{P}^0_u}{\mathbb{P}^0_0} \biggr|\bigl(Y^{0,u}_{\check{U}}\bigr)
\Big|\Pi \biggr]\\
 && \quad \leq E \biggl[\int^1_0 \biggl|
\frac{\partial_{t}(\mathbb
{P}^0_{tu})}{\mathbb{P}^0_{tu}} \biggr|\bigl(Y^{0,u}_{\check{U}}\bigr)\,\mathrm{d}t \Big|
\Pi \biggr]
\\
&&\quad \leq E \biggl[b_n^{-{1}/{2}}|u| \biggl|\sum
_k\frac{\check
{p}^{0,(1)}_{k,u}}{\check{p}^0_{k,u}}\bigl(Y^{0,u}_{\check
{U}^{k-1}},Y^{0,u}_{\check{U}^k}
\bigr) \biggr|
\\
&&\hspace*{12pt}\qquad {}+b_n^{-1}|u|^2\int^1_0\!\!\int^1_{t_1} \biggl|\sum_k
\partial_{\sigma} \biggl(\frac{\partial_{\sigma}p}{p} \biggr) \bigl(Y^{0,u}_{\check
{U}^k};Y^{0,u}_{\check{U}^{k-1}},0,
\check{U}^{k-1},\check{U}^k,\sigma ^n_{t_2u}
\bigr) \biggr|\,\mathrm{d}t_2\,\mathrm{d}t_1 \Big| \Pi \biggr]
\\
&&\quad \leq Cb_n^{-1/2}(\ell_{1,n}+
\ell_{2,n})^{1/2}+Cb_n^{-1}(\ell
_{1,n}+\ell_{2,n}).
\end{eqnarray*}
Hence, by Lemma~\ref{cond-conv-to-conv} and the assumptions, for any
$\varepsilon>0$ there exists $M>0$ such  that
$\sup_{n\geq n_u}P[|\log(\mathbb{P}^0_u/\mathbb
{P}^0_0)|(Y^{0,u}_{\check{U}})>M]<\varepsilon$.
Similarly, we obtain $\sup_{n\geq n_u}P[|\log(\mathbb{P}^0_u/\break \mathbb
{P}^0_0)| (Y^{0,0}_{\check{U}})>M]<\varepsilon$.
\end{pf}

We define
\[
A^n_M(\bar{u})=\Bigl\{(x,y)\in\mathbb{R}^{L_0(\bar{u})+1}
\times\mathbb {R}^{L_0(\bar{u})+1};\sup_{0\leq r\leq1}\bigl|\log\bigl(
\mathbb{P}^r_u/\mathbb {P}^0_0
\bigr)\bigr|(x,y)\leq M\Bigr\}
\]
for $\bar{u}\in\mathcal{U}$ and $M>0$.

\begin{lemma}\label{synchro-tight-lemma2}
Let $u\in\mathbb{R}^d$. Assume $[A1']$ and $[A3']$. Then for any
$\varepsilon>0$, there exists $M>0$ such that
\[
\sup_{n\geq n_u,r} \biggl\{P\bigl[Y^{0,0}_{\check{U}}\in
\bigl(A^n_M\bigr)^c(\Pi)\bigr]\vee P
\bigl[Y^{r,u}_{\check{U}}\in\bigl(A^n_M
\bigr)^c(\Pi)\bigr]\vee E \biggl[ \biggl|\frac{\partial
_r\mathbb{P}^r_u}{\mathbb{P}^r_u}
\biggr|1_{(A^n_M)^c}\bigl(Y^{r,u}_{\check
{U}}\bigr) \biggr] \biggr\}<
\varepsilon.
\]
\end{lemma}

\begin{pf}
Fix $\varepsilon>0$. Then for $r\in[0,1]$, we obtain
\[
E \biggl[\sup_{0\leq r'\leq1} \biggl|\log\frac{\mathbb{P}^{r'}_u}{\mathbb
{P}^0_u} \biggr|
\bigl(Y^{r,u}_{\check{U}}\bigr) \Big|\Pi \biggr]\leq E \biggl[\int
^1_0 \biggl|\frac{\partial_r\mathbb{P}^{r'}_u}{\mathbb{P}^{r'}_u} \biggr|\bigl(Y^{r,u}_{\check{U}}
\bigr)\,\mathrm{d}r' \Big|\Pi \biggr] =\int^1_0E
\biggl[ \biggl|\frac{\partial_r\mathbb{P}^{r'}_u}{\mathbb
{P}^{r'}_u} \biggr|\bigl(Y^{r,u}_{\check{U}}\bigr) \Big|\Pi
\biggr]\,\mathrm{d}r'.
\]

On the other hand,
by Lemmas \ref{density-est2} and \ref{density-est}, we have
\begin{eqnarray*}
E \biggl[ \biggl|\frac{\partial_r\mathbb{P}^{r'}_u}{\mathbb{P}^{r'}_u} \biggr|\bigl(Y^{r,u}_{\check{U}}\bigr) \Big|\Pi
\biggr] &=&E \biggl[ \biggl|\frac{\partial_r\mathbb{P}^{r}_u}{\mathbb{P}^{r}_u}+\int^{r'}_r
\partial_r \biggl(\frac{\partial_r\mathbb{P}^{r_1}_u}{\mathbb
{P}^{r_1}_u} \biggr)\,\mathrm{d}r_1 \biggr|
\bigl(Y^{r,u}_{\check{U}}\bigr) \Big|\Pi \biggr]
\\
&\leq& E \Biggl[\sum_{k=1}^{L_0} \biggl(
\frac{\partial_r\check
{p}^r_{k,u}}{\check{p}^r_{k,u}} \biggr)^2\bigl(Y^{r,u}_{\check
{U}^{k-1}},Y^{r,u}_{\check{U}^k}
\bigr) \Big|\Pi \Biggr]^{{1}/{2}}\\
&&{} +\sup_{r'}E \Biggl[\sum
_{k=1}^{L_0} \biggl|\partial_r \biggl(
\frac{\partial
_r\check{p}^{r'}_{k,u}}{\check{p}^{r'}_{k,u}} \biggr) \biggr|\bigl(Y^{r,u}_{\check{U}^{k-1}},Y^{r,u}_{\check{U}^k}
\bigr) \Big|\Pi \Biggr]
\\
&
\leq&C\sqrt{T}+\mathit{CT}.
\end{eqnarray*}

Hence, by Lemma~\ref{cond-conv-to-conv}, for any $\varepsilon>0$ there
exists $M_1>0$ such that
\[
\sup_{0\leq r\leq1}P \biggl[\sup_{0\leq r'\leq1} \biggl|\log
\frac{\mathbb
{P}^{r'}_u}{\mathbb{P}^0_u} \biggr|\bigl(Y^{r,u}_{\check{U}}\bigr)>
\frac
{M_1}{2} \biggr]<\frac{\varepsilon}{2}.
\]

Therefore, Lemma~\ref{synchro-tight-lemma1} yields
%
\begin{eqnarray}
&&\hspace*{-18pt}\sup_{n\geq n_u,r}P \bigl[Y^{r,u}_{\check{U}}
\in\bigl(A_M^n\bigr)^c(\Pi) \bigr]
\nonumber
\\
&&\hspace*{-18pt} \quad \leq \sup_{n\geq n_u,r}P \biggl[ \biggl|\log\frac{\mathbb{P}^0_u}{\mathbb
{P}^0_0} \biggr|
\bigl(Y^{r,u}_{\check{U}}\bigr)>\frac{M}{2}  \mbox{ and } \sup
_{0\leq r'\leq1} \biggl|\log\frac{\mathbb{P}^{r'}_u}{\mathbb{P}^0_u} \biggr|\bigl(Y^{r,u}_{\check{U}}
\bigr)\leq\frac{M_1}{2} \biggr]+\frac{\varepsilon}{2}
\nonumber
\\[-8pt]
\label{synchro-tight}
\\[-8pt]
\nonumber
&&\hspace*{-18pt}\quad \leq \sup_{n\geq n_u,r}E \biggl[\frac{\mathbb{P}^r_u}{\mathbb
{P}^0_u}
\bigl(Y^{0,u}_{\check{U}}\bigr), \biggl|\log\frac{\mathbb{P}^0_u}{\mathbb
{P}^0_0} \biggr|
\bigl(Y^{0,u}_{\check{U}}\bigr)>\frac{M}{2}  \mbox{ and } \sup
_{0\leq r'\leq1} \biggl|\log\frac{\mathbb{P}^{r'}_u}{\mathbb{P}^0_u} \biggr|\bigl(Y^{0,u}_{\check{U}}
\bigr)\leq\frac{M_1}{2} \biggr]+\frac{\varepsilon}{2}
\\
&&\hspace*{-18pt}\quad \leq  \mathrm{e}^{M_1/2}\sup_{n\geq n_u}P\biggl[\biggl|\log \frac{\mathbb{P}^0_u}{\mathbb {P}^0_0} \biggr|\bigl(Y^{0,u}_{\check{U}}\bigr) > \frac{M}{2}\biggr]+
\frac{\varepsilon}{2}<
\varepsilon\nonumber
\end{eqnarray}
for sufficiently large $M>0$.

Moreover, by Lemma~\ref{p-est}, we obtain
\[
\sup_{n\geq n_u,r}E\bigl[\bigl|\partial_r
\mathbb{P}^r_u/\mathbb {P}^r_u\bigr|^2
\bigl(Y^{r,u}_{\check{U}}\bigr)\bigr]=\sup_{n\geq n_u,r}E
\biggl[\sum_k\bigl(\partial _r
\check{p}^r_{k,u}/\check{p}^r_{k,u}
\bigr)^2\bigl(Y^{r,u}_{\check
{U}^{k-1}},Y^{r,u}_{\check{U}^k}
\bigr)\biggr]<\infty.
\]
Hence, by (\ref{synchro-tight}), we have $\sup_{n\geq n_u,r}E[|\partial
_r\mathbb{P}^r_u/\mathbb{P}^r_u|1_{(A^n_{M})^c}(Y^{r,u}_{\check
{U}})]<\varepsilon$
for sufficiently large $M>0$.

On the other hand, there exists $M_2>0$ such that $\sup_{n\geq
n_u}P[|\log(\mathbb{P}^0_u/\mathbb{P}^0_0)|(Y^{0,0}_{\check
{U}})>M_2]<\varepsilon/2$ by Lemma~\ref{synchro-tight-lemma1}.
Therefore,\vspace*{-2pt} we have
\begin{eqnarray*}
P\bigl[Y^{0,0}_{\check{U}}\in\bigl(A^n_{M}
\bigr)^c(\Pi)\bigr]&\leq&P \biggl[\sup_{0\leq
r\leq1} \biggl|\log
\frac{\mathbb{P}^r_u}{\mathbb{P}^0_0} \biggr|\bigl(Y^{0,0}_{\check{U}}\bigr)>M,  \biggl|\log
\frac{\mathbb{P}^0_u}{\mathbb
{P}^0_0} \biggr|\bigl(Y^{0,0}_{\check{U}}\bigr)\leq
M_2 \biggr]+\frac{\varepsilon}{2}
\\[-2pt]
&\leq& \mathrm{e}^{M_2}P \biggl[\sup_{0\leq r\leq1} \biggl|\log
\frac{\mathbb
{P}^r_u}{\mathbb{P}^0_0} \biggr|\bigl(Y^{0,u}_{\check{U}}\bigr)>M \biggr]+
\frac
{\varepsilon}{2}<\varepsilon\vspace*{-3pt}
\end{eqnarray*}
for sufficiently large\vspace*{-1pt} $M>0$ by (\ref{synchro-tight}).
\end{pf}

Let $\mathsf{Z}^n=\{\mathsf{Z}^n_t\}_{0\leq t\leq T}$ and $\mathsf{Z}^{n,r}=\{
\mathsf{Z}^{n,r}_t(\bar{u})\}_{0\leq t\leq T}$ be two-dimensional
continuous $\mathbf{F}$-adapted processes
satisfying that $\mathsf{Z}^{n,r}_0=\mathsf{Z}^n_0$ for $n\in\mathbb{N}$,
$\bar{u}\in\mathcal{U}$ and $0\leq r\leq1$,
$(t,\bar{u},\omega)\mapsto Z^{n,r}_t(\bar{u})(\omega)$ is
measurable, and $(Z^n_t,Z^{n,r}_t(\bar{u}))_{t,r,\bar{u}}$ are
independent of $\sigma((\Pi_n)_n)$.
Let the distributions of $\mathsf{Z}^n_{\check{u}}$ and $\mathsf{Z}^{n,r}_{\check{u}}(\bar{u})$ be given by $F_n(z_0,\bar{z},\hat
{z})P_{Z^n_0}(\mathrm{d}z_0)\,\mathrm{d}\bar{z}\,\mathrm{d}\hat{z}$
and $F^r_n(z_0,\bar{z},\hat{z})P_{Z^n_0}(\mathrm{d}z_0)\,\mathrm{d}\bar{z}\,\mathrm{d}\hat{z}$, respectively,
for some positive-valued Borel functions $F_n$ and $F^r_n$.
Let $K^n_M=K^n_M(\bar{u})$ be a Borel set in $\mathbb{R}^{2(L_0+1)}$,
$\bar{K}^n_M=\bar{K}^n_M(z_0,\bar{z},\bar{u})=\{\hat{z};(z_0,\bar
{z},\hat{z})\in K^n_M\}$,
$\bar{F}_n(z_0,\bar{z})=\int F_n(z_0,\bar{z},\hat{z})\,\mathrm{d}\hat{z}, \bar
{F}_{n,M}(z_0,\bar{z})= \int_{\bar{K}^n_M(z_0,\bar{z})} F_n(z_0,\bar
{z},\hat{z})\,\mathrm{d}\hat{z}$
and $\bar{F}^r_n$ and $\bar{F}^r_{n,M}$ be\vspace*{1pt} defined similarly.

The following proposition is a key result to deduce properties of
density ratios in the nonsynchronous\vspace*{-6pt} scheme.

%
\begin{proposition}\label{synch-to-nonsynch}
Suppose $F^r_n$ can be continuously differentiable with respect to $r$
and  $\int\partial_rF^r_n 1_{(K^n_M)^c}\,\mathrm{d}\hat{z}$ exists and is
continuous with respect to $r$ for each $n,z_0,\bar{z}$ and $\hat{z}$.
\begin{enumerate}[2.]
\item[1.] Suppose for any $\varepsilon>0$, there exists $M_1>0$ such\vspace*{-2pt} that
\begin{eqnarray*}
&& \sup_nP \biggl[\sup_r \biggl|\log
\frac{F^r_n}{F_n}\biggr|\bigl(\mathsf{Z}^n_{\check
{U}}\bigr)>M \biggr]<
\varepsilon \quad\mbox{and}\\[-2pt]
&& \sup_{n,r} \biggl\{P\bigl[\mathsf{Z}^{n,0}_{\check{U}}\in \bigl(K^n_M
\bigr)^c(\Pi)\bigr]\vee E \biggl[ \biggl|\frac{\partial_rF^r_n}{F^r_n}
\biggr|1_{(K_M^n)^c}\bigl(\mathsf{Z}^{n,r}_{\check{U}}\bigr) \biggr]
\biggr\}<\varepsilon
\end{eqnarray*}
for $M\geq M_1$.
Then for any $\varepsilon,\eta>0$, there exists $M_2>0$ such\vspace*{-2pt} that
\[
\sup_nP\Bigl[\sup_r\bigl|\log\bigl(
\bar{F}^r_n/\bar{F}^r_{n,M}
\bigr)\bigr|\bigl(\mathsf{Z}^n_{\Pi}\bigr)\geq \eta\Bigr] <
\varepsilon\vspace*{-3pt}
\]
for $M\geq M_2$.

\item[2.]
\vspace*{-3pt}
\begin{eqnarray*}
&& \sup_nP \biggl[\sup_r \biggl|\log
\frac{\bar{F}^r_n}{\bar{F}^r_{n,M}} \biggr|\bigl(Z^n_{\Pi}\bigr)\geq\eta \Big|\Pi
\biggr]\\[-3pt]
&&\quad \leq \frac{\mathrm{e}^{M'}}{1-\mathrm{e}^{-\eta}} \biggl\{P\bigl[Z^{n,0}_{\check{U}}
\in \bigl(K^n_M\bigr)^c|\Pi\bigr]+\sup
_rE \biggl[ \biggl|\frac{\partial_rF^r_n}{F^r_n} \biggr|1_{(K^n_M)^c}
\bigl(Z^{n,r}_{\check{U}}\bigr) \Big|\Pi \biggr] \biggr\}
\\[-3pt]
&&\qquad {}+P \biggl[\sup_r \biggl|\log\frac{F^r_n}{F_n} \biggr|
\bigl(Z^n_{\check
{U}}\bigr)>M' \Big|\Pi \biggr]\vspace*{-3pt}
\end{eqnarray*}
for any $\eta, M,M'>0$.
\end{enumerate}
\end{proposition}

\begin{pf}
We first prove 1. By the assumptions, for any $\varepsilon,\eta>0$, there
exist $M_1,M_2>0$ such that
\begin{eqnarray*}
 \sup_nP \biggl[\sup_r \biggl|\log
\frac{F^r_n}{F_n} \biggr|\bigl(\mathsf{Z}^n_{\check
{U}}
\bigr)>M_1 \biggr] &<& \frac{\varepsilon}{3}\quad \mbox{and}\\
\sup
_{n,r} \biggl\{ P\bigl[\mathsf{Z}^{n,0}_{\check{U}}
\in\bigl(K^n_M\bigr)^c(\Pi)\bigr]\vee E
\biggl[ \biggl|\frac
{\partial_rF^r_n}{F^r_n} \biggr|1_{(K_M^n)^c}\bigl(\mathsf{Z}^{n,r}_{\check
{U}}
\bigr) \biggr] \biggr\} &< &\frac{\varepsilon\eta'}{3\mathrm{e}^{M_1}}
\end{eqnarray*}
for $M\geq M_2$, where $\eta'=1-\mathrm{e}^{-\eta}$. Hence, we obtain
\begin{eqnarray*}
&&P\biggl[\sup_r\biggl|\log\frac{\bar{F}^r_n}{\bar{F}^r_{n,M}}\biggr|\bigl(\mathsf{Z}^n_{\Pi}
\bigr)\geq\eta \biggr]
\\
&& \quad \leq  P\biggl[\sup_r\biggl|1-\frac{\bar{F}^r_{n,M}}{\bar{F}^r_n}\biggr|\bigl(\mathsf{Z}^n_{\Pi}
\bigr)\geq \eta', \mathsf{Z}^n_{\check{U}}\in
L^n_{M_1}\biggr]+\frac{\varepsilon}{3}
\\
&& \quad \leq \frac{1}{\eta'}E \biggl[\int\sup_r \biggl|1-
\frac{\bar
{F}^r_{n,M}}{\bar{F}^r_n} \biggr|1_{L^n_{M_1}}(z_0,\bar{z},\hat
{z})F_n(z_0,\bar{z},\hat{z})P_{Z^n_0}(\mathrm{d}z_0)\,\mathrm{d}
\bar{z}\,\mathrm{d}\hat{z} \Big|_{\bar
{u}=\Pi} \biggr]+\frac{\varepsilon}{3}
\\
&&\quad \leq\frac{1}{\eta'}E \biggl[\int\sup_r
\frac{\int
F_n1_{L^n_{M_1}}(z_0,\bar{z},\hat{z})\,\mathrm{d}\hat{z}}{\bar{F}^r_n}\sup_r\bigl|\bar {F}^r_n-
\bar{F}^r_{n,M}\bigr|P_{Z^n_0}(\mathrm{d}z_0)\,\mathrm{d}\bar{z}
\Big|_{\bar{u}=\Pi} \biggr]+\frac{\varepsilon}{3},
\end{eqnarray*}
where $L_M^n=\{(z_0,\bar{z},\hat{z});\sup_r|\log(F^r_n/F_n)|(z_0,\bar
{z},\hat{z})\leq M\}$.

Since
\[
\sup_r\frac{\int F_n1_{L^n_{M_1}}(z_0,\bar{z},\hat{z})\,\mathrm{d}\hat{z}}{\bar
{F}^r_n}=\sup_r
\frac{1}{\bar{F}^r_n}\int\frac
{F_n}{F^r_n}F^r_n1_{L^n_{M_1}}(z_0,
\bar{z},\hat{z})\,\mathrm{d}\hat{z}\leq \mathrm{e}^{M_1},
\]
we obtain
\begin{eqnarray*}
&&P\Bigl[\sup_r\bigl|\log\bigl(\bar{F}^r_n/
\bar{F}^r_{n,M}\bigr)\bigr|\bigl(\mathsf{Z}^n_{\Pi}
\bigr)\geq\eta \Bigr]
\\
&& \quad \leq  \frac{\mathrm{e}^{M_1}}{\eta'}E \biggl[\int \biggl\{\bigl|\bar{F}^0_n-
\bar {F}^0_{n,M}\bigr|+\int^1_0\bigl|
\partial_r\bigl(\bar{F}^r_n-
\bar{F}^r_{n,M}\bigr)\bigr|\,\mathrm{d}r \biggr\} P_{Z^n_0}(\mathrm{d}z_0)\,\mathrm{d}
\bar{z} \Big|_{\bar{u}=\Pi} \biggr]+\frac{\varepsilon}{3}
\\
&& \quad \leq \frac{\mathrm{e}^{M_1}}{\eta'} \biggl\{P\bigl[\mathsf{Z}^{n,0}_{\check{U}}
\in \bigl(K^n_M\bigr)^c(\Pi)\bigr]+\sup
_rE \biggl[ \biggl|\frac{\partial_rF^r_n}{F^r_n} \biggr|1_{(K^n_M)^c}\bigl(\mathsf{Z}^{n,r}_{\check{U}}\bigr) \biggr] \biggr\}+\frac{\varepsilon
}{3}<
\varepsilon
\end{eqnarray*}
for $M\geq M_2$. Hence, we obtain 1.

The result in 2. is proved by a similar argument as above.
\end{pf}

Let
\[
\bar{A}^n_M(z_0,\bar{z})=\Bigl\{\hat{z};
\sup_r\bigl|\log\bigl(\mathbb{P}^r_u/
\mathbb {P}^0_0\bigr) (z_0,\bar{z},\hat{z})\bigr|
\leq M\Bigr\},\qquad  \bar{\mathbb {P}}^r_{M,u}(z_0,
\bar{z})=\int_{\bar{A}^n_M(z_0,\bar{z})}\mathbb {P}^r_u(z_0,
\bar{z},\hat{z})\,\mathrm{d}\hat{z}
\]
for $M>0$.

\begin{lemma}\label{pre-nonsynchro-tight-lemma}
Let $u\in\mathbb{R}^d$. Assume $[A1']$ and $[A3']$. Then for any
$\varepsilon, \eta>0$, there exists $M'>0$ such that
\[
\sup_{n\geq n_u}P \biggl[\sup_r \biggl|\log
\frac{\bar{\mathbb{P}}^r_u}{\bar
{\mathbb{P}}^r_{M,u}} \biggr|(Y_{\Pi})\geq\eta \biggr] < \varepsilon,\qquad \sup
_{n\geq n_u}P \biggl[ \biggl|\log\frac{\bar{\mathbb{P}}^0_0}{\bar
{\mathbb{P}}^0_{M,0}} \biggr|(Y_{\Pi})
\geq\eta \biggr] < \varepsilon
\]
for $M\geq M'$.
\end{lemma}

\begin{pf}
The results are obtained by using Proposition~\ref{synch-to-nonsynch}
and Lemma~\ref{synchro-tight-lemma2}.
The first inequality is obtained by setting
$\mathsf{Z}^n=Y$, $\mathsf{Z}^{n,r}=Y^{r,u}$ and $K^n_M=A^n_M$ in Proposition~\ref{synch-to-nonsynch}.
For the second inequality, set\vspace*{-2pt} $\mathsf{Z}^n=\mathsf{Z}^{n,r}=Y$ and $K^n_M=A^n_M$.
\end{pf}

\begin{proposition}\label{logp-conv2}
Let $u\in\mathbb{R}^d$. Assume $[A1']$ and $[A3']$.
Then $\{\sup_{0\leq r\leq1}|\log(\bar{\mathbb{P}}^r_u/\break \bar{\mathbb
{P}}^0_0)| (Y_{\Pi})\}_{n\geq n_u}$  is tight.
\end{proposition}

\begin{pf}\label{euler-app}
 We easily obtain the result by Lemma~\ref
{pre-nonsynchro-tight-lemma} and an estimate
\[
\sup_r \biggl|\log\frac{\bar{\mathbb{P}}^r_{M,u}}{\bar{\mathbb
{P}}^0_{M,0}} \biggr|\leq\sup
_r \biggl|\log\frac{1}{\bar{\mathbb
{P}}^0_{M,0}}\int_{\bar{A}^n_M}
\frac{\mathbb{P}^r_u}{\mathbb
{P}^0_0}\mathbb{P}^0_0\,\mathrm{d}\hat{z} \biggr|\leq M
\]
for sufficiently large $M>0$.
\end{pf}

The following lemma is similarly proved and used later.

%
\begin{lemma}\label{Pvu-tightness}
Let $u\in\mathbb{R}^d$. Assume $[A1']$ and $[A3']$.
Then $\{\sup_{0\leq v\leq1}|\log(\mathbb{P}^0_{vu}/\break\mathbb
{P}^0_0)| (Y_{\check{U}})\}_{n\geq n_u}$ and
$\{\sup_{0\leq v\leq1}|\log(\bar{\mathbb{P}}^0_{vu}/\bar{\mathbb
{P}}^0_0)|(Y_{\Pi})\}_{n\geq n_u}$  are  tight.
\end{lemma}


\section{The proof of LAMN property}\label{LAMN-proof-section}

In this section, we will complete the proof of the LAMN property of $\{
P_{\sigma,n}\}_{\sigma,n}$.

It is essential in the proof to replace $\bar{\mathbb{P}}^0_u$ in (\ref
{ratio-equ1}) by the function $\int\exp(\sum_k\tilde{f}_k^u)\,\mathrm{d}\hat{z}$
below so that coefficient $b$ is predictable and does not depend on
$\hat{z}$.
For this purpose, we use It\^{o}'s rule and martingale
properties and estimate the difference.
However, the proof is technically complicated because the function $\log
\bar{\mathbb{P}}^0_u$ contains a $\mathrm{d}\hat{z}$-integral of an exponential function.
This integral is far more difficult to handle than a simple function of
increments of the process, which appears in  synchronous
sampling models of Gobet \cite{gob}.
We estimate the difference step by step in Lemmas \ref
{coeff-move-lemma1} and \ref{coeff-move-lemma2}.
The function $\log\int\exp(\sum_k\tilde{f}_k^u)\,\mathrm{d}\hat{z}$ can be
rewritten in a simple function of increments of the process as seen in
Lemma~\ref{induction-lemma}.
Then the proof is completed by  proving asymptotic equivalence
of the replaced likelihood ratio and the quasi-likelihood ratio
$H_n(\sigma)-H_n(\sigma_{\ast})$.

In the following, we assume that $[A2]$ holds true.
Let $\mu_{(k)}(z,\sigma)=\mu(\check{u}^{k-1},z_{k-1},\sigma)$,
$b_{(k)}(z,\sigma)=b(\check{u}^{k-1},z_{k-1},\sigma)$, $\tilde
{b}_{(k)}(z,\sigma)=b(\check{u}^{k-1},x_{k_1(i(k))},y_{k_2(j(k))},\sigma)$,
\begin{eqnarray*}
f_{(k)}(z,\sigma)&=&-\tfrac{1}{2}\bigl(z_k-z_{k-1}-
\Delta\check{u}^k \mu _{(k)}(z,\sigma)\bigr)^{\star}
\bigl(\Delta\check{u}^kb_{(k)}b^{\star
}_{(k)}(z,
\sigma)\bigr)^{-1}\\
&&{}\times\bigl(z_k-z_{k-1}-\Delta
\check{u}^k \mu _{(k)}(z,\sigma)\bigr)
\\
&&{}-\tfrac{1}{2}\log\det\bigl(\Delta\check{u}^k
b_{(k)}b^{\star
}_{(k)}(z,\sigma)\bigr)-\log(2\uppi),
\\
\tilde{f}_{(k)}(z,\sigma)&=&-\tfrac{1}{2}(z_k-z_{k-1})^{\star}
\bigl(\Delta \check{u}^k\tilde{b}_{(k)}\tilde{b}^{\star}_{(k)}(z,
\sigma )\bigr)^{-1}(z_k-z_{k-1})\\
&&{}-\tfrac{1}{2}
\log\det\bigl(\Delta\check{u}^k\tilde {b}_{(k)}
\tilde{b}^{\star}_{(k)}(z,\sigma)\bigr)-\log(2\uppi)
\end{eqnarray*}
for $z=(z_k)_{k=0}^{L_0(\bar{u})}=((x_k)_{k=0}^{L_0(\bar
{u})},(y_k)_{k=0}^{L_0(\bar{u})})\in\mathbb{R}^{2L_0(\bar{u})+2}$, and
let $\mu^{u}_k(z)=\mu_{(k)}(z,\sigma^n_u)$, $b^{u}_k(z)=b_{(k)}(z,\sigma
^n_u)$, $\tilde{b}^{u}_k(z)=\tilde{b}_{(k)}(z,\sigma^n_u)$,
$f^u_k(z)=f_{(k)}(z,\sigma^n_u)$, $\tilde{f}^u_k(z)=\tilde
{f}_{(k)}(z,\sigma^n_u)$,
$f^{u,(1)}_k(z)= \partial_{\sigma}f_{(k)}(z,\sigma^n_u)$, $\tilde
{f}^{u,(1)}_k(z)=\partial_{\sigma}\tilde{f}_{(k)}(z,\sigma^n_u)$.

Then we obtain $\bar{\mathbb{P}}^1_u=\int\exp(\sum_kf^u_k(z))\,\mathrm{d}\hat{z}$.

Moreover, let $\kappa$ be a positive constant satisfying
\[
\delta_2\vee(\delta_1+\delta_3)<\kappa<
\biggl(\frac
{1}{4}-\frac{(3\delta_1+2\delta_3)\vee(\delta_1+\delta_2)}{2} \biggr)\wedge \biggl(
\frac{1}{6}-\frac{\delta_1}{2} \biggr),
\]
$h=h_n=[b_n^{\kappa}]$ and
\[
\acute{f}^{k,u}_{k'}(z)= \left\{
\begin{array}{l@{\qquad}l} \tilde{f}^u_{k'}(z), &
\bigl|k-k'\bigr|\leq h,
\\
\log\check{p}^0_{k',u}(z) & \mbox{otherwise},
\end{array}\right.
\]
where $\{\delta_j\}_{j=1}^3$  appears in $[A2]$.
Then we obtain
%
\begin{eqnarray}
\log\frac{\bar{\mathbb{P}}^0_u}{\bar{\mathbb{P}}^0_0}(Y_{\Pi})& =& \int
^1_0\frac{\partial_{v}(\bar{\mathbb{P}}^0_{vu})}{\bar{\mathbb
{P}}^0_{vu}}\,\mathrm{d}v(Y_{\Pi})
\nonumber
\\[-8pt]
\label{P1-equ}\\[-8pt]
\nonumber
&=& b_n^{-{1}/{2}}u\int^1_0
\frac{\int\sum_k(\check
{p}^{0,(1)}_{k,vu}/\check{p}^0_{k,vu})\exp(\sum_{k'}\log\check
{p}^0_{k',vu})(z)\,\mathrm{d}\hat{z}}{\int\exp(\sum_{k'}\log\check
{p}^0_{k',vu})(z)\,\mathrm{d}\hat{z}}\,\mathrm{d}v(Y_{\Pi}).
\end{eqnarray}

If we have asymptotic equivalence of $\log(\bar{\mathbb{P}}^0_u/\bar
{\mathbb{P}}^0_0)(Y_{\Pi})$ and
%
\begin{eqnarray}
&& \log \biggl(\int\exp \biggl(\sum_k
\tilde{f}^u_k \biggr)\,\mathrm{d}\hat{z}\Big/\int\exp \biggl(\sum
_k\tilde{f}^0_k \biggr)\,\mathrm{d}
\hat{z} \biggr) (Y_{\Pi})
\nonumber
\\[-8pt]
\label{likelihood-ratio1}\\[-8pt]
\nonumber
&& \quad =b_n^{-1/2}u\int
^1_0\frac{\int\sum_k\tilde{f}^{vu,(1)}_k\exp(\sum_{k'}\tilde{f}^{vu}_{k'})\,\mathrm{d}\hat{z}}{\int\exp(\sum_{k'}\tilde
{f}^{vu}_{k'})\,\mathrm{d}\hat{z}}\,\mathrm{d}v(Y_{\Pi}),
\end{eqnarray}
then Lemma~\ref{induction-lemma} gives a simple  asymptotic
representation of $\log(\bar{\mathbb{P}}^0_u/\bar{\mathbb
{P}}^0_0)(Y_{\Pi})$ using the increments of processes.
However, it is difficult to estimate directly the difference of these
two quantities since
$\exp(\sum_{k'}\log\check{p}^{0}_{k',vu})-\exp(\sum_{k'}\tilde
{f}^{vu}_{k'})$ is not asymptotically negligible.
So we first prove asymptotic equivalence of $\log(\bar{\mathbb
{P}}^0_u/\bar{\mathbb{P}}^0_0)(Y_{\Pi})$ and
%
\begin{equation}
\label{likelihood-ratio2}
b_n^{-1/2}u\int^1_0
\frac{\int\sum_k\tilde{f}^{vu,(1)}_k\exp(\sum_{k'}\acute{f}^{k,vu}_{k'})\,\mathrm{d}\hat{z}}{\int\exp(\sum_{k'}{\acute{f}}^{k,vu}_{k'})\,\mathrm{d}\hat{z}}\,\mathrm{d}v(Y_{\Pi})
\end{equation}
in Lemmas \ref{coeff-move-lemma1} and \ref{coeff-move-lemma2}.
Then we prove asymptotic equivalence of (\ref{likelihood-ratio1}) and
(\ref{likelihood-ratio2}) in Lemma~\ref{integration-lemma},
using a simpler expression of (\ref{likelihood-ratio2}) obtained by
calculating $\mathrm{d}\hat{z}$-integral partially by the virtue of Lemma~\ref
{induction-lemma}.

We start with preparation of several lemmas.
The first one is proved similarly to Lemma~5 in Ogihara and Yoshida
\cite{ogi-yos}, so we omit details.

%
\begin{lemma}\label{doubleThetaSum-est}
Assume $[A2]$ and $[A3']$. Then
\[
b_n^{-{1}/{2}+\delta}\sum_{p_1,p_2=0}^{\infty}
\frac{\sum_{l_1,l_2}|\theta_{p_1,l_1}\cap\theta
_{p_2,l_2}|}{(p_1+1)^5(p_2+1)^5}\to^p 0
\]
as $n\to\infty$ for any $\delta$ satisfying $0<\delta<1/2-(3\delta
_1+2\delta_3)\vee(\delta_1+\delta_2)$.
\end{lemma}

\begin{lemma}\label{cutoff-lemma}
Let $u\in\mathbb{R}^d$. Assume $[A1'],[A2]$ and $[A3']$.
Then for any $\varepsilon,\eta>0$, there exists $M'>0$ such that
\[
\sup_{n\geq n_u}P \biggl[\sup_{0\leq v\leq1} \biggl|\log
\frac{\bar{\mathbb
{P}}^0_{vu}}{\tilde{\mathbb{P}}^0_{M,vu}} \biggr|(Y_{\Pi})\geq\eta \biggr]<\varepsilon
\]
for $M\geq M'$, where
\[
B^n_M= \Bigl\{z\in\mathbb{R}^{2L_0+2};\sup
_{k',v'}\bigl|\tilde {f}^{v'u}_{k'}-\log
\check{p}^0_{k',v'u}(z_{k'-1},z_{k'})\bigr|\leq
Mb_n^{-1/3-\kappa} \Bigr\},
\]
$\bar{B}^n_M(z_0,\bar{z})=\{\hat{z};(z_0,\bar{z},\hat{z})\in B^n_M\}$
and $\tilde{\mathbb{P}}^0_{M,u}(z_0,\bar{z})=\int_{\bar{B}^n_M(z_0,\bar
{z})}\mathbb{P}^0_u(z)\,\mathrm{d}\hat{z}$.
\end{lemma}

\begin{pf}
We will apply 2. of Proposition~\ref{synch-to-nonsynch}.
By using the Burkholder--Davis--Gundy inequality and Lemma~\ref{density-est2}, we have
%
\begin{eqnarray}
&& E\biggl[\sup_v\biggl|\log\frac{\mathbb{P}^0_{vu}}{\mathbb{P}^0_0}
\biggr|(Y_{\check{U}})\Big|\Pi\biggr]
\nonumber
\\
&& \quad =E \biggl[\sup_v \biggl|\int^v_0
\frac{\partial_s(\mathbb
{P}^0_{su})}{\mathbb{P}^0_{su}}\,\mathrm{d}s \bigg|(Y_{\check{U}}) \Big|\Pi \biggr]\nonumber \\
&& \label{synchro-k-est}\quad \leq
b_n^{-{1}/{2}}|u|\int^1_0E
\biggl[ \biggl|\sum_k\frac{\check
{p}^{0,(1)}_{k,vu}}{\check{p}^0_{k,vu}}(Y_{\check{U}^{k-1}},Y_{\check
{U}^k})
\biggr| \Big|\Pi \biggr]\,\mathrm{d}v
\\
\nonumber
&&\quad \leq b_n^{-{1}/{2}}|u|E \biggl[ \biggl|\sum
_k\frac{\check
{p}^{0,(1)}_{k,0}}{\check{p}^0_{k,0}}(Y_{\check{U}^{k-1}},Y_{\check
{U}^k}) \biggr|
\Big|\Pi \biggr] \\
&&\qquad {}+b_n^{-{1}/{2}}|u|\int^1_0\!\!
\int^v_0E \biggl[ \biggl|\sum
_k\partial_v \biggl(\frac{\check{p}^{0,(1)}_{k,v_2u}}{\check{p}^0_{k,v_2u}} \biggr)
(Y_{\check
{U}^{k-1}},Y_{\check{U}^k}) \biggr| \Big|\Pi \biggr]\,\mathrm{d}v_2\,\mathrm{d}v
\nonumber\\
&& \quad \leq Cb_n^{-1/2}|u|(\ell_{1,n}+
\ell_{2,n})^{1/2}+Cb_n^{-1}|u|^2(
\ell _{1,n}+\ell_{2,n}).\nonumber
\end{eqnarray}

On the other hand, for any $\varepsilon>0$, Lemmas \ref{density-est} and
\ref{density-est2} yield
%
\begin{eqnarray}
&&\sup_vP\bigl[Y^{0,vu}_{\check{U}}
\in\bigl(B^n_M\bigr)^c|\Pi\bigr]
\nonumber
\\
&& \quad \leq \frac{b_n^{q/3+q\kappa}}{M^q}\sup_vE \biggl[\sup
_{k',v'} \biggl|\tilde{f}^{v'u}_{k'}-f^{v'u}_{k'}+
\int^1_0\frac{\partial_r\check
{p}^r_{k',v'u}}{\check{p}^r_{k',v'u}}\,\mathrm{d}r \biggr|^q
\bigl(Y^{0,vu}_{\check
{U}}\bigr) \Big|\Pi \biggr]
\nonumber
\\
\label{BnMc-est}
&&\quad \leq C_q\frac{b_n^{q/3+q\kappa}}{M^q} \biggl(r_n^{{q}/{2}}(
\ell _{1,n}+\ell_{2,n})\\
&&\hspace*{54pt}\qquad {}+\sup_{r,v}E
\biggl[\sum_{k'} \biggl( \biggl|\frac{\partial_r\check
{p}^r_{k',0}}{\check{p}^r_{k',0}}
\biggr|^q +\int^1_0 \biggl|
\partial_v \biggl(\frac{\partial_r\check
{p}^r_{k',v'u}}{\check{p}^r_{k',v'u}} \biggr) \biggr|^q\,\mathrm{d}v'
\biggr) \bigl(Y^{0,vu}_{\check{U}}\bigr) \Big|\Pi \biggr] \biggr)\nonumber
\\
&& \quad \leq C_qb_n^{q/3+q\kappa}M^{-q}r_n^{q/2}(
\ell_{1,n}+\ell_{2,n})\nonumber
\end{eqnarray}
for any $q>0$ and $M>0$.

By (\ref{synchro-k-est}), (\ref{BnMc-est}) and 2. of
Proposition~\ref{synch-to-nonsynch}, we obtain
\begin{eqnarray*}
&& P \biggl[\sup_v \biggl|\log\frac{\bar{\mathbb{P}}^0_{vu}}{\tilde{\mathbb
{P}}^0_{M,vu}}
\biggr|(Y_{\Pi})\geq\eta \Big|\Pi \biggr] \\
&& \quad \leq \frac{C_q\mathrm{e}^{M'}}{1-\mathrm{e}^{-\eta}} \biggl\{
\biggl(\frac{1}{M^q}+\frac
{M_2}{M^q} \biggr)b_n^{{q}/{3}+q\kappa}r_n^{{q}/{2}}(
\ell _{1,n}+\ell_{2,n})+\frac{1}{M_2}\sup
_vE \biggl[ \biggl|\frac{\partial
_v\mathbb{P}^0_{vu}}{\mathbb{P}^0_{vu}} \biggr|^2 \Big|\Pi
\biggr] \biggr\}
\\
&&\qquad {}+\frac{C_q}{M'}\bigl(1+b_n^{-1}|u|^2(
\ell_{1,n}+\ell_{2,n})\bigr)
\end{eqnarray*}
for any $M,M',M_2>0$.

Hence, we have
\[
\sup_nP \biggl[\sup_v \biggl|\log
\frac{\bar{\mathbb{P}}^0_{vu}}{\tilde
{\mathbb{P}}^0_{M,vu}} \biggr|(Y_{\Pi})\geq\eta \biggr] = \sup
_nE \biggl[P \biggl[\sup_v \biggl|\log
\frac{\bar{\mathbb
{P}}^0_{vu}}{\tilde{\mathbb{P}}^0_{M,vu}} \biggr|(Y_{\Pi})\geq\eta \Big|\Pi \biggr]\wedge1 \biggr]<
\varepsilon
\]
for sufficiently large $M>0$.
\end{pf}

Similarly to (\ref{P1-equ}), we obtain
%
\begin{eqnarray}
\log\frac{\tilde{\mathbb{P}}^0_{M,u}}{\tilde{\mathbb{P}}^0_{M,0}}(Y_{\Pi
})& =& \int
^1_0\frac{\partial_v(\tilde{\mathbb{P}}^0_{M,vu})}{\tilde{\mathbb
{P}}^0_{M,vu}}\,\mathrm{d}v(Y_{\Pi})
\nonumber
\\[-8pt]
\label{tildeP1-equ}\\[-8pt]
\nonumber
&=& b_n^{-{1}/{2}}u\int^1_0
\frac{\int_{\bar{B}^n_M} \sum_k(\check
{p}^{0,(1)}_{k,vu}/\check{p}^0_{k,vu})\exp(\sum_{k'}\log\check
{p}^0_{k',vu})(z)\,\mathrm{d}\hat{z}}{\tilde{\mathbb{P}}^0_{M,vu}}\,\mathrm{d}v(Y_{\Pi}).
\end{eqnarray}

Let $\bar{\mathbb{P}}^{2,u}_k(g)(z_0,\bar{z})=\int g(z)\exp(\sum_{k'}\acute{f}_{k'}^{k,u}(z))\,\mathrm{d}\hat{z}$ for an integrable function $g$.
For $1\leq k\leq L_0(\bar{u})$ and $p\in\mathbb{Z}_+$, let $\tilde
{\theta}(p,k;\bar{u})$ be $\theta(p,l;\bar{u})$,
where an\vspace*{1pt} integer $l$ satisfies $1\leq l\leq L^1$ and $[\check
{u}^{k-1},\check{u}^k)\subset[s^{i-1},s^i)$.
Let $\tilde{\theta}_{p,k}=\tilde{\theta}(p,k;\Pi)$.

The following lemma is the first step to replace $\bar{\mathbb{P}}^0_u$
by $\int\exp(\sum_k\tilde{f}^u_k)\,\mathrm{d}\hat{z}$.

%
\begin{lemma}\label{coeff-move-lemma1}
Let $u\in\mathbb{R}^d$. Assume $[A1'],[A2]$ and $[A3']$. Then for any
$\varepsilon,\eta>0$, there exist $M'>0$ and $\{N_M\}_{M\geq M'}\subset
\mathbb{N}$ such that
\[
P \biggl[ \biggl|\log\frac{\bar{\mathbb{P}}^0_u}{\bar{\mathbb
{P}}^0_0}(Y_{\Pi})-b_n^{-{1}/{2}}u
\int^1_0\sum_k
\frac{\bar{\mathbb
{P}}^{2,vu}_k(\tilde{f}^{vu,(1)}_k1_{B^n_M})}{\tilde{\mathbb
{P}}^{0}_{M,vu}}\,\mathrm{d}v(Y_{\Pi}) \biggr|\geq\eta \biggr]<\varepsilon
\]
for $M\geq M'$ and $n\geq N_M$.
\end{lemma}

\begin{pf}
Fix $\varepsilon,\eta\in(0,1)$.
By Lemmas \ref{Pvu-tightness} and \ref{cutoff-lemma}, there exists
$M'>0$ such that $\sup_{n\geq n_u}P[Y_{\Pi}\in(\mathcal{K}^1_M)^c(\Pi
)]<\varepsilon/2$ for $M\geq M'$,
where
\[
\mathcal{K}^1_M(\bar{u})=\Bigl\{(z_0,
\bar{z});\sup_{0\leq v\leq1}\bigl|\log\bigl(\bar {\mathbb{P}}^0_{vu}/
\bar{\mathbb{P}}^0_0\bigr)\bigr|(z_0,\bar{z})\leq M
\mbox{ and }  \sup_{0\leq v\leq1} \bigl|\log\bigl(\bar{
\mathbb{P}}^0_{vu}/\tilde {\mathbb{P}}^0_{M,vu}
\bigr)\bigr|(z_0,\bar{z})\leq1 \Bigr\}.
\]

Therefore by (\ref{tildeP1-equ}), Lemmas \ref{cond-conv} and \ref
{cutoff-lemma}, it is sufficient to show that
\begin{eqnarray*}
\Phi_n &=& E \biggl[ \biggl|b_n^{-{1}/{2}}u
\int
^1_0\sum_k
\biggl(\int_{\bar
{B}^n_M}\biggl\{ \bigl(\check{p}^{0,(1)}_{k,vu}/\check{p}^0_{k,vu}\bigr)\exp\biggl(\sum_{k'}\log\check{p}^0_{k',vu}\biggr)(z)\\
&&\hspace*{104pt}{}-\tilde{f}^{vu,(1)}_k
\exp\biggl(\sum_{k'}\acute{f}^{k,vu}_{k'}\biggr)\biggr\}\,\mathrm{d}\hat{z}\biggr)\Big/\bigl(
\tilde{\mathbb{P}}^0_{M,vu}\bigr)\,\mathrm{d}v \biggr|\\
&&\hspace*{16pt} {}\times 1_{\mathcal{K}^1_M(\Pi)}(Y_{\Pi
}) \Big|\Pi \biggr]
\\
&\to^p & 0
\end{eqnarray*}
as $n\to\infty$ for any $M>0$.

By the definition of $\mathcal{K}^1_M$ and the relation $|\exp
(x)-1-x|\leq Cx^2$ for $|x|\leq3M$, we obtain
%
\begin{eqnarray}
\Phi_n&\leq&\mathrm{e}^{M+1}|u|\nonumber\\
&&{}\times\sup
_vE \biggl[ \biggl|b_n^{-{1}/{2}}\sum
_k \biggl\{\frac{\check{p}^{0,(1)}_{k,vu}}{\check{p}^0_{k,vu}}
-\tilde {f}^{vu,(1)}_k
\exp \biggl(\sum_{k';|k'-k|\leq h}\bigl(\tilde{f}^{vu}_{k'}-
\log \check{p}^0_{k',vu}\bigr) \biggr) \biggr\}
\biggr|\nonumber\\
&&\hspace*{22pt}\qquad{}\times 1_{B^n_M}\bigl(Y^{0,vu}_{\check
{U}}\bigr) \Big|\Pi \biggr]
\nonumber
\\
& \leq &  C\sup_vE \biggl[ \biggl|b_n^{-{1}/{2}}
\sum_k \biggl\{\frac{\check
{p}^{0,(1)}_{k,vu}}{\check{p}^0_{k,vu}}-
\tilde{f}^{vu,(1)}_k \biggl(1+\sum
_{k';|k'-k|\leq h}\bigl(\tilde{f}^{vu}_{k'}-\log
\check{p}^0_{k',vu}\bigr) \biggr) \biggr\} \biggr|
\bigl(Y^{0,vu}_{\check{U}}\bigr) \Big|\Pi \biggr]\nonumber\\
&&{}+\mathrm{o}_p(1)
\nonumber
\\[-8pt]
\label{Phi-ineq1}\\[-8pt]
\nonumber
&\leq&C\sup_vE \biggl[ \biggl|b_n^{-{1}/{2}}
\sum_k \biggl\{\frac{\check
{p}^{0,(1)}_{k,vu}}{\check{p}^0_{k,vu}}-
\frac{\check
{p}^{1,(1)}_{k,vu}}{\check{p}^1_{k,vu}} \biggr\} \biggr|\bigl(Y^{0,vu}_{\check
{U}}\bigr) \Big|\Pi
\biggr]
\\[-2pt]
&& {}+C\sup_vE \biggl[ \biggl|b_n^{-{1}/{2}}\sum
_k\tilde{f}^{vu,(1)}_k\sum
_{k';|k'-k|\leq h}\bigl(\log\check{p}^1_{k',vu}-
\log\check {p}^0_{k',vu}\bigr) \biggr|\bigl(Y^{0,vu}_{\check{U}}
\bigr) \Big|\Pi \biggr]\nonumber
\\[-2pt]
&& {}+C\sup_vE \biggl[ \biggl|b_n^{-{1}/{2}}\sum
_k \biggl\{ f^{vu,(1)}_k-
\tilde{f}^{vu,(1)}_k \biggl(1+\sum
_{k';|k'-k|\leq h}\bigl(\tilde {f}^{vu}_{k'}-f^{vu}_{k'}
\bigr) \biggr) \biggr\} \biggr|\bigl(Y^{0,vu}_{\check{U}}\bigr) \Big|\Pi
\biggr]\nonumber\\[-2pt]
&&{}+\mathrm{o}_p(1)
\nonumber
\\[-2pt]
&=&\Phi_{n,1}+\Phi_{n,2}+\Phi_{n,3}+\mathrm{o}_p(1).\nonumber
\end{eqnarray}

The quantity $\Phi_{n,1}$ is\vspace*{-2pt} estimated as
\begin{eqnarray*}
\Phi_{n,1}&\leq&C\sup_{r,v}E \biggl[
\biggl|b_n^{-{1}/{2}}\sum_k
\partial_r \biggl(\frac{\check{p}^{r,(1)}_{k,vu}}{\check
{p}^r_{k,vu}} \biggr) \biggr|\bigl(Y^{0,vu}_{\check{U}}
\bigr) \Big|\Pi \biggr]
\\[-2pt]
&\leq&Cb_n^{-{1}/{2}}\sum_k\sup
_{r,v}E \biggl[ \biggl|E \biggl[\partial _r \biggl(
\frac{\check{p}^{r,(1)}_{k,vu}}{\check{p}^r_{k,vu}} \biggr) \bigl(Y^{0,vu}_{\check{u}^{k-1}},Y^{0,vu}_{\check{u}^k}
\bigr) \Big|\mathcal {F}_{\check{u}^{k-1}} \biggr] \biggr| \biggr] \bigg|_{\bar{u}=\Pi}
\\[-2pt]
&&{}+C\sup_{r,v}E \biggl[ \biggl|b_n^{-{1}/{2}}\sum
_k \biggl(\partial_r \biggl(
\frac{\check{p}^{r,(1)}_{k,vu}}{\check{p}^r_{k,vu}} \biggr) \bigl(Y^{0,vu}_{\check{u}^{k-1}},Y^{0,vu}_{\check{u}^k}
\bigr) \\[2pt]
&&\hspace*{97pt}{}-E \biggl[\partial_r \biggl(\frac{\check{p}^{r,(1)}_{k,vu}}{\check
{p}^r_{k,vu}} \biggr)
\bigl(Y^{0,vu}_{\check{u}^{k-1}},Y^{0,vu}_{\check
{u}^k}\bigr) \Big|
\mathcal{F}_{\check{u}^{k-1}} \biggr] \biggr) \biggr| \biggr] \bigg|_{\bar{u}=\Pi}
\\[-2pt]
&\leq& Cb_n^{-{1}/{2}}\sum_k
\sup_{r,v}E \biggl[ \biggl|\int\partial _r \biggl(
\frac{\check{p}^{r,(1)}_{k,vu}}{\check{p}^r_{k,vu}} \biggr)\check {p}^0_{k,vu}(z_{k-1},z_k)\,\mathrm{d}z_k
\Big|_{z_{k-1}=Y^{0,vu}_{\check
{u}^{k-1}}} \biggr| \biggr] \bigg|_{\bar{u}=\Pi} + \mathrm{o}_p(1).
\end{eqnarray*}
Then we have $\Phi_{n,1}=\mathrm{o}_p(1)$ since
\begin{eqnarray*}
&& \biggl|\int\partial_r \biggl(\frac{\check{p}^{r,(1)}_{k,vu}}{\check
{p}^r_{k,vu}} \biggr)
\check{p}^0_{k,vu}(z_{k-1},z_k)\,\mathrm{d}z_k
\biggr|
\\[-2pt]
&&\quad = \biggl|\int\partial_r \biggl(\frac{\check{p}^{r,(1)}_{k,vu}}{\check
{p}^r_{k,vu}} \biggr) \biggl(
\check{p}^r_{k,vu}-\int^r_0
\partial_r\check {p}^{r'}_{k,vu}\,\mathrm{d}r'
\biggr) (z_{k-1},z_k)\,\mathrm{d}z_k \biggr|
\\[-2pt]
&& \quad \leq  \biggl|\int\partial_r \biggl(\frac{\check{p}^{r,(1)}_{k,vu}}{\check
{p}^r_{k,vu}} \biggr)
\check{p}^r_{k,vu}(z_{k-1},z_k)\,\mathrm{d}z_k
\biggr|\\[-2pt]
&& \qquad  {}+ \biggl\{\int \biggl(\partial_r \biggl(\frac{\check
{p}^{r,(1)}_{k,vu}}{\check{p}^r_{k,vu}}
\biggr) \biggr)^2\check {p}^r_{k,vu}\,\mathrm{d}z_k
\biggr\}^{{1}/{2}} \sup_{r'} \biggl(\int \biggl(
\frac{\partial_r\check{p}^{r'}_{k,vu}}{\check
{p}^r_{k,vu}} \biggr)^2\check{p}^r_{k,vu}\,\mathrm{d}z_k
\biggr)^{{1}/{2}},
\end{eqnarray*}
$E[|\check{p}^{r,(1)}_{k,vu}/\check{p}^r_{k,vu}-\check
{p}^{1,(1)}_{k,vu}/\check{p}^1_{k,vu}|^p(Y^{r,vu}_{\check
{u}})]^{1/p}=\mathrm{O}((\Delta\check{u}^k)^{1/2})$ and
\begin{eqnarray*}
&&\int\frac{\partial_r\check{p}^r_{k,vu}}{\check{p}^r_{k,vu}}\frac
{\check{p}^{1,(1)}_{k,vu}}{\check{p}^1_{k,vu}}\check {p}^r_{k,vu}(z_{k-1},z_k)\,\mathrm{d}z_k
\\
&& \quad =-\frac{1}{2}E_{z_{k-1}} \biggl[\frac{\delta(\mathcal{B}^r\partial
_r\mathbf{Y}^{r,u,k,z_{k-1}}_{\Delta\check{u}^k})}{\Delta\check{u}^k}
\\
&&\qquad \hspace*{42pt}{}\times\partial_{\sigma} \biggl(\bigl(\Delta\bar{\mathbf{Y}}^k
\bigr)^{\star}\frac
{(bb^{\star})^{-1}(\check{u}^{k-1},z_{k-1},\sigma^n_{vu})}{\Delta\check
{u}^k}\Delta\bar{\mathbf{Y}}^k\\
&&\hspace*{70pt}\qquad {}+\log
\det\bigl(bb^{\star}\bigr) \bigl(\check {u}^{k-1},z_{k-1},
\sigma^n_{vu}\bigr) \biggr) \biggr]
\\
&& \quad =\mathrm{O}\bigl(\Delta\check{u}^k\bigr),
\end{eqnarray*}
where $\Delta\bar{\mathbf{Y}}^k=(\mathbf{Y}^r_{\Delta\check
{u}^k}-z_{k-1}-\Delta\check{u}^k\mu(\check{u}^{k-1},z_{k-1},\sigma^n_{vu}))$.

Similarly, $\Phi_{n,2}$ is estimated as
\begin{eqnarray*}
\Phi_{n,2}&\leq&C\sup_{r,v}E \biggl[
\biggl|b_n^{-{1}/{2}}\sum_k\tilde
{f}^{vu,(1)}_k\sum_{|k-k'|\leq h}
\frac{\partial_r\check
{p}^r_{k',vu}}{\check{p}^r_{k',vu}} \biggr|\bigl(Y^{0,vu}_{\check{U}}\bigr) \Big|\Pi \biggr]
\\
&\leq & C\sup_{r,v}E \biggl[ \biggl|b_n^{-{1}/{2}}
\sum_k\bigl(\tilde {f}^{vu,(1)}_k-E
\bigl[\tilde{f}^{vu,(1)}_k|\mathcal{F}_{\check
{u}^{k-1}}\bigr]
\bigr)\sum_{|k-k'|\leq h}\frac{\partial_r\check
{p}^r_{k',vu}}{\check{p}^r_{k',vu}} \biggr|
\bigl(Y^{0,vu}_{\check{u}}\bigr) \biggr] \bigg|_{\bar{u}=\Pi}\\
&&{}+\mathrm{o}_p
\bigl(b_n^{-{1}/{2}}b_nr_nb_n^{\kappa}
\bigr)
\\
&\leq& C\sup_{r,v}E \biggl[ \biggl|b_n^{-{1}/{2}}
\sum_k\bigl(\tilde {f}^{vu,(1)}_k-E
\bigl[\tilde{f}^{vu,(1)}_k|\mathcal{F}_{\check
{u}^{k-1}}\bigr]
\bigr)\sum_{k'; k-h\leq k'\leq k}\frac{\partial_r\check
{p}^r_{k',vu}}{\check{p}^r_{k',vu}} \biggr|
\bigl(Y^{0,vu}_{\check{u}}\bigr) \biggr] \bigg|_{\bar{u}=\Pi}
\\
&&{}+C\sup_{r,v}E \biggl[ \biggl|b_n^{-{1}/{2}}\sum
_{k'} \biggl(\sum_{k;k'-h\leq k<k'}
\bigl(\tilde{f}^{vu,(1)}_k-E\bigl[\tilde{f}^{vu,(1)}_k|
\mathcal {F}_{\check{u}^{k-1}}\bigr]\bigr) \biggr)\\
&&\hspace*{51pt}{}\times\frac{\partial_r\check{p}^r_{k',vu}}{\check
{p}^r_{k',vu}} \biggr|
\bigl(Y^{0,vu}_{\check{u}}\bigr) \biggr] \bigg|_{\bar{u}=\Pi
}\\
&&{}+\mathrm{o}_p(1)
\\
&=&\mathrm{o}_p(1),
\end{eqnarray*}
by the Burkholder--Davis--Gundy inequality.

Finally, we will prove $\Phi_{n,3}=\mathrm{o}_p(1)$.
Let $t^0_k=S^{n,i(k;\Pi)}\wedge T^{n,j(k;\Pi)}$, $\hat
{b}^v_k=b(t^0_k,Y^{0,vu}_{t^0_k},\sigma^n_{vu})$ and
\begin{eqnarray*}
\mathcal{A}^1_k&=&-\Delta W_k^{\star}
\bigl(\hat{b}^v_k\bigr)^{\star}\bigl(\partial
_{\sigma}\bigl(\bigl(b_{(k)}b^{\star}_{(k)}
\bigr)^{-1}\bigr)-\partial_{\sigma}\bigl(\bigl(\tilde
{b}_{(k)}\tilde{b}^{\star}_{(k)}\bigr)^{-1}
\bigr)\bigr) \bigl(Y^{0,vu}_{\check{U}},\sigma ^n_{vu}
\bigr)\hat{b}^v_k\Delta W_k/\bigl(2\Delta
\check{U}^{k}\bigr)
\\
&&{}-\frac{1}{2}\partial_{\sigma}\log\frac{\det(b_{(k)}b^{\star
}_{(k)})}{\det(\tilde{b}_{(k)}\tilde{b}^{\star
}_{(k)})}
\bigl(Y^{0,vu}_{\check{U}},\sigma^n_{vu}
\bigr) +\partial_{\sigma}\bigl(\mu^{\star}_{(k)}
\bigl(b_{(k)}b^{\star
}_{(k)}\bigr)^{-1}\bigr)
\bigl(Y^{0,vu}_{\check{U}},\sigma^n_{vu}
\bigr)\hat{b}^v_k\Delta W_k,
\\
\mathcal{A}^2_k&=&-\Delta W_k^{\star}
\bigl(\hat{b}^v_k\bigr)^{\star}\partial
_{\sigma}\bigl(\bigl(\hat{b}^v_k\bigl(
\hat{b}^v_k\bigr)^{\star}\bigr)^{-1}
\bigr)\hat{b}^v_k\Delta W_k/\bigl(2\Delta
\check{U}^{k}\bigr)-\frac{1}{2}\partial_{\sigma}\log\det
\bigl(\hat {b}^v_k\bigl(\hat{b}^v_k
\bigr)^{\star}\bigr),
\\
\mathcal{A}^3_{k,k'}&=& \biggl\{-\frac{1}{2}\operatorname{tr}
\bigl(\bigl(\hat {b}^v_{k'}\bigr)^{\star}\bigl(
\bigl(\tilde{b}^{vu}_{k'}\bigl(\tilde{b}^{vu}_{k'}
\bigr)^{\star
}\bigr)^{-1}-\bigl(b^{vu}_{k'}
\bigl(b^{vu}_{k'}\bigr)^{\star}\bigr)^{-1}
\bigr) \bigl(Y^{0,vu}_{\check
{U}}\bigr)\hat{b}^v_{k'}
\bigr)\\
&&\hspace*{5pt}{}-\frac{1}{2}\log\frac{\det(\tilde
{b}_{k'}^{vu}(\tilde{b}_{k'}^{vu})^{\star})}{\det
(b_{k'}^{vu}(b_{k'}^{vu})^{\star})}\bigl(Y^{0,vu}_{\check{U}}
\bigr) \biggr\}1_{\{
|k-k'|\leq h\} },
\\
\mathcal{A}^4_{k,k'}&=&-\int_{\check{u}^{k'-1}}^{\check
{u}^{k'}}(W_t-W_{\check{u}^{k'-1}})^{\star}
\bigl(\hat{b}^v_{k'}\bigr)^{\star}
\frac
{(\tilde{b}^{vu}_{k'}(\tilde{b}^{vu}_{k'})^{\star
})^{-1}-(b^{vu}_{k'}(b^{vu}_{k'})^{\star})^{-1}}{\Delta\check
{u}^{k'}}\hat{b}^v_{k'}\,\mathrm{d}W_t
\bigg|_{\bar{u}=\Pi}1_{\{|k-k'|\leq h\} }.
\end{eqnarray*}

Then $\mathcal{A}^1_k,\mathcal{A}^2_k,\mathcal{A}^3_{k,k'},\mathcal
{A}^4_{k,k'}$ satisfy $E[|f^{vu,(1)}_k-\tilde{f}^{vu,(1)}_k-\mathcal
{A}^1_k||\Pi]\leq C|\tilde{\theta}_{1,k}|$,
$E[|\tilde{f}^{vu,(1)}_k-\mathcal{A}^2_k|^2|\Pi]^{1/2}\leq C|\tilde
{\theta}_{1,k}|^{1/2}$ and
$E[|(\tilde{f}^{vu}_{k'}-f^{vu}_{k'}-\mu^{\star}_{(k')}(b_{(k')}b^{\star
}_{(k')})^{-1}b_{(k')}\Delta W_{k'})1_{\{|k-k'|\leq h\} }-\mathcal
{A}^3_{k,k'}-\mathcal{A}^4_{k,k'}|^2|\Pi]^{1/2}\leq C|\tilde{\theta
}_{1,k'}|1_{\{|k-k'|\leq h\} }$.\vspace*{1pt}

Hence, we obtain
%
\begin{eqnarray}
\Phi_{n,3}&\leq&C\sup_vE \biggl[
\biggl|b_n^{-{1}/{2}}\sum_k \biggl(
\mathcal{A}^1_k-\mathcal{A}^2_k
\sum_{k';|k'-k|\leq h}\bigl(\mathcal {A}^3_{k,k'}+
\mathcal{A}^4_{k,k'}\bigr) \biggr) \biggr| \Big|\Pi \biggr]
\nonumber\\
&&\label{Phi-est}{}+\mathrm{O}_p\bigl(b_n^{-{1}/{2}}r_nb_n^{\kappa}(
\ell_{1,n}+\ell_{2,n})\bigr)+\mathrm{o}_p(1)
\\
\nonumber
&=&C\sup_vE \biggl[ \biggl|b_n^{-{1}/{2}}\sum
_k \biggl(\mathcal {A}^1_k-
\mathcal{A}^2_k\sum_{k';|k'-k|\leq h}
\bigl(\mathcal {A}^3_{k,k'}+\mathcal{A}^4_{k,k'}
\bigr) \biggr) \biggr| \Big|\Pi \biggr] +\mathrm{o}_p(1),
\end{eqnarray}
where we use the Cauchy--Schwarz inequality,
tightness of $\{b_n^{-1}(\ell_{1,n}+\ell_{2,n})\}_n$ and
$r_nb_n^{1/2+\kappa}=\mathrm{o}_p(1)$ by  the definition of $\kappa$.

Moreover, by using the Burkholder--Davis--Gundy inequality, we obtain
\begin{eqnarray*}
&& \sup_vE \biggl[ \biggl|b_n^{-{1}/{2}}\sum
_k\mathcal{A}^1_k \biggr| \Big|\Pi
\biggr]\\[1pt]
 && \quad =\sup_vE \biggl[ \biggl|\frac{b_n^{-1/2}}{2}\sum
_k \biggl(\operatorname{tr}\bigl(\bigl(\hat {b}^v_k
\bigr)^{\star}\bigl(\partial_{\sigma} \bigl(\bigl(b_{(k)}b^{\star
}_{(k)}
\bigr)^{-1}\bigr)-\partial_{\sigma}\bigl(\bigl(
\tilde{b}_{(k)}\tilde{b}^{\star
}_{(k)}
\bigr)^{-1}\bigr)\bigr) \bigl(Y^{0,vu}_{\check{U}},
\sigma^n_{vu}\bigr)\hat{b}^v_k
\bigr)
\\[1pt]
&&\qquad \hspace*{80pt}{}+\partial_{\sigma}\log\frac{\det(b_{(k)}b^{\star}_{(k)})}{\det
(\tilde{b}_{(k)}\tilde{b}^{\star}_{(k)})}\bigl(Y^{0,vu}_{\check{U}},
\sigma ^n_{vu}\bigr) \biggr) \biggr| \Big|\Pi \biggr]+\mathrm{o}_p(1).
\end{eqnarray*}

Furthermore, by applying It\^{o}'s formula to
$(\partial_{\sigma} ((b_{(k)}b^{\star}_{(k)})^{-1})-\partial_{\sigma
}((\tilde{b}_{(k)}\tilde{b}^{\star}_{(k)})^{-1}))$
and  $\partial_{\sigma}\log(\det(b_{(k)}b^{\star}_{(k)})/\det(\tilde
{b}_{(k)}\tilde{b}^{\star}_{(k)}))$,
1. of Lemma~\ref{tech-lemma} in the \hyperref[app]{Appendix} and Lemma~\ref
{doubleThetaSum-est}, we have
%
\begin{eqnarray}
\sup_vE \biggl[ \biggl|b_n^{-{1}/{2}}
\sum_k\mathcal{A}^1_k \biggr| \Big|
\Pi \biggr] &\leq&Cb_n^{-{1}/{2}} \biggl(\sum
_{l_1,l_2}\sum_k\Delta\check
{U}^k1_{\{ [\check{U}^{k-1},\check{U}^k)\subset\tilde{\theta
}_{1,l_1}\cap\tilde{\theta}_{1,l_2} \} } \biggr)^{{1}/{2}}\nonumber\\[1pt]
&&\label{A1-est}{}+Cb_n^{-{1}/{2}}
\sum_k|\tilde{\theta}_{1,k}|+\mathrm{o}_p(1)
\\[1pt]
\nonumber
&\leq&Cb_n^{-{1}/{2}} \biggl(\sum
_{l_1,l_2}|\tilde{\theta }_{1,l_1}\cap\tilde{
\theta}_{1,l_2}| \biggr)^{{1}/{2}}+\mathrm{o}_p(1)=\mathrm{o}_p(1).
\end{eqnarray}

We can see $\mathcal{A}^3_{k,k'}$ can be decomposed as $\mathcal
{A}^3_{k,k'}=\sum_{\tilde{k}}\hat{\mathcal{A}}^3_{k,k',\tilde
{k}}+\mathrm{O}_p(|\tilde{\theta}_{1,k'}|)$,
where $\{\sum_{\tilde{k}\leq l}\hat{\mathcal{A}}^3_{k,k',\tilde
{k}}|_{\Pi=\bar{u}}\}_{l=0}^{L_0(\bar{u})}$ is a martingale for any
$\bar{u}\in\mathcal{U}$,
$E[|\sum_{\tilde{k};\tilde{k}<l}\hat{\mathcal{A}}^3_{k,k',\tilde
{k}}|^4|\Pi]^{1/4}\leq C|\tilde{\theta}_{1,k'}|^{1/2}$ for any $l$
and $E[|\hat{\mathcal{A}}^3_{k,k',\tilde{k}}|^4|\Pi]^{1/4}\leq C|\Delta
\check{U}^{\tilde{k}}|^{1/2}1_{\{ [\check{U}^{\tilde{k}-1},\check
{U}^{\tilde{k}})\subset\tilde{\theta}_{1,k'}\} }$.
Hence by Lemmas \ref{doubleThetaSum-est} and \ref{tech-lemma}, we obtain
%
\begin{eqnarray}
&&\hspace*{-18pt}\sup_vE \biggl[ \biggl|b_n^{-{1}/{2}}
\sum_k\mathcal{A}^2_k\sum
_{k';|k'-k|\leq h}\mathcal{A}^3_{k,k'} \biggr| \Big|
\Pi \biggr]
\nonumber
\\
&& \hspace*{-18pt}\quad \leq\sup_vE \biggl[ \biggl|b_n^{-{1}/{2}}
\sum_k\mathcal{A}^2_k\sum
_{k';|k'-k|\leq h}\hat{\mathcal{A}}^3_{k,k',k}
\biggr| \Big|\Pi \biggr] +\mathrm{O}_p \biggl(b_n^{-{1}/{2}}\sum
_k\sum_{k';|k'-k|\leq h}|\tilde{
\theta }_{1,k'}| \biggr)
\nonumber
\\
&&\hspace*{-18pt}\qquad {}+Cb_n^{-{1}/{2}}b_n^{\kappa} \biggl\{
\sum_{l_1,l_2;|l_1-l_2|\leq
2h}|\tilde{\theta}_{1,l_1}|^{{1}/{2}}|
\tilde{\theta}_{1,l_2}|^{{1}/{2}} +\sum
_{l_1,l_2}\sum_k\Delta
U^k1_{\{ [U^{k-1},U^k)\subset\tilde{\theta
}_{1,l_1}\cap\tilde{\theta}_{1,l_2}\} } \biggr\}^{{1}/{2}}
\nonumber
\\[-8pt]
\label{A2A3-est}\\[-8pt]
\nonumber
&& \hspace*{-18pt}\quad \leq\sup_vE \biggl[ \biggl|b_n^{-{1}/{2}}
\sum_k\mathcal{A}^2_k\sum
_{k';|k'-k|\leq h}\hat{\mathcal{A}}^3_{k,k',k}
\biggr| \Big|\Pi \biggr] +\mathrm{O}_p\bigl(b_n^{-{1}/{2}}b_n^{\kappa}r_n(
\ell_{1,n}+\ell_{2,n})\bigr)
\\
&&\hspace*{-18pt}\qquad {}+Cb_n^{-{1}/{2}}b_n^{\kappa} \biggl
\{r_n^{{1}/{2}}b_n^{{\kappa}/{2}}(
\ell_{1,n}+\ell_{2,n})^{{1}/{2}} + \biggl(\sum
_{l_1,l_2}|\tilde{\theta}_{1,l_1}\cap\tilde{\theta
}_{1,l_2}| \biggr)^{{1}/{2}} \biggr\}
\nonumber
\\
&&\hspace*{-18pt} \quad \leq\sup_vE \biggl[ \biggl|b_n^{-{1}/{2}}
\sum_k\mathcal{A}^2_k\sum
_{k';|k'-k|\leq h}\hat{\mathcal{A}}^3_{k,k',k}
\biggr| \Big|\Pi \biggr]+\mathrm{o}_p(1).\nonumber
\end{eqnarray}

Moreover, by using Lemma~\ref{tech-lemma} with relations $E[|\mathcal
{A}^2_k|^4|\Pi]^{1/4}\leq C$
and $E[|\mathcal{A}^4_{k,k'}|^4|\Pi]^{1/4}\leq C|\tilde{\theta
}_{1,k'}|^{1/2}$,\vspace*{-3pt} we obtain
%
\begin{eqnarray}
&& \sup_vE \biggl[ \biggl|b_n^{-{1}/{2}}
\sum_k\mathcal{A}^2_k\sum
_{k';0<|k'-k|\leq h}\mathcal{A}^4_{k,k'} \biggr| \Big|
\Pi \biggr]\nonumber\\[-3pt]
 &&\label{A2A4-est} \quad \leq Cb_n^{-{1}/{2}}b_n^{\kappa}
\biggl(\sum_{l_2,k'_2;|l_2-k'_2|\leq2h}|\tilde{\theta}_{1,l_2}|^{{1}/{2}}|
\tilde{\theta}_{1,k'_2}|^{{1}/{2}} \biggr)^{{1}/{2}}
\\[-4pt]
\nonumber
&&\quad  \leq Cb_n^{-{1}/{2}+\kappa}\bigl((\ell_{1,n}+
\ell_{2,n})\cdot \bigl(4b_n^{\kappa}+1
\bigr)r_n\bigr)^{1/2}=\mathrm{o}_p(1).
\end{eqnarray}

By (\ref{Phi-est})--(\ref{A2A4-est}), we obtain
\begin{eqnarray*}
\Phi_{n,3} &\leq &  C \biggl\{\sup_vE \biggl[
\biggl|b_n^{-{1}/{2}}\sum_k
\mathcal{A}^2_k\sum_{k';|k'-k|\leq h}
\hat{\mathcal{A}}^3_{k,k',k} \biggr| \Big|\Pi \biggr] \\
&&\hspace*{12pt}{}+\sup
_vE \biggl[ \biggl|b_n^{-{1}/{2}}\sum
_k\mathcal{A}^2_k\mathcal
{A}^4_{k,k} \biggr| \Big|\Pi \biggr] \biggr\} +\mathrm{o}_p(1).
\end{eqnarray*}
By using It\^{o}'s formula, the Burkholder--Davis--Gundy
inequality and Lemma~\ref{tech-lemma} 1. similarly,
we obtain $\Phi_{n,3}=\mathrm{o}_p(1)$.
\end{pf}

We proceed to the second step. We will prove
\[
\log\frac{\bar{\mathbb{P}}^0_u}{\bar{\mathbb{P}}^0_0}(Y_{\Pi
})-b_n^{-{1}/{2}}u\int
^1_0\sum_k
\frac{\bar{\mathbb
{P}}^{2,vu}_k(\tilde{f}^{vu,(1)}_k)}{\bar{\mathbb
{P}}^{2,vu}_k(1)}\,\mathrm{d}v(Y_{\Pi})\to^p 0
\]
as $n\to\infty$ under $[A1']$, $[A2]$ and $[A3']$.
To this end, we need to estimate
\begin{eqnarray*}
&&b_n^{-{1}/{2}}u\int^1_0
\sum_k\frac{\bar{\mathbb
{P}}^{2,vu}_k(\tilde{f}^{vu,(1)}_k1_{B^n_M})}{\tilde{\mathbb
{P}}^0_{M,vu}}\,\mathrm{d}v(Y_{\Pi})
-b_n^{-{1}/{2}}u\int^1_0\sum
_k\frac{\bar{\mathbb{P}}^{2,vu}_k(\tilde
{f}^{vu,(1)}_k1_{B^n_M})}{\bar{\mathbb{P}}^{2,vu}_k(1_{B^n_M})}\,\mathrm{d}v(Y_{\Pi
})
\\
&& \quad =b_n^{-{1}/{2}}u\int^1_0
\frac{1}{\tilde{\mathbb{P}}^0_{M,vu}}\sum_k\frac{\bar{\mathbb{P}}^{2,vu}_k(\tilde{f}^{vu,(1)}_k1_{B^n_M})}{\bar
{\mathbb{P}}^{2,vu}_k(1_{B^n_M})}
\\
&&\hspace*{44pt}\qquad{}\times\biggl(\exp \biggl(\sum_{|k-k'|\leq
h}\bigl(
\tilde{f}^{vu}_{k'}-\log\check{p}^0_{k',vu}
\bigr) \biggr)-1 \biggr)1_{B^n_M}\mathbb{P}^0_{vu}\,\mathrm{d}
\hat{z}\,\mathrm{d}v(Y_{\Pi}).
\end{eqnarray*}
If $\bar{\mathbb{P}}^{2,vu}_k(\tilde{f}^{vu,(1)}_k1_{B^n_M})/\bar
{\mathbb{P}}^{2,vu}_k(1_{B^n_M})$ has \textit{good} properties,
we can apply techniques from It\^{o} calculus similarly
to the proof of Lemma~\ref{coeff-move-lemma1}.
So we will investigate properties of this quantity.

We prepare some additional lemmas.

%
\begin{lemma}\label{P2P1-tight}
Let $u\in\mathbb{R}^d$. Assume $[A1'],[A2]$ and $[A3']$. Then $\{\sup_{k,v}|\log(\bar{\mathbb{P}}^{2,vu}_k(1)/\break\bar{\mathbb{P}}^0_0)|(Y_{\Pi
})\}_n$ is tight.
\end{lemma}

\begin{pf}
Let $\mathbb{P}^{2,vu}_k=\exp(\sum_{k'}\acute{f}^k_{k'})$,
$\mathcal{K}'_2(\bar{u})=\{z;\sup_{k,v}|\tilde{f}^{vu}_k-\log\check
{p}^0_{k,vu}|(z)\leq b_n^{-2\kappa}\}$,
and $Y^{2,vu,k}_{\check{U}}$ be random variables with the $\Pi
$-conditional distribution $\mathbb{P}^{2,vu}_k(z)\,\mathrm{d}\hat{z}\,\mathrm{d}\bar
{z}P_{Y_0}(\mathrm{d}z_0)$.
For any $q>0$, let $q'\geq2(q+1)/(1-\delta_1-4\kappa)$. Then we obtain
%
\begin{eqnarray}
&&P\bigl[Y_{\check{U}}\in\bigl(\mathcal{K}'_2
\bigr)^c(\Pi)|\Pi\bigr]
\nonumber
\\
&& \quad \leq b_n^{2\kappa q'}\sum_kE
\Bigl[\sup_v\bigl|\tilde{f}^{vu}_k-\log
\check {p}^0_{k,vu}\bigr|^{q'}(Y_{\check{U}}) |\Pi
\Bigr]
\nonumber
\\
&&\label{P2P1-tight-est1}
\quad \leq Cb_n^{2\kappa q'}\sum_kE
\Bigl[\sup_v\bigl|\tilde {f}^{vu}_k-f^{vu}_k\bigr|^{q'}(Y_{\check{U}})
|\Pi \Bigr]\\
&&\nonumber\qquad{}+Cb_n^{2\kappa
q'}\sup_r\sum
_kE \biggl[\sup_v \biggl|
\frac{\partial_r\check
{p}^r_{k,vu}}{\check{p}^r_{k,vu}} \biggr|^{q'}(Y_{\check{U}}) \Big|\Pi \biggr]
\\
&&\quad \leq Cb_n^{2\kappa q'}r_n^{q'/2}(
\ell_{1,n}+\ell_{2,n}) =\mathrm{O}_p\bigl(b_n^{-q}
\bigr).\nonumber
\end{eqnarray}

Similarly, we have $\sup_vP[Y^{0,vu}_{\check{U}}\in(\mathcal
{K}'_2)^c(\Pi)|\Pi]=\mathrm{O}_p(b_n^{-q})$
and $\sup_{k,v}P[Y^{2,vu,k}_{\check{U}}\in (\mathcal{K}'_2)^c(\Pi)|\Pi
]=\mathrm{O}_p(b_n^{-q})$ for any $q>0$.

Hence, we have
\begin{eqnarray*}
&&P \biggl[\sup_{k,v} \biggl|\frac{\bar{\mathbb{P}}^{2,vu}_k(1_{(\mathcal
{K}'_2)^c})}{\bar{\mathbb{P}}^0_{vu}}
\biggr|(Y_{\Pi})>\frac{b_n^{-1}}{2}, \sup_v \biggl|\log
\frac{\bar{\mathbb{P}}^0_{vu}}{\bar{\mathbb
{P}}^0_0} \biggr|(Y_{\Pi})\leq M' \Big|\Pi \biggr]
\\
&& \quad \leq 2\mathrm{e}^{M'}b_n\int\sup_{k,v}
\mathbb{P}^{2,vu}_k1_{(\mathcal
{K}'_2)^c}\,\mathrm{d}\hat{z}\,\mathrm{d}
\bar{z}P_{Y_0}(\mathrm{d}z_0) \bigg|_{\bar{u}=\Pi}
\\
&&\quad \leq 2\mathrm{e}^{M'}b_n\sum_k
\int \biggl(\mathbb{P}^{2,0}_k+\int^1_0\bigl|
\partial _v\mathbb{P}^{2,vu}_k\bigr|\,\mathrm{d}v
\biggr)1_{(\mathcal{K}'_2)^c}\,\mathrm{d}\hat{z}\,\mathrm{d}\bar {z}P_{Y_0}(\mathrm{d}z_0)
\bigg|_{\bar{u}=\Pi}
\\
&&\quad \leq  2\mathrm{e}^{M'}b_n\sum_kP
\bigl[Y^{2,0,k}_{\check{U}}\in\bigl(\mathcal {K}'_2
\bigr)^c|\Pi\bigr] + 2\mathrm{e}^{M'}b_n\sum
_k\sup_vE \biggl[ \biggl|\frac{\partial
_v\mathbb{P}^{2,vu}_k}{\mathbb{P}^{2,vu}_k}
\biggr|1_{(\mathcal
{K}'_2)^c}\bigl(Y^{2,vu,k}_{\check{U}}\bigr) \Big|\Pi \biggr]
\\
&& \quad = \mathrm{O}_p\bigl(b_n^{-q}\bigr)
\end{eqnarray*}
for any $q>0$ and $M'>0$.

Therefore, for any $\varepsilon>0$, there exists $N>0$ such that
%
\begin{equation}
\label{P2P1-tight-est2}
P \biggl[\sup_{k,v} \biggl|\frac{\bar{\mathbb{P}}^{2,vu}_k(1_{(\mathcal
{K}'_2)^c})}{\bar{\mathbb{P}}^0_{vu}}
\biggr|(Y_{\Pi})>\frac
{b_n^{-1}}{2} \biggr]<\varepsilon
\end{equation}
for $n\geq N$ by Lemma~\ref{Pvu-tightness}.

Moreover, for any $\delta>0$ and $M'>0$, we have
\begin{eqnarray*}
&&P \biggl[\sup_{k,v} \biggl|\frac{\bar{\mathbb{P}}^{2,vu}_k(1_{\mathcal
{K}'_2})}{\bar{\mathbb{P}}^0_{vu}}-1
\biggr|(Y_{\Pi})>\delta, \sup_v \biggl|\log
\frac{\bar{\mathbb{P}}^0_{vu}}{\bar{\mathbb{P}}^0_0} \biggr|(Y_{\Pi})\leq M' \Big|\Pi \biggr]
\\
&&\quad \leq \frac{\mathrm{e}^{M'}}{\delta}\int\sup_{k,v}\bigl|\mathbb
{P}^{2,vu}_k1_{\mathcal{K}'_2}-\mathbb{P}^0_{vu}\bigr|\,\mathrm{d}
\hat{z}\,\mathrm{d}\bar {z}P_{Y_0}(\mathrm{d}z_0) \bigg|_{\bar{u}=\Pi}
\\
&&\quad \leq \frac{\mathrm{e}^{M'}}{\delta}\int \biggl\{\sup_k\bigl|\mathbb
{P}^{2,0}_k1_{\mathcal{K}'_2}-\mathbb{P}^0_0\bigr|+
\int^1_0\sup_k\bigl|\partial
_v\mathbb{P}^{2,vu}_k1_{\mathcal{K}'_2}-
\partial_v\mathbb {P}^0_{vu}\bigr|\,\mathrm{d}v \biggr\}\,\mathrm{d}
\hat{z}\,\mathrm{d}\bar{z}P_{Y_0}(\mathrm{d}z_0) \bigg|_{\bar{u}=\Pi}
\\
&&\quad \leq \frac{\mathrm{e}^{M'}}{\delta} \biggl\{P\bigl[Y_{\check{U}}\in\bigl(\mathcal
{K}'_2\bigr)^c(\Pi)|\Pi
\bigr]+\mathrm{O}_p\bigl(b_n^{-\kappa}\bigr)+\sup
_vE \biggl[ \biggl|\frac{\partial
_v\mathbb{P}^0_{vu}}{\mathbb{P}^0_{vu}} \biggr|1_{(\mathcal
{K}'_2)^c}
\bigl(Y^{0,vu}_{\check{U}}\bigr) \Big|\Pi \biggr]
\\
&&\qquad\hspace*{20pt} {}+\sup_v\int\sup_k \biggl|\sum
_{k''}\partial_v\acute{f}^{k,vu}_{k''}
\exp \biggl(\sum_{|k'-k|\leq h}\bigl(\tilde{f}^{vu}_{k'}-
\log\check {p}^0_{k',vu}\bigr) \biggr)
-\sum
_{k'}\frac{\partial_v\check
{p}^0_{k',vu}}{\check{p}^0_{k',vu}} \biggr|\\
&&\hspace*{57pt}\qquad{}\times 1_{\mathcal{K}'_2}\mathbb
{P}^0_{vu}\,\mathrm{d}\hat{z}\,\mathrm{d}\bar{z}P_{Y_0}(\mathrm{d}z_0)
\Big|_{\bar{u}=\Pi} \biggr\} \\
&& \quad =\mathrm{o}_p(1),
\end{eqnarray*}
by (\ref{P2P1-tight-est1}). Hence,
%
\begin{equation}
\label{P2P1-tight-est3}
P \biggl[\sup_{k,v} \biggl|\frac{\bar{\mathbb{P}}^{2,vu}_k(1_{\mathcal
{K}'_2})}{\bar{\mathbb{P}}^0_{vu}}-1
\biggr|(Y_{\Pi})>\delta \biggr]<\varepsilon
\end{equation}
for sufficiently large $n$.

Lemmas \ref{density-est} and \ref{Pvu-tightness}, (\ref{P2P1-tight-est2}) and (\ref{P2P1-tight-est3}) complete the proof.
\end{pf}

Let $\mathcal{K}''_{2,M}=\{(z_0,\bar{z});\inf_{k,v}\bar{\mathbb
{P}}^{2,vu}_k(1_{B^n_M})(z_0,\bar{z})>0\}$ and
\[
\mathcal{K}^2_M(\bar{u})=\mathcal{K}^1_M(
\bar{u})\cap \biggl\{(z_0,\bar {z});\sup_{k,v} \biggl|
\frac{\bar{\mathbb{P}}^{2,vu}_k(\tilde
{f}^{vu,(1)}_k1_{B^n_{M}})}{\bar{\mathbb
{P}}^{2,vu}_k(1_{B^n_{M}})}-\frac{\bar{\mathbb{P}}^{2,vu}_k(\tilde
{f}^{vu,(1)}_k)}{\bar{\mathbb{P}}^{2,vu}_k(1)} \biggr|(z_0,\bar{z})\leq
b_n^{-1} \biggr\} \cap\mathcal{K}''_{2,M}
\]
for $u\in\mathbb{R}^d$, $\bar{u}\in\mathcal{U}$, $1\leq k\leq L_0(\bar
{u})$ and $M>0$.

\begin{lemma}\label{L2M-est}
Let $u\in\mathbb{R}^d$. Assume $[A1'],[A2]$ and $[A3']$.
Then for any $\varepsilon>0$, there exists $M'>0$ and $\{N_M\}_{M\geq
M'}\subset\mathbb{N}$ such that $\sup_{n\geq N_M}P[Y_{\Pi}\in(\mathcal
{K}^2_M)^c(\Pi)]<\varepsilon$ for $M\geq M'$.
\end{lemma}

\begin{pf}
By (\ref{P2P1-tight-est2}),
for any $\varepsilon>0$, there exists $N'_1\in\mathbb{N}$ such that $P[\sup_{k,v}(\bar{\mathbb{P}}^{2,vu}_k(1_{B^n_M}1_{(\mathcal{K}'_2)^c})/\allowbreak \bar
{\mathbb{P}}^0_{vu})(Y_{\Pi})>b_n^{-1}/2]<\varepsilon$ for $n\geq N'_1$
and $M>0$.
Moreover,
by (\ref{BnMc-est}), we have
\[
P \biggl[\sup_{k,v} \biggl|\frac{\bar{\mathbb
{P}}^{2,vu}_k(1_{B^n_M}1_{\mathcal{K}'_2})}{\bar{\mathbb
{P}}^0_{vu}}-1
\biggr|(Y_{\Pi})>\delta \biggr]<\varepsilon
\]
for any $\delta,\varepsilon>0$ and sufficiently large $n$ and $M$,
similarly to the derivation of (\ref{P2P1-tight-est3}).

Therefore, there exist $N'_2\in\mathbb{N}$ and $M_2>0$ such that
\[
P\Bigl[\inf_{k,v}\bar{\mathbb{P}}^{2,vu}_k(1_{B^n_M})
(Y_{\Pi})>0\Bigr]>1-\varepsilon
\]
and
%
\begin{equation}
\label{logP2k-est}
P\Bigl[\sup_{k,v}\bigl|\log\bigl(\bar{
\mathbb{P}}^{2,vu}_k(1_{B^n_M})/\bar{\mathbb
{P}}^0_{vu}\bigr)\bigr|(Y_{\Pi})>\delta\Bigr]<
\varepsilon
\end{equation}
for $M>M_2$ and $n\geq N'_2$.
Moreover, we have $\sup_{k,v}|\bar{\mathbb{P}}^{2,vu}_k(\tilde
{f}_k^{vu,(1)})/\bar{\mathbb{P}}^0_{vu}|(Y_{\Pi})=\mathrm{O}_p(b_n^2)$.

Since
\begin{eqnarray*}
&&\frac{\bar{\mathbb{P}}^{2,vu}_k(\tilde{f}^{vu,(1)}_k1_{B^n_M})}{\bar
{\mathbb{P}}^{2,vu}_k(1_{B^n_M})}-\frac{\bar{\mathbb{P}}^{2,vu}_k(\tilde
{f}^{vu,(1)}_k)}{\bar{\mathbb{P}}^{2,vu}_k(1)}
\\
&& \quad =-\frac{\bar{\mathbb{P}}^0_{vu}}{\bar{\mathbb
{P}}^{2,vu}_k(1_{B^n_M})}\frac{\bar{\mathbb{P}}^{2,vu}_k(\tilde
{f}^{vu,(1)}_k1_{(B^n_M)^c})}{\bar{\mathbb{P}}^0_{vu}}\\
&&\qquad {}+\frac{\bar{\mathbb{P}}^0_{vu}}{\bar{\mathbb
{P}}^{2,vu}_k(1_{B^n_M})}
\frac{\bar{\mathbb{P}}^0_{vu}}{\bar{\mathbb
{P}}^{2,vu}_k(1)}\frac{\bar{\mathbb{P}}^{2,vu}_k(\tilde
{f}^{vu,(1)}_k)}{\bar{\mathbb{P}}^0_{vu}} \frac{\bar{\mathbb{P}}^{2,vu}_k(1_{(B^n_M)^c})}{\bar{\mathbb
{P}}^0_{vu}},
\end{eqnarray*}
there exist $M',M_1>0$ and $\{N_M\}_{M\geq M'}\subset\mathbb{N}$ such that
\begin{eqnarray*}
&&P\bigl[Y_{\Pi}\in\bigl(\mathcal{K}^2_M
\bigr)^c(\Pi)\bigr]
\\
&& \quad \leq P\bigl[Y_{\Pi}\in\bigl(\mathcal{K}^1_M
\bigr)^c(\Pi) \cup\bigl(\mathcal{K}''_{2,M}
\bigr)^c\bigr] \\
&&\qquad{}+P \biggl[\sup_{k,v} \biggl|2
\frac{\bar{\mathbb{P}}^{2,vu}_k(|\tilde
{f}^{vu,(1)}_k|1_{(B^n_M)^c})}{\bar{\mathbb{P}}^0_{vu}} +4b_n^3\frac{\bar{\mathbb{P}}^{2,vu}_k(1_{(B^n_M)^c})}{\bar{\mathbb
{P}}^0_{vu}}
\biggr|(Y_{\Pi})>b_n^{-1} \biggr]+\varepsilon
\\
&& \quad \leq P \biggl[\sup_{k,v} \biggl|2\frac{\bar{\mathbb{P}}^{2,vu}_k(|\tilde
{f}^{vu,(1)}_k|1_{(B^n_M)^c}1_{\mathcal{K}'_2})}{\bar{\mathbb{P}}^0_{vu}}
+4b_n^3\frac{\bar{\mathbb{P}}^{2,vu}_k(1_{(B^n_M)^c}1_{\mathcal
{K}'_2})}{\bar{\mathbb{P}}^0_{vu}} \biggr|(Y_{\Pi})>
\frac
{b_n^{-1}}{2},\\
&&\hspace*{34pt} \sup_v \biggl|\log\frac{\bar{\mathbb{P}}^0_{vu}}{\bar
{\mathbb{P}}^0_0}
\biggr|(Y_{\Pi})\leq M_1 \biggr] +5\varepsilon
\\
&& \quad \leq E \biggl[ \biggl(4b_n\mathrm{e}^{M_1}\int\sup
_{k,v} \bigl\{\bigl(2\bigl|\tilde {f}^{vu,(1)}_k\bigr|+4b_n^3
\bigr)\mathbb{P}^0_{vu} \bigr\} 1_{(B^n_M)^c}1_{\mathcal{K}'_2}\,\mathrm{d}
\hat{z}\,\mathrm{d}\bar{z}P_{Y_0}(\mathrm{d}z_0) \bigg|_{\bar
{u}=\Pi} \biggr)
\wedge1 \biggr] +5\varepsilon
\\
&&\quad \leq \mathit{CE} \biggl[ \biggl(b_n^4P\bigl[Y_{\check{U}}
\in\bigl(B^n_M\bigr)^c|\Pi\bigr] +
b_n^{-2}(\ell_{1,n}+\ell_{2,n})
\\
&&\hspace*{25pt}\qquad {}+ \sup_vE \biggl[ \biggl|\frac{\partial_v\mathbb{P}^0_{vu}}{\mathbb
{P}^0_{vu}} \biggr|
\biggl(b_n^4+b_n\sum
_k\sup_v\bigl|\tilde{f}^{vu,(1)}_k\bigr|
\biggr)1_{(B^n_M)^c}\bigl(Y^{0,vu}_{\check{U}}\bigr) \Big|\Pi
\biggr] \biggr)\wedge1 \biggr]+5\varepsilon\\
&&\quad < 6\varepsilon
\end{eqnarray*}
for $M\geq M'$ and $n\geq N_M$, by (\ref{logP2k-est}) and similar
arguments to (\ref{P2P1-tight-est1}) and (\ref{P2P1-tight-est2}).
\end{pf}

In the following,\vspace*{1pt} we see that the integral $\int\exp(\sum_{\nu\leq
k'\leq\chi}\tilde{f}^{vu}_{k'})\,\mathrm{d}\hat{z}$  has a simple
representation of a function of increments
for $1\leq\nu\leq \chi\leq L_0(\bar{u})$.
To see this, we will define some notation related to the observation
times and increments of processes in the interval $(\check{u}^{\nu
-1},\check{u}^{\chi}]$.

For $\bar{u}=((s^i)^i,(t^j)^j)\in\mathcal{U}$ and $1\leq\nu\leq \chi
\leq L_0(\bar{u})$,\vspace*{-2pt}
let
\begin{eqnarray*}
\mathcal{I}_{\nu,\chi} &=& \bigl(\bigl\{k_1(i)\bigr\}_i
\cap[\nu,\chi-1]\bigr)\cup\{\nu-1,\chi \},\\[-2pt]
 \mathcal{J}_{\nu,\chi} &=& \bigl(\bigl
\{k_2(j)\bigr\}_j\cap[\nu,\chi-1]\bigr)\cap\{ \nu-1,\chi
\},
\\[-2pt]
\mathcal{I}^-_{\nu,\chi}& =& \mathcal{I}_{\nu,\chi}\setminus\{\chi\},\qquad
\mathcal{J}_{\nu,\chi}^-=\mathcal{J}_{\nu,\chi}\setminus\{\chi\}.
\end{eqnarray*}
Moreover, for $z=(x_k,y_k)_{k=0}^{L_0(\bar{u})}$, define
\begin{eqnarray*}
\hat{z}_{\nu,\chi} &=& \bigl((x_k)_{\nu\leq k\leq\chi-1,k\notin\{k_1(i)\}
_i},(y_k)_{\nu\leq k\leq\chi-1,k\notin\{k_2(j)\}_j}
\bigr),
\\
x(\mathcal{I}_{\nu,\chi}) &=& \{x_{i_k}-x_{i_{k-1}}
\}_{k=1}^{L^1_{\nu,\chi
}},\qquad  y(\mathcal{J}_{\nu,\chi})=
\{y_{j_k}-y_{j_{k-1}}\} _{k=1}^{L^2_{\nu,\chi}},
\end{eqnarray*}
where $\mathcal{I}_{\nu,\chi}=\{i_k\}_{k=0}^{L^1_{\nu,\chi}}, \mathcal
{J}_{\nu,\chi}=\{j_k\}_{k=0}^{L^2_{\nu,\chi}}$ and $\nu-1=i_0<\cdots
<i_{L^1_{\nu,\chi}}=\chi$, $\nu-1=j_0<\cdots<j_{L^2_{\nu,\chi}}=\chi$.
For $p=1,2$, $k\in\mathcal{I}_{\nu,\chi}$ and $l\in\mathcal{J}_{\nu,\chi
}$, let $\check{b}^p_k=\check{b}^{n,p}_{k,vu}=b^p(\check
{u}^{k-1},x_{i(k)},y_{j(k)},\sigma^n_{vu})$, $\tilde{I}_{\nu,\chi
}^k=[s^{i-1},s^i)\cap[\check{u}^{\nu-1},\check{u}^{\chi})$, $\tilde
{J}_{\nu,\chi}^k=[t^{j-1},t^j)\cap[\check{u}^{\nu-1},\check{u}^{\chi})$,
where $i,j$ satisfy $[\check{u}^{k-1},\check{u}^k)\subset
[s^{i-1},s^i)$ and $[\check{u}^{l-1},\check{u}^l)\subset[t^{j-1},t^j)$.
Let $L_{\nu,\chi}=L^1_{\nu,\chi}+L^2_{\nu,\chi}$, $\tilde
{K}_{k'}=[\check{u}^{k'-1},\check{u}^{k'})$ and
\begin{eqnarray*}
S_{\nu,\chi}&=& \left( %
\begin{array}{cc}\ds \operatorname{diag}\biggl(
\biggl\{\sum_{k'}\bigl|\check{b}^1_{k'}\bigr|^2\bigl|
\tilde{I}^{i_k}_{\nu,\chi
}\cap\tilde{K}_{k'}\bigr|\biggr
\}_{1\leq k\leq L^1_{\nu,\chi}}\biggr)
\\
\ds\biggl\{\sum_{k'}\check{b}^1_{k'}
\cdot\check{b}^2_{k'}\bigl|\tilde{I}_{\nu,\chi
}^{i_k}
\cap\tilde{J}^{j_l}_{\nu,\chi}\cap\tilde{K}_{k'}\bigr|\biggr
\}_{1\leq
l\leq L^2_{\nu,\chi},1\leq k\leq L^1_{\nu,\chi}} \end{array}\right.\\
&&\hspace*{6pt} {}\left.\begin{array}{cc}\ds\biggl\{\sum_{k'}
\check{b}^1_{k'}\cdot\check{b}^2_{k'}\bigl|
\tilde{I}_{\nu,\chi
}^{i_k}\cap\tilde{J}^{j_l}_{\nu,\chi}
\cap\tilde{K}_{k'} \bigr|\biggr\}_{1\leq
k\leq L^1_{\nu,\chi},1\leq l\leq L^2_{\nu,\chi}}   \\ \ds\operatorname{diag}\biggl(\biggl\{\sum_{k'}\bigl|
\check{b}^2_{k'}\bigr|^2\bigl|\tilde{J}^{j_l}_{\nu,\chi
}
\cap\tilde{K}_{k'}\bigr|\biggr\}_{1\leq l\leq L^2_{\nu,\chi}}\biggr) \end{array}
 \right).
\end{eqnarray*}

Let $\varphi(x;V)$ be the density function of $N(0,V)$ for a
symmetric, positive definite matrix $V$.
The following lemma enables us to calculate integrals of exponential
functions of $\tilde{f}^{vu}_k$.

%
\begin{lemma}\label{induction-lemma}
Let $u\in\mathbb{R}^d,\bar{u}\in\mathcal{U},n\in\mathbb{N}$ and $1\leq
\nu\leq\chi\leq L_0(\bar{u})$. Assume $[A1]$. Then $\det S_{\nu,\chi
}>0$ and
%
\begin{equation}
\label{induction-lemma-eq}
\int\exp \biggl(\sum_{\nu\leq k'\leq\chi}
\tilde{f}^{vu}_{k'}(z) \biggr)\,\mathrm{d}\hat{z}_{\nu,\chi} =
\varphi\bigl(\bigl(x(\mathcal{I}_{\nu,\chi})^{\star},y(
\mathcal{J}_{\nu
,\chi})^{\star}\bigr)^{\star};S_{\nu,\chi}
\bigr). 
\end{equation}
\end{lemma}

\begin{pf}
We see $\det S_{\nu,\chi}>0$ by  a similar argument to the
proof of Proposition~1 in Ogihara and Yoshida \cite{ogi-yos},
so we omit the details.

We prove (\ref{induction-lemma-eq}) by induction on $\chi$.
The results obviously hold true for $\chi=\nu$.

 Let $\chi>\nu$ and assume the results hold for $\chi-1$.
We give the proof only for the case $\check{u}^{\chi-1} \notin(s^i)_i$ and $\check{u}^{\chi-1}\in(t^j)_j$. The other cases are
proved similarly.

By the induction assumption, we obtain
\begin{eqnarray*}
\label{induction-Lemma-eq1} 
&& \int\exp \biggl(\sum_{\nu\leq k'\leq\chi}
\tilde{f}_{k'}^{vu} \biggr)\,\mathrm{d}\hat {z}_{\nu,\chi} \\
&& \quad =\int
\varphi\bigl(\bigl(x(\mathcal{I}_{\nu,\chi-1})^{\star},y(\mathcal
{J}_{\nu,\chi-1})^{\star}\bigr)^{\star};S_{\nu,\chi-1}\bigr)
\varphi\bigl(z_{\chi
}-z_{\chi-1};\tilde{b}^{vu}_{\chi}
\bigl(\tilde{b}^{vu}_{\chi}\bigr)^{\star}\Delta
\check{u}^{\chi}\bigr)\,\mathrm{d}x_{\chi-1}.
\end{eqnarray*}

Let $Z_1$ and $Z_2$ be random variables independent of each other,
satisfying $Z_1\sim N(0,S_{\nu,\chi-1})$ and $Z_2\sim N(0,\tilde
{b}^{vu}_{\chi}(\tilde{b}^{vu}_{\chi})^{\star}\Delta\check{u}^{\chi})$.
Moreover, let $\mathcal{D}$ be an $(L_{\nu,\chi-1}+1)\times(L_{\nu,\chi
-1}+2)$ matrix  with $\mathcal{D}_{pq}=\delta_{p,q}$ for $1\leq p,q\leq
L_{\nu,\chi-1}$,
$\mathcal{D}_{pq}=1$ for $(p,q)=(L^1_{\nu,\chi-1},L_{\nu,\chi-1}+1)$
or
$(L_{\nu,\chi-1}+1,L_{\nu,\chi-1}+2)$, and $\mathcal{D}_{pq}=0$ in
other cases.
Then the covariance matrix of $\mathcal{D}(Z_1^{\star},Z_2^{\star
})^{\star}$ is
\[
\mathcal{D}\left(
\begin{array}{c@{\quad}c} S_{\nu,\chi-1} & 0
\\
0 & \tilde{b}^{vu}_{\chi}\bigl(\tilde{b}^{vu}_{\chi}
\bigr)^{\star}\Delta\check {u}^{\chi} \end{array} %
 \right)\mathcal{D}^{\star}=S_{\nu,\chi}.
\]
Hence, we obtain the result by considering relations between densities
of $Z_1,Z_2$ and  $\mathcal{D}(Z_1^{\star},\allowbreak Z_2^{\star})^{\star}$.
\end{pf}

\begin{remark}
We emphasize that we can prove the above lemma because $\tilde
{b}^{vu}_{\chi}$ does not depend on~$\hat{z}_{\nu,\chi}$.
\end{remark}

We  now give another representation of $\bar{\mathbb
{P}}^{2,vu}_k(\tilde{f}^{vu,(1)}_k)/\bar{\mathbb
{P}}^{2,vu}_k(1)(z_0,\bar{z})$ consisting of a quadratic form of increments.
This representation is useful to apply It\^{o}'s rule and
martingale properties.

Let $\Theta(n,k,1;\bar{u})=\{i;1\leq i\leq L^1, s^{i-1}>\check
{u}^{(k-h-1)\vee0},s^i<\check{u}^{(k+h)\wedge L_0}\}$,
$\Theta(n,k,2;\bar{u})=\{j;1\leq j\leq L^2, t^{j-1}>\check
{u}^{(k-h-1)\vee0},t^j<\check{u}^{(k+h)\wedge L_0}\}$, $\mathbf{M}=\sharp(\Theta(n,k,1;\bar{u}))+\sharp(\Theta(n,k,2;\bar{u}))$
and
\[
\mathcal{Z}_k= \biggl( \biggl(\frac{x_{k_1(i)}-x_{k_1(i-1)}}{\sqrt
{s^i-s^{i-1}}}
\biggr)^{\star}_{i\in\Theta(n,k,1;\bar{u}) }, \biggl(\frac
{y_{k_2(j)}-y_{k_2(j-1)}}{\sqrt{t^j-t^{j-1}}}
\biggr)^{\star}_{j\in\Theta
(n,k,2;\bar{u}) } \biggr)^{\star}
\]
for $\bar{u}=((s^i)_{i=0}^{L^1},(t^j)_{j=0}^{L^2})\in\mathcal{U}$ and
$z=((x_k)_{k=0}^{L_0(\bar{u})},(y_k)_{k=0}^{L_0(\bar{u})})\in\mathbb
{R}^{2(L_0(\bar{u})+1)}$.

\begin{lemma}\label{integration-lemma}
Let $u\in\mathbb{R}^d$. Assume $[A1'],[A2]$ and $[A3']$.
Then\vspace*{1pt} there exist an $\mathbb{R}^d\otimes\mathbb{R}^{\mathbf{M}}\otimes
\mathbb{R}^{\mathbf{M}}$-valued function $\mathcal{Q}^{k,v}_1(z_0,\bar
{z},\bar{u})$,
$\mathbb{R}^d$-valued functions $\{\mathcal{Q}^{k,v}_{p}(z_0,\bar
{z},\bar{u})\}_{p=2}^4$
$(v\in[0,1]$, $\bar{u}=((s^i)_i,(t^j)_j)\in\mathcal{U}$, $1\leq k\leq
L_0(\bar{u}))$ and a constant $C>0$ such that
$\mathcal{Q}_1^{k,v}(Y^{r,v'u}_{\bar{u}},\bar{u})$ and
$\mathcal{Q}^{k,v}_2(Y^{r,v'u}_{\bar{u}},\bar{u})$ are
$\mathcal{F}_{\inf\tilde{\theta}((k-h-1)\vee1,2;\bar{u})}$-measurable,
\begin{eqnarray*}
\sup_{n,v,\bar{u},k,z_0,\bar{z}}\bigl(\bigl\Vert\mathcal{Q}_1^{k,v}(z_0,
\bar {z};\bar{u})\bigr\Vert\vee\bigl|\mathcal{Q}_2^{k,v}(z_0,
\bar{z};\bar {u})\bigr|\bigr)&\leq&C,
\\
\sup_{r,k,v,v'}E\bigl[\bigl|\mathcal{Q}^{k,v}_3
\bigl(Y^{r,v'u}_{\Pi};\Pi \bigr)\bigr|^q|\Pi
\bigr]^{1/q}&\leq&Cr_n^{1/2}b_n^{3\kappa/2}
\qquad \mbox{a.s.},
\end{eqnarray*}
$\sup_{k,v}|\mathcal{Q}^{k,v}_4(Y_{\Pi};\Pi)|=\mathrm{o}_p(b_n^{-q})$ and
%
\begin{equation}
\label{integration-lemma-eq}
\frac{\bar{\mathbb{P}}^{2,vu}_k(\tilde{f}^{vu,(1)}_k)}{\bar{\mathbb
{P}}^{2,vu}_k(1)}(z_0,\bar{z}) =
\mathcal{Z}_k^{\star}\mathcal{Q}_1^{k,v}(z_0,
\bar{z},\bar{u})\mathcal {Z}_k1_{\{\mathcal{Z}_k\neq\varnothing\} }+\sum
_{p=2}^4\mathcal {Q}^{k,v}_p(z_0,
\bar{z},\bar{u})
\end{equation}
for $v,v',r\in[0,1]$, $n\geq n_u$, $q>0$, $\bar
{u}=((s^i)_i,(t^j)_j)\in\mathcal{U}$ and $1\leq k\leq L_0(\bar{u})$.
Moreover,
%
\begin{equation}
\label{integration-lemma-state2} \sup_{k,v} \biggl|\frac{\bar{\mathbb{P}}^{2,vu}_k(\tilde
{f}_k^{vu,(1)})}{\bar{\mathbb{P}}^{2,vu}_k(1)} -
\frac{\int\tilde{f}^{vu,(1)}_k\exp(\sum_{k'}\tilde{f}^{vu}_{k'})\,\mathrm{d}\hat
{z}}{\int\exp(\sum_{k'}\tilde{f}^{vu}_{k'})\,\mathrm{d}\hat{z}} \biggr|(Y_{\Pi
})=\mathrm{o}_p\bigl(b_n^{-q}
\bigr)
\end{equation}
for any $q>0$.
\end{lemma}

\begin{pf}
We only consider the case that $L^1_{k-h,k-1}\wedge L^2_{k-h,k-1}
\wedge L^1_{k+1,k+h}\wedge L^2_{k+1,k+h}>1$,
$k\geq h+1, k+h\leq L_0(\bar{u})$, $\check{u}^{k-1},\check
{u}^k\notin(s^i)_i$ and $\check{u}^{k-1},\check{u}^k\in(t^j)_j$.
Other cases are proved in  a similar way.\vspace*{1pt}

The proof is rather complicated. It is divided in several
steps.

\begin{step}
In this step, we will have an expression of $\bar{\mathbb
{P}}^{2,vu}_k(\tilde{f}^{vu,(1)}_k)/\bar{\mathbb{P}}^{2,vu}_k(1)$
similar to (\ref{integration-lemma-eq})
by using elementary formulas (Lemma~\ref{induction-lemma} and Lemma~\ref
{dist-comb-lemma} in the \hyperref[app]{Appendix}) of Gaussian distributions.\vspace*{1pt}

Let\vspace*{1.5pt} $L_k=L_{k-h,k-1}$, $L^j_k=L^j_{k-h,k-1}$ for $j=1,2$,
$\mathcal{I}_k=\mathcal{I}_{k-h,k-1}$, $\mathcal{J}_k=\mathcal
{J}_{k-h,k-1}$ and $S_k=S_{k-h,k-1}$.
Then Lemma~\ref{induction-lemma} yields
%
\begin{eqnarray}
&&\int\tilde{f}^{vu,(1)}_k\exp \biggl(\sum
_{k';k'\leq k}\acute {f}^{k,vu}_{k'}
\biggr) (z)\,\mathrm{d}\hat{z}_{1,k}
\nonumber
\\
&&\quad =\int \biggl(-\frac{1}{2}(z_k-z_{k-1})^{\star}
\frac{\partial_{\sigma
}((\tilde{b}_{(k)}\tilde{b}^{\star}_{(k)})^{-1})(z,\sigma
^n_{vu})}{\Delta\check{u}^k }(z_k-z_{k-1})\nonumber\\
&&\label{integration-lemma-equ1}\hspace*{15pt} \qquad{}-\frac{1}{2}
\partial_{\sigma
}\log\det\bigl(\tilde{b}_{(k)}
\tilde{b}^{\star}_{(k)}\bigr) \bigl(z,\sigma^n_{vu}
\bigr) \biggr)
\\
&&\hspace*{19pt}\quad {}\times\varphi\bigl(z_k-z_{k-1};\Delta
\check{u}^k\tilde {b}^{vu}_k\bigl(
\tilde{b}^{vu}_k\bigr)^{\star}\bigr)\varphi\bigl(
\bigl(x(\mathcal{I}_k)^{\star
},y(\mathcal{J}_k)^{\star}
\bigr)^{\star};S_k\bigr)\,\mathrm{d}x_{k-1}\nonumber\\
&&\hspace*{19pt}\quad {}\times\exp \biggl(\sum
_{k'\leq k-h-1}f^{vu}_{k'} \biggr)\,\mathrm{d}
\hat{z}_{1,k-h}.
\nonumber
\end{eqnarray}

Let
$\mathcal{M}_1=((S^{-1}_k)_{ij})_{i,j\neq L^1_k}$, $\mathcal
{M}_2=((S^{-1}_k)_{i,L^1_k})_{i\neq L^1_k}$,
$\mathcal{M}_3=(S^{-1}_k)_{L^1_k,L^1_k}$, $\tilde{\mathcal
{Z}}_2=x(\mathcal{I}_k)_{L^1_k}$,
$\tilde{\mathcal{Z}}_1=(x(\mathcal{I}_k)^{\star},y(\mathcal{J}_k)^{\star
})^{\star}\setminus\tilde{\mathcal{Z}}_2$,
$v_1=((\tilde{b}_k^{vu}(\tilde{b}_k^{vu})^{\star}\Delta\check
{u}^k)^{-1})_{11}$ and
\[
\mathcal{M}_4=
\pmatrix{
\mathcal{M}_3(v_1+\mathcal{M}_3)^{-1}
& -\bigl(\bigl(\tilde{b}_k^{vu}\bigl(\tilde
{b}_k^{vu}\bigr)^{\star}\Delta\check{u}^k
\bigr)^{-1}\bigr)_{12}(v_1+\mathcal
{M}_3)^{-1}
\cr\noalign{\vspace*{2pt}}
0 & 1 },
\]
where $x(\mathcal{I}_k)=\{x(\mathcal{I}_k)_i\}_{i=1}^{L^1_{k-h,k-1}}$
and $y(\mathcal{J}_k)=\{y(\mathcal{J}_k)_j\}_{j=1}^{L^2_{k-h,k-1}}$.
Then we obtain
%
\begin{eqnarray}
&&\bigl(x(\mathcal{I}_{k})^{\star},y(
\mathcal{J}_{k})^{\star
}\bigr)S_k^{-1}
\bigl(x(\mathcal{I}_{k})^{\star},y(\mathcal{J}_{k})^{\star}
\bigr)^{\star
}
\nonumber
\\
\label{integration-lemma-equ2}&&\quad =\tilde{\mathcal{Z}}_1^{\star}\mathcal{M}_1
\tilde{\mathcal {Z}}_1+2\tilde{\mathcal{Z}}_1^{\star}
\mathcal{M}_2\tilde{\mathcal {Z}}_2+\tilde{
\mathcal{Z}}_2^{\star}\mathcal{M}_3\tilde{
\mathcal{Z}}_2\nonumber
\\[-8pt]\\[-8pt]
&&\quad = \bigl(\tilde{\mathcal{Z}}_2+\mathcal{M}_3^{-1}
\mathcal{M}_2^{\star}\tilde {\mathcal{Z}}_1
\bigr)^{\star}\mathcal{M}_3\bigl(\tilde{\mathcal{Z}}_2+
\mathcal {M}_3^{-1}\mathcal{M}_2^{\star}
\tilde{\mathcal{Z}}_1\bigr)\nonumber\\
&&\qquad{}+\tilde{\mathcal {Z}}_1^{\star}
\mathcal{M}_1\tilde{\mathcal{Z}}_1 - \tilde{\mathcal
{Z}}_1^{\star}\mathcal{M}_2\mathcal{M}_3^{-1}
\mathcal{M}_2^{\star}\tilde {\mathcal{Z}}_1.\nonumber
\end{eqnarray}

By (\ref{integration-lemma-equ1}), (\ref{integration-lemma-equ2}) and Lemma~\ref{dist-comb-lemma} in the
\hyperref[app]{Appendix}, we have
\begin{eqnarray*}
&& \int\tilde{f}^{vu,(1)}_k\exp \biggl(\sum
_{k'\leq k}\acute {f}^{k,vu}_{k'} \biggr) (z)\,\mathrm{d}
\hat{z}_{1,k}\\
&& \quad  =\int \biggl\{-\tilde{\mathcal{Z}}_3^{\star}
\mathcal{M}_4^{\star
}\frac{\partial_{\sigma}((\tilde{b}_{(k)}\tilde{b}^{\star
}_{(k)})^{-1})}{2\Delta\check{u}^k} \mathcal{M}_4
\tilde{\mathcal{Z}}_3 +\Upsilon_6 \biggr\} \exp \biggl(
\sum_{k'\leq k}\acute{f}^{k,vu}_{k'}
\biggr) (z)\,\mathrm{d}\hat{z}_{1,k},
\end{eqnarray*}
where $k'_1=\max\mathcal{I}^-_{k-h,k-1}$, $\tilde{\mathcal
{Z}}_3=(x_k-x_{k'_1}+\mathcal{M}_3^{-1}\mathcal{M}_2^{\star}\tilde
{\mathcal{Z}}_1,y_k-y_{k-1})^{\star}$
and $\Upsilon_6= -\partial_{\sigma}((\tilde{b}_{(k)}\* \tilde{b}^{\star
}_{(k)})^{-1})_{11}(\Delta\check{u}^k)^{-1}(v_1+\mathcal{M}_3)^{-1}/2
-\partial_{\sigma}\log\det(\tilde{b}_{(k)}\tilde{b}^{\star}_{(k)})/2$.

Moreover, a similar argument yields
\begin{eqnarray*}
&&\int\tilde{f}^{vu,(1)}_k\exp \biggl(\sum
_{k'}\acute{f}^{k,vu}_{k'} \biggr) (z)\,\mathrm{d}
\hat{z}
\\
&&\quad =\int \biggl\{-\tilde{\mathcal{Z}}_5^{\star}\mathcal
{M}_4^{\star}\frac{\partial_{\sigma}((\tilde{b}_{(k)}\tilde{b}^{\star
}_{(k)})^{-1})}{2\Delta\check{u}^k} \mathcal{M}_4
\tilde{\mathcal{Z}}_5 +\tilde{\mathcal{Q}}^{k,v}_2(z_0,
\bar{z};\bar{u}) \biggr\} \exp \biggl(\sum_{k'}
\acute{f}^{k,vu}_{k'} \biggr) (z)\,\mathrm{d}\hat{z},
\end{eqnarray*}
where $\mathcal{P}=\{L^1_{k-h,k},L_{k-h,k}\}$,
$\mathcal{M}'_2=((S^{-1}_{k-h,k})_{ij})_{i\notin\mathcal{P},j\in
\mathcal{P}}$,
$\mathcal{M}'_3=((S^{-1}_{k-h,k})_{ij})_{i,j\in\mathcal{P}}$,
$\tilde{\mathcal{Z}}'_1=(x(\mathcal{I}_{k-h,k})^{\star},y(\mathcal
{J}_{k-h,k})^{\star})^{\star}\setminus(x(\mathcal
{I}_{k-h,k})_{L^1_{k-h,k}},y(\mathcal{J}_{k-h,k})_{L^2_{k-h,k}})^{\star}$, 
$\mathcal{M}_5=((S^{-1}_{k+1,k+h})_{i1})_{i\neq1}$,  $\mathcal
{M}_6 = (S^{-1}_{k+1,k+h})_{11}$, $\tilde{\mathcal{Z}}_4=(x(\mathcal
{I}_{k+1,k+h})^{\star},y(\mathcal{J}_{k+1,k+h})^{\star})^{\star
}\setminus x(\mathcal{I}_{k+1,k+h})_1$,
$k'_3=\min(\mathcal{I}_{k+1,k+h}\setminus\{k\})$,\vspace*{-1pt}
\begin{eqnarray*}
\tilde{\mathcal{Z}}_5&=&\pmatrix{
\mathcal{M}_6\bigl(\bigl(\mathcal{M}'_3
\bigr)_{11}+\mathcal{M}_6\bigr)^{-1} & -\bigl(
\mathcal {M}'_3\bigr)_{12}\bigl(\bigl(
\mathcal{M}'_3\bigr)_{11}+
\mathcal{M}_6\bigr)^{-1}
\cr\noalign{\vspace*{2pt}}
0 & 1} \\
&&{}\times\biggl\{ \pmatrix{
x_{k'_3}-x_{k'_1}+\mathcal{M}_6^{-1}
\mathcal{M}_5^{\star}\tilde{\mathcal {Z}}_4
\cr\noalign{\vspace*{2pt}}
y_k-y_{k-1}}
+\bigl(\mathcal{M}'_3\bigr)^{-1}\bigl(
\mathcal{M}'_2\bigr)^{\star}\tilde{\mathcal
{Z}}'_1 \biggr\} \\
&&{}+\bigl(\mathcal{M}_3^{-1}
\mathcal{M}_2^{\star}\tilde{\mathcal{Z}}_1,0
\bigr)^{\star
}-\bigl(\mathcal{M}'_3
\bigr)^{-1}\bigl(\mathcal{M}'_2
\bigr)^{\star}\tilde{\mathcal{Z}}'_1,
\nonumber
\end{eqnarray*}
and\vspace*{-1pt}
\[
\tilde{\mathcal{Q}}^{k,v}_2(z_0,\bar{z};
\bar{u})=\Upsilon_6-\bigl(\mathcal {M}_4^{\star}
\partial_{\sigma}\bigl(\bigl(\tilde{b}_{(k)}\tilde{b}^{\star
}_{(k)}
\bigr)^{-1}\bigr) \mathcal{M}_4\bigr)_{11}
\bigl((\mathcal{M}_3)_{11}+\mathcal {M}_6
\bigr)^{-1}/\bigl(2\Delta\check{u}^k\bigr).
\]

Let $\hat{\mathcal{Z}}_5$ be a vector obtained by substituting $0$
for $x(\mathcal{I}_k)_1,y(\mathcal{J}_k)_1,x(\mathcal
{I}_{k-h,k})_1,y(\mathcal{J}_{k-h,k})_1,\allowbreak x(\mathcal
{I}_{k+1,k+h})_{L^1_{k+1,k+h}}$ and $y(\mathcal
{J}_{k+1,k+h})_{L^2_{k+1,k+h}}$ in $\tilde{\mathcal{Z}}_5$.
Then since $\mathcal{M}_4$, $\hat{\mathcal{Z}}_5$ and $\tilde{\mathcal
{Q}}^{k,v}_2$ do not depend on $\hat{z}$, we obtain\vspace*{-1pt}
%
\begin{eqnarray}
\frac{\bar{\mathbb{P}}^{2,vu}_k(\tilde{f}^{vu,(1)}_k)}{\bar{\mathbb
{P}}^{2,vu}_k(1)}(z_0,\bar{z}) &=& -\hat{
\mathcal{Z}}_5^{\star}\mathcal{M}_4^{\star}
\frac
{\partial_{\sigma}((\tilde{b}_{(k)}\tilde{b}^{\star
}_{(k)})^{-1})}{2\Delta\check{u}^k} \mathcal{M}_4\hat{\mathcal{Z}}_5
\nonumber
\\[-8pt]
\label{integration-lemma-equ4}\\[-8pt]
\nonumber
&&{}+
\tilde{\mathcal{Q}}^{k,v}_2(z_0,\bar{z};\bar{u})
+ \tilde{\mathcal {Q}}^{k,v}_4(z_0,\bar{z};
\bar{u}),
\end{eqnarray}
where $\Upsilon_7=(\hat{\mathcal{Z}}_5+\tilde{\mathcal
{Z}}_5)^{\star}\mathcal{M}_4^{\star}\partial_{\sigma}((\tilde
{b}_{(k)}\tilde{b}^{\star}_{(k)})^{-1})\mathcal{M}_4(\hat{\mathcal
{Z}}_5-\tilde{\mathcal{Z}}_5)/(2\Delta\check{u}^k)$
and\vspace*{-1pt}
\[
\tilde{\mathcal{Q}}^{k,v}_4(z_0,\bar{z};
\bar{u})=\int\Upsilon_7\exp \biggl(\sum
_{k'}\acute{f}^{k,vu}_{k'}\biggr) (z)\,\mathrm{d}
\hat{z}\Big/\int\exp\biggl(\sum_{k'}\acute
{f}^{k,vu}_{k'}\biggr) (z)\,\mathrm{d}\hat{z}.
\]
\end{step}

\begin{step}
We will prove $\sup_{v,k}|\tilde{\mathcal
{Q}}^{k,v}_4(Y_{\Pi};\Pi)|=\mathrm{o}_p(b_n^{-q})$  for any $q>0$ in
this step.
We follow the approach in Section~2 of Ogihara and Yoshida \cite{ogi-yos}.

Let\vspace*{-1pt}
\begin{eqnarray*}
D'_k &=& \operatorname{diag} \biggl( \biggl(\sum
_{k'}\bigl|\check{b}^1_{k'}\bigr|^2\bigl|
\tilde {I}^{i_p}_{k-h,k-1}\cap\tilde{K}_{k'}\bigr|
\biggr)_{p=1}^{L^1_k}, \biggl(\sum_{k'}\bigl|
\check{b}^2_{k'}\bigr|^2\bigl|\tilde{J}^{j_q}_{k-h,k-1}
\cap\tilde {K}_{k'}\bigr| \biggr)_{q=1}^{L_k^2} \biggr),
\\
G_k&=&  \biggl\{\frac{|\tilde{I}^{i_p}_{k-h,k-1}\cap\tilde
{J}^{j_q}_{k-h,k-1}|}{|\tilde{I}^{i_p}_{k-h,k-1}|^{1/2}|\tilde{J}^{j_q}_{k-h, k-1}|^{1/2}} \biggr\}_{1\leq
p\leq L^1_k,1\leq q\leq L^2_k},\\
\tilde{G}_k &=&  \biggl\{\frac{\sum_{k'}\check{b}^1_{k'}\cdot
\check{b}^2_{k'}|\tilde{I}^{i_p}_{k-h,k-1}\cap\tilde
{J}^{j_q}_{k-h,k-1}\cap\tilde{K}_{k'}|}{
((D'_k)_{p,p})^{1/2}((D'_k)_{q+L^1_k,q+L^1_k})^{1/2}} \biggr\}_{1\leq
p\leq L^1_k, 1\leq q\leq L^2_k}.
\end{eqnarray*}
Then we obtain\vspace*{-1pt}
%
\begin{equation}
\label{S-inv-rep}
S_k^{-1}=\bigl(D'_k
\bigr)^{-{1}/{2}}\pmatrix{
 \bigl(
\mathcal{E}-\tilde{G}_k\tilde{G}_k^{\star}
\bigr)^{-1} & -\bigl(\mathcal{E}-\tilde {G}_k
\tilde{G}_k^{\star}\bigr)^{-1}\tilde{G}_k
\cr\noalign{\vspace*{2pt}}
-\tilde{G}_k^{\star}\bigl(\mathcal{E}-\tilde{G}_k
\tilde{G}_k^{\star}\bigr)^{-1} & \bigl(\mathcal{E}-
\tilde{G}_k^{\star}\tilde{G}_k\bigr)^{-1}
} \bigl(D'_k\bigr)^{-{1}/{2}},
\end{equation}
by a standard formula for block matrices.

Moreover, the argument in Lemma~2 of Ogihara and Yoshida \cite{ogi-yos} yields
%
\begin{equation}
\label{tildeG-est} \Vert\tilde{G}_k\Vert\vee\bigl\Vert
\tilde{G}_k^{\star} \bigr\Vert\leq\tilde{\rho}(\bar{u})\Vert
G_k \Vert\vee \bigl\Vert G_k^{\star} \bigr\Vert
\leq\tilde{\rho}(\bar{u}),
\end{equation}
where
\begin{eqnarray*}
\tilde{\rho}(\bar{u})& =&\tilde{\rho}(z_0,\bar{z},\bar{u})\\
 &=&\sup
\biggl\{\frac{|\check{b}^1_{k_1}\cdot\check
{b}^2_{k_2}|}{|\check{b}^1_{k_3}||\check{b}^2_{k_4}|};v\in[0,1],1\leq k\leq L_0(\bar{u})  \mbox{ and  there  exist }  l_1,l_2  \mbox{ such  that }
\\
&&\qquad \tilde{I}^{l_1}_{k-h,k+h}\cap\tilde{J}^{l_2}_{k-h,k+h}
\neq \varnothing, \tilde{K}_{k_1},\tilde{K}_{k_3}\subset
\tilde{I}^{l_1}_{k-h,k+h}  \mbox{ and }  \tilde{K}_{k_2},
\tilde{K}_{k_4}\subset\tilde{J}^{l_2}_{k-h,k+h} \biggr\}.
\end{eqnarray*}
Let $\tilde{l}_1(k;\bar{u})=\min\{l\in\mathbb{Z}_+;((G_kG_k^{\star
})^lG_k)_{1,L^2_k}>0\}$,
then (\ref{S-inv-rep}), (\ref{tildeG-est}) and relations $(\mathcal
{E}-\tilde{G}_k\tilde{G}_k^{\star})^{-1}=\sum_{l=0}^{\infty}(\tilde
{G}_k\tilde{G}_k^{\star})^l$ and $(\mathcal{E}-\tilde{G}_k^{\star}\tilde
{G}_k)^{-1}=\sum_{l=0}^{\infty}(\tilde{G}_k^{\star}\tilde{G}_k)^l$ yield
%
\begin{eqnarray}
&& \bigl|\bigl(D'_k\bigr)^{1/2}_{11}(
\mathcal{M}_2)_1\bigl(D'_k
\bigr)^{1/2}_{L^1_k,L^1_k}\bigr| \vee\bigl|\bigl(D'_k
\bigr)^{1/2}_{L^1_k+1,L^1_k+1}(\mathcal {M}_2)_{L^1_k}
\bigl(D'_k\bigr)^{1/2}_{L^1_k,L^1_k}\bigr|
\nonumber
\\[-8pt]
\label{integration-lemma-equ5}\\[-8pt]
\nonumber
&& \quad \leq C
\tilde{\rho}(\bar{u})^{2\tilde{l}_1-1}/\bigl(1-\tilde{\rho}(\bar{u})^2
\bigr)\nonumber
\end{eqnarray}
if $\tilde{\rho}(\bar{u})<1$.

On the other hand, we have
\[
\mathcal{M}_3^{-1}=(S_k)_{L^1_k,L^1_k} -
\bigl((S_k)_{i,L^1_k}\bigr)^{\star}_{i\neq L^1_k}
\bigl(\bigl((S_k)_{ij}\bigr)_{i,j\neq
L^1_k}
\bigr)^{-1} \bigl((S_k)_{i,L^1_k}\bigr)_{i\neq L^1_k}
\]
by a standard formula for block matrices,
and hence
%
\begin{equation}
\label{integration-lemma-equ6}
\bigl|\bigl(D'_k\bigr)^{-1/2}_{L^1_k,L^1_k}
\mathcal {M}_3^{-1}\bigl(D'_k
\bigr)^{-1/2}_{L^1_k,L^1_k}\bigr| \leq C\bigl(1-\tilde{\rho}(
\bar{u})^2\bigr)^{-1}
\end{equation}
if $\tilde{\rho}(\bar{u})<1$.

Moreover we have
\[
v_1^{-1}+\mathcal{M}_3^{-1}\geq
\bigl(1-\tilde{\rho}^2\bigr) \bigl(\bigl(\tilde{b}^{vu}_k
\bigl(\tilde{b}^{vu}_k\bigr)^{\star}\Delta\check
{u}^k\bigr)_{11}+(S_k)_{L^1_k,L^1_k}\bigr)
\geq\bigl(1-\tilde{\rho}^2\bigr) (S_{k-h,k})_{L^1_k,L^1_k},
\]
and consequently we obtain $(v_1+\mathcal{M}_3)^{-1}\leq C(1-\tilde{\rho
}^2)^{-1}\Delta\check{u}^k(\check{u}^{k-1}-\check{u}^{k'_1})/(\check
{u}^k-\check{u}^{k'_1})$.
Similarly we have $((\mathcal{M}'_3)_{11}+\mathcal{M}_6)^{-1}\leq
C(1-\tilde{\rho}^2)^{-1}(\check{u}^k-\check{u}^{k'_1})(\check
{u}^{k'_3}-\check{u}^k)/(\check{u}^{k'_3}-\check{u}^{k'_1})$.

Therefore, we obtain
%
\begin{equation}
\label{integration-lemma-equ8} \sup_{k,v} \bigl|\tilde{\mathcal{Q}}^{k,v}_4(Y_{\Pi};
\Pi) \bigr| \leq C\sup_k\frac{\tilde{\rho}(Y_{\Pi};\Pi)^{2(\tilde{l}_1(k;\Pi)\wedge
\tilde{l}_2(k;\Pi))-1}}{(1-\tilde{\rho}(Y_{\Pi};\Pi)^2)^6}\times \mathrm{O}_p
\bigl((\ell_{1,n}+\ell_{2,n})^3\bigr)
\end{equation}
on $\{\tilde{\rho}(Y_{\Pi};\Pi)<1\}$  by Lemma~\ref
{P2P1-tight}, (\ref{integration-lemma-equ5}), (\ref{integration-lemma-equ6})
and similar estimates for $\{\mathcal{M}_l\}_{2\leq l\leq6}$ and $\{
\mathcal{M}'_l\}_{l=2,3}$,\vspace*{-2pt} where
\[
G'_k= \biggl\{\frac{|\tilde{I}^{i_p}_{k+1,k+h}\cap\tilde
{J}^{j_q}_{k+1,k+h}|}{|\tilde{I}^{i_p}_{k+1,k+h}|^{1/2}|\tilde
{J}^{j_q}_{k+1,k+h}|^{1/2}} \biggr
\}_{1\leq p\leq L^1_{k+1,k+h},1\leq
q\leq L^2_{k+1,k+h}}
\]
and $\tilde{l}_2(k;\bar{u})=\min\{l\in\mathbb{Z}_+;((G'_k(G'_k)^{\star
})^lG'_k)_{1,L^2_{k+1,k+h}}>0\}$.

By the definitions of $\tilde{l}_1$, we\vspace*{-1pt} obtain
%
\begin{equation}
\label{tildeL-est1} \bigl(2\tilde{l}_1(k;\bar{u})+2\bigr)\max
_{i,j}\bigl(\bigl|s^i-s^{i-1}\bigr|\vee
\bigl|t^j-t^{j-1}\bigr|\bigr)\geq\bigl|\check{u}^{k-1}-
\check{u}^{k-h-1}\bigr|
\end{equation}
for any $k$.
Moreover, since the numbers of elements of $(s^i)_i\cap[\check
{u}^{k-h-1},\check{u}^{k-1}]$ or $(t^j)_j\cap[\check
{u}^{k-h-1},\check{u}^{k-1}]$
is equal to or greater than $(h+1)/2$, $[A2]$ and $[A3']$\vspace*{-1pt} yields
%
\begin{equation}
\label{tildeL-est2}
\liminf_{n\to\infty}P \Bigl[\inf_k\bigl|U^{k-1}-U^{k-h-1}\bigr|>b_n^{-1-\delta_3}
\bigl((h+1)/2-1\bigr) \Bigr] \geq \liminf_{n\to\infty}P[
\mathbf{A}_n] =1,
\end{equation}
where\vspace*{-2pt}
\begin{eqnarray*}
\mathbf{A}_n &=&  \bigcap_{|j_2-j_1|\geq b_n^{\delta_2}} \biggl[
\biggl\{\frac
{|S^{n,j_2}-S^{n,j_1}|}{|j_2-j_1|}>b_n^{-1-\delta_3}  \mbox{ if }
j_1\vee j_2 \leq\ell_{1,n} \biggr\}\\[-2pt]
&&\hspace*{45pt}{}\cap \biggl
\{\frac{|T^{n,j_2}-T^{n,j_1}|}{|j_2-j_1|}>b_n^{-1-\delta_3}  \mbox{ if }  j_1
\vee j_2 \leq\ell_{2,n} \biggr\} \biggr].
\end{eqnarray*}
By (\ref{tildeL-est1}) and (\ref{tildeL-est2}), we have\vspace*{-2pt}
%
\begin{equation}
\label{integration-lemma-equ9} \lim_{n\to\infty} P \Bigl[\inf_k
\tilde{l}_1(k;\Pi)>b_n^{\kappa-\delta
_1-\delta_3}/5 \Bigr]=1.
\end{equation}
Similarly, we obtain\vspace*{-2pt}
%
\begin{equation}
\label{integration-lemma-equ10} \lim_{n\to\infty} P \Bigl[\inf_k
\tilde{l}_2(k;\Pi)>b_n^{\kappa-\delta
_1-\delta_3}/5 \Bigr]=1.
\end{equation}

Let $\bar{\rho}=\sup_{t,x,y,\sigma}|b^1|^{-1}|b^2|^{-1}|b^1\cdot
b^2|(t,x,y,\sigma)$.
Then by virtue of $[A1']$ and the relation $\det(bb^{\star
})=|b^1|^2|b^2|^2-(b^1\cdot b^2)^2$, we obtain $\bar{\rho}<1$.
Moreover, we\vspace*{-2pt} obtain
%
\begin{equation}
\label{tildeRho-est} \lim_{n\to\infty}P\bigl[\tilde{\rho}(Y_{\Pi};
\Pi)>1-(1-\bar{\rho})/2\bigr]=0,
\end{equation}
since $r_n\to^p0$ and $b(t,x,y,\sigma)$ is continuous with respect to
$(t,x,y,\sigma)$.

By (\ref{integration-lemma-equ8}), (\ref{integration-lemma-equ9}), (\ref
{integration-lemma-equ10}) and (\ref{tildeRho-est}), we have
$\sup_{v,k}|\tilde{\mathcal{Q}}^{k,v}_4(Y_{\Pi};\Pi)|=\mathrm{o}_p(b_n^{-q})$
for any $q>0$.
Furthermore, we can\vspace*{-2pt} write
\[
\frac{\bar{\mathbb{P}}^{2,vu}_k(\tilde{f}^{vu,(1)}_k)}{\bar{\mathbb
{P}}^{2,vu}_k(1)}(z_0,\bar{z}) =\mathcal{Z}_k^{\star}
\tilde{\mathcal{Q}}^{k,v}_1(z_0,\bar{z};\bar
{u})\mathcal{Z}_k+\tilde{\mathcal{Q}}^{k,v}_2(z_0,
\bar{z};\bar {u})+\tilde{\mathcal{Q}}^{k,v}_4(z_0,
\bar{z};\bar{u}),
\]
where $\sup_{v,k}\Vert\tilde{Q}^{k,v}_1(z_0,\bar{z};\bar
{u})\Vert\leq C(1-\tilde{\rho}^2(\bar{u}))^{-6}$
and $\sup_{k,v}|\tilde{\mathcal{Q}}_2^{k,v}(z_0,\bar{z};\bar{u})|\leq
C(1-\tilde{\rho}(\bar{u})^2)^{-3}$.
\end{step}

\begin{step}
We now complete the proof.

Let $\mathcal{Q}^{k,v}_p$ be obtained by substituting the same values
in $\tilde{\mathcal{Q}}^{k,v}_p$ as
\[
\mathcal{Q}_p^{k,v}(z_0,\bar{z};\bar{u})=
\tilde{\mathcal {Q}}_p^{k,v}\bigl(z_0,
\bigl((x_{\hat{k}_1})_{i=1}^{L^1},(y_{\hat
{k}_2})_{j=1}^{L^2}
\bigr);\bar{u}\bigr) 
\]
for $p=1,2$, where $\hat{k}_1,\hat{k}_2$ are  the maximum
integers satisfying $\check{u}^{\hat{k}_1}=s^i$, $\check{u}^{\hat{k}_2}=t^j$,
$s^i\vee t^j\leq\inf(\tilde{\theta}(k-h-1,2;\bar{u}))$ for some $i$
and $j$.
Then we have
\[
\sup_{n,v,\bar{u},k,z_0,\bar{z}}\bigl(\bigl\Vert\mathcal{Q}^{k,v}_1(z_0,
\bar {z};\bar{u})\bigr\Vert\vee\bigl|\mathcal{Q}_2^{k,v}(z_0,
\bar{z};\bar {u})\bigr|\bigr)\leq C\bigl(1-\bar{\rho}^2
\bigr)^{-6}\leq C.
\]
Therefore, by setting
\begin{eqnarray*}
\mathcal{Q}^{k,v}_3&=&\bigl(\mathcal{Z}_k^{\star}
\bigl(\tilde{\mathcal {Q}}^{k,v}_1-\mathcal{Q}^{k,v}_1
\bigr)\mathcal{Z}_k+\tilde{\mathcal {Q}}^{k,v}_2-
\mathcal{Q}^{k,v}_2\bigr)1_{\{\tilde{\rho}(\bar{u})\leq
1-(1-\bar{\rho})/2\} },
\\
\mathcal{Q}^{k,v}_4&=&\tilde{\mathcal{Q}}^{k,v}_4+
\bigl(\mathcal{Z}_k^{\star
}\bigl(\tilde{\mathcal{Q}}^{k,v}_1-
\mathcal{Q}^{k,v}_1\bigr)\mathcal{Z}_k+\tilde {
\mathcal{Q}}_2^{k,v}-\mathcal{Q}_2^{k,v}
\bigr)1_{\{\tilde{\rho}(\bar{u})>
1-(1-\bar{\rho})/2\} },
\end{eqnarray*}
we obtain $\sup_{k,v}|\mathcal{Q}^{k,v}_4(Y_{\Pi};\Pi)|=\mathrm{o}_p(b_n^{-q})$,
$\sup_{n,r,k,v,v'}(r_n^{-1/2}b_n^{-3\kappa/2}E[|\mathcal{Q}^{k,v}_3(Y^{r,v'u}_{\Pi};\break\Pi)|^q| \Pi]^{1/q})\leq C$ a.s.
by (\ref{tildeRho-est}).

Furthermore, a similar argument for
\[
\int\tilde{f}^{vu,(1)}_k\exp\biggl(\sum
_{k'}\tilde{f}^{vu}_{k'}\biggr)\,\mathrm{d}
\hat{z}\Big/\int \exp\biggl(\sum_{k'}
\tilde{f}^{vu}_{k'}\biggr)\,\mathrm{d}\hat{z}
\]
yields
\begin{eqnarray*}
&& \sup_{k,v} \biggl|\frac{\bar{\mathbb{P}}^{2,vu}_k(\tilde
{f}_k^{vu,(1)})}{\bar{\mathbb{P}}^{2,vu}_k(1)} -\frac{\int\tilde{f}^{vu,(1)}_k\exp(\sum_{k'}\tilde{f}^{vu}_{k'})\,\mathrm{d}\hat
{z}}{\int\exp(\sum_{k'}\tilde{f}^{vu}_{k'})\,\mathrm{d}\hat{z}}
\biggr|(Y_{\Pi})\\
&&\quad  \leq \mathrm{O}_p\bigl((\ell_{1,n}+
\ell_{2,n})^3\cdot\tilde{\rho}(\Pi )^{2\inf_k(\tilde{l}_1\wedge\tilde{l}_2(k;\Pi))-1}
\bigr)+\mathrm{o}_p\bigl(b_n^{-q}\bigr)=\mathrm{o}_p
\bigl(b_n^{-q}\bigr)
\end{eqnarray*}
for any $q>0$.
\end{step}
\upqed\end{pf}

The following lemma enables us to replace $\bar{\mathbb{P}}^0_u$ and
$\bar{\mathbb{P}}^0_0$ in $\log(\bar{\mathbb{P}}^0_u/\bar{\mathbb{P}}^0_0)$
by the function $\int\exp(\sum_k\tilde{f}_k^u)\,\mathrm{d}\hat{z}$ and $\int\exp
(\sum_k\tilde{f}_k^0)\,\mathrm{d}\hat{z}$, respectively.

%
\begin{lemma}\label{coeff-move-lemma2}
Let $u\in\mathbb{R}^d$. Assume $[A1'],[A2]$ and $[A3']$. Then
\[
\log\frac{\bar{\mathbb{P}}^0_u}{\bar{\mathbb{P}}^0_0}(Y_{\Pi
})-b_n^{-{1}/{2}}u\int
^1_0\sum_k
\frac{\bar{\mathbb
{P}}^{2,vu}_k(\tilde{f}^{vu,(1)}_k)}{\bar{\mathbb
{P}}^{2,vu}_k(1)}\,\mathrm{d}v(Y_{\Pi})\to^p 0
\]
as $n\to\infty$.
\end{lemma}

\begin{pf}
Let
\[
\mathcal{A}^5_{k,v}(z_0,\bar{z}) =
\mathcal{Z}_k^{\star}\mathcal{Q}_1^{k,v}(z_0,
\bar{z},\bar{u})\mathcal {Z}_k1_{\{\mathcal{Z}_k\neq\varnothing\} }+\mathcal{Q}^{k,v}_2(z_0,
\bar {z},\bar{u}) +\mathcal{Q}^{k,v}_3(z_0,
\bar{z},\bar{u})
\]
and
\[
\mathcal{K}^3_M(\bar{u})=\mathcal{K}^2_M(
\bar{u})\cap\Bigl\{(z_0,\bar {z});\sup_{k,v}\bigl|
\mathcal{Q}^{k,v}_4(z_0,\bar{z};\bar{u})\bigr|\leq
b_n^{-1}\Bigr\}
\]
for $M>0$ and $\bar{u}\in\mathcal{U}$.

By Lemmas \ref{coeff-move-lemma1}, \ref{L2M-est} and \ref
{integration-lemma}  and the definition of $\mathcal{K}^2_M$,
it is sufficient to show that
\[
\Phi'_n=\sup_vE \biggl[
\biggl|b_n^{-{1}/{2}}\sum_k \biggl(
\frac{\bar
{\mathbb{P}}^{2,vu}_k(\tilde{f}^{vu,(1)}_k1_{B^n_M})}{\tilde{\mathbb
{P}}^0_{M,vu}}-\frac{\bar{\mathbb{P}}^{2,vu}_k(\tilde
{f}^{vu,(1)}_k1_{B^n_M})}{\bar{\mathbb{P}}^{2,vu}_k(1_{B^n_M})} \biggr)1_{\mathcal{K}^3_M(\Pi)}
\biggr|(Y_{\Pi}) \Big|\Pi \biggr]\to^p 0
\]
as $n\to\infty$ for any $M>0$.

Fix $M>0$. Then Lemma~\ref{integration-lemma} and the definition of
$\mathcal{K}^3_M$ yield
\begin{eqnarray*}
\Phi'_n&=&\sup_vE \biggl[
\biggl|b_n^{-{1}/{2}}\sum_k
\frac{\bar{\mathbb
{P}}^{2,vu}_k(\tilde{f}^{vu,(1)}_k1_{B^n_M})}{\bar{\mathbb
{P}}^{2,vu}_k(1_{B^n_M})} \biggl(\frac{\bar{\mathbb
{P}}^{2,vu}_k(1_{B^n_M})}{\tilde{\mathbb{P}}^0_{M,vu}}-1 \biggr)1_{\mathcal{K}^3_M(\Pi)}
\biggr|(Y_{\Pi}) \Big|\Pi \biggr]
\\
&\leq&\sup_vE \biggl[ \biggl|b_n^{-{1}/{2}}
\sum_k\mathcal {A}^5_{k,v}(Y_{\Pi})
\biggl(\frac{\bar{\mathbb
{P}}^{2,vu}_k(1_{B^n_M})}{\tilde{\mathbb{P}}^0_{M,vu}}-1 \biggr)1_{\mathcal{K}^3_M(\Pi)} \biggr|(Y_{\Pi}) \Big|\Pi
\biggr]
\\
&&{}+2\sup_vE \biggl[b_n^{-{3}/{2}}\sum
_k \biggl|\frac{\bar{\mathbb
{P}}^{2,vu}_k(1_{B^n_M})}{\tilde{\mathbb{P}}^0_{M,vu}}-1 \biggr|1_{\mathcal{K}^3_M(\Pi)}(Y_{\Pi})
\Big|\Pi \biggr].
\end{eqnarray*}

The second term of the right-hand side in the above inequality is equal
to or smaller than
\[
Cb_n^{-3/2}\sup_v\sum
_k\int\bigl|\mathbb{P}^{2,vu}_k1_{B^n_M}-
\mathbb{P}^0_{vu}1_{B^n_M}\bigr|\,\mathrm{d}\hat{z}\,\mathrm{d}
\bar{z}P_{Y_0}(\mathrm{d}z_0) \bigg|_{\bar{u}=\Pi
}=\mathrm{o}_p(1).
\]
Hence, by a similar argument to the proof of Lemma~\ref{coeff-move-lemma1}, we obtain
\begin{eqnarray*}
\Phi'_n &\leq& \mathrm{e}^{M+1}\\
&&{}\times\sup
_v E \biggl[ \biggl|b_n^{-{1}/{2}}\sum
_k \mathcal{A}^5_{k,v}
\bigl(Y^{0,vu}_{\Pi}\bigr)
\biggl\{\exp \biggl(\sum
_{k';|k'-k|\leq
h}\bigl(\tilde{f}^{vu}_{k'}-\log
\check{p}^0_{k',vu}\bigr) \biggr)-1 \biggr\}\\
&&\hspace*{42pt}{}\times 1_{B^n_M}
\biggr|\bigl(Y^{0,vu}_{\check{U}}\bigr) \Big|\Pi \biggr]\\
&&{}+\mathrm{o}_p(1)
\\
& \leq & C\sup_v E \biggl[ \biggl|b_n^{-{1}/{2}}
\sum_k\mathcal {A}^5_{k,v}
\bigl(Y^{0,vu}_{\Pi}\bigr)\sum_{k';|k'-k|\leq h}
\bigl(\tilde {f}^{vu}_{k'}-\log\check{p}^0_{k',vu}
\bigr)1_{B^n_M} \biggr|\bigl(Y^{0,vu}_{\check
{U}}\bigr) \Big|\Pi
\biggr] \\
&&{}+\mathrm{O}_p\bigl(b_n^{-{1}/{2}}b_nb_n^{\kappa}b_n^{2\kappa
}b_n^{-{2}/{3}-2\kappa}
\bigr)+\mathrm{o}_p(1)
\\
&\leq& C\sup_vE \biggl[ \biggl|b_n^{-{1}/{2}}
\sum_k\mathcal {A}^5_{k,v}
\bigl(Y^{0,vu}_{\Pi}\bigr)\sum_{k';|k'-k|\leq h}
\bigl(\log\check {p}^1_{k',vu}-\log\check{p}^0_{k',vu}
\bigr) \biggr|\bigl(Y^{0,vu}_{\check{U}}\bigr) \Big|\Pi \biggr]
\\
&&{}+C\sup_vE \biggl[ \biggl|b_n^{-{1}/{2}}\sum
_k\mathcal {A}^5_{k,v}
\bigl(Y^{0,vu}_{\Pi}\bigr)\sum_{k';|k'-k|\leq h}
\biggl(\sum_{\tilde
{k}}\hat{\mathcal{A}}^3_{k,k',\tilde{k}}+
\mathcal{A}^4_{k,k'} \biggr) \biggr| \Big|\Pi \biggr]
+\mathrm{o}_p(1)\\
&=& \Phi'_{1,n}+\Phi'_{2,n}+\mathrm{o}_p(1).
\end{eqnarray*}

Furthermore, let $\tilde{\mathcal{A}}^5_{k,v}(z_0,\bar{z})=\mathcal
{Z}_k^{\star}\mathcal{Q}_1^{k,v}(z_0,\bar{z},\bar{u})\mathcal{Z}_k1_{\{
\mathcal{Z}_k\neq\varnothing\} }$,
then Lemma~\ref{integration-lemma} and the Burkholder--Davis--Gundy
inequality yield
\[
\Phi'_{1,n}\leq\sup_{r,v} E \biggl[
\biggl|b_n^{-{1}/{2}}\sum_k\tilde {
\mathcal{A}}^5_{k,v}\bigl(Y^{0,vu}_{\Pi}
\bigr)\sum_{k';|k'-k|\leq h}\frac
{\partial_r \check{p}^r_{k',vu}}{\check{p}^r_{k',vu}}
\bigl(Y^{0,vu}_{\check
{U}}\bigr) \biggr| \Big|\Pi \biggr]+\mathrm{o}_p(1),
\]
and
\[
\Phi'_{2,n}\leq C\sup_vE \biggl[
\biggl|b_n^{-{1}/{2}}\sum_k\tilde {
\mathcal{A}}^5_{k,v}\bigl(Y^{0,vu}_{\Pi}
\bigr)\sum_{k';|k'-k|\leq h} \biggl(\sum
_{\tilde{k}}\hat{\mathcal{A}}^3_{k,k',\tilde{k}}+
\mathcal {A}^4_{k,k'} \biggr) \biggr| \Big|\Pi \biggr]
+\mathrm{o}_p(1).
\]

Let $\acute{\mathcal{Z}}^j=\{\acute{\mathcal{Z}}^j_k\}_{k=1}^{L_0(\bar
{u})}$  $(j=1,2)$,
\[
\acute{\mathcal{Z}}^1_k=\frac{Y^{0,vu,1}_{\check
{U}^k}-Y^{0,vu,1}_{\check{U}^{k-1}}}{|\tilde{\theta
}_{0,k}|^{1/2}}, \qquad\acute{\mathcal{Z}}^2_k=\frac{Y^{0,vu,2}_{\check
{U}^k}-Y^{0,vu,2}_{\check{U}^{k-1}}}{|\inf\{T^j;T^j\geq\check{U}^k\}
-\sup\{T^j;T^j\leq\check{U}^{k-1}\}|^{1/2}},
\]
and $\check{\mathcal{Q}}^{k,v,j_1,j_2}=\{(\check{\mathcal
{Q}}^{k,v,j_1,j_2})_{l_1,l_2}\}_{l_1,l_2}$ be a certain symmetric\vspace*{1pt}
matrix $(1\leq j_1,j_2\leq2)$ satisfying
$\mathcal{Z}_k^{\star}\mathcal{Q}^{k,v}_1\mathcal{Z}_k1_{\{\mathcal
{Z}_k\neq\varnothing\} }=\sum_{j_1,j_2=1}^2(\acute{\mathcal
{Z}}^{j_1})^{\star}\check{\mathcal{Q}}^{k,v,j_1,j_2}\acute{\mathcal{Z}}^{j_2}$.
Moreover, let
\begin{eqnarray*}
\tilde{\mathcal{X}}^1_{l,k}&=&\sum
_{j_1,j_2=1}^2\bigl(\check{\mathcal {Q}}^{l,v,j_1,j_2}
\bigr)_{k,k}\acute{\mathcal{Z}}^{j_1}_k\acute{\mathcal{Z}}^{j_2}_k +2\sum_{l_2<k}
\sum_{j_1,j_2=1}^2\bigl(\check{\mathcal
{Q}}^{l,v,j_1,j_2}\bigr)_{k,l_2}\acute{\mathcal{Z}}^{j_1}_k
\acute{\mathcal {Z}}^{j_2}_{l_2},
\\
\tilde{\mathcal{X}}^2_{l,k',\tilde{k}}&=& \hat{\mathcal
{A}}^3_{l,k',\tilde{k}}+\mathcal{A}^4_{l,k'}1_{\{k'=\tilde{k}\} }.
\end{eqnarray*}

Then we have
\begin{eqnarray*}
\Phi'_{1,n}&\leq&C\sup_{r,v}E \biggl[
\biggl|b_n^{-{1}/{2}}\sum_{\tilde
{k}}\sum
_{k;|\tilde{k}-k|\leq h} \biggl(\sum_{k';|k'-k|\leq h,k'\leq
\tilde{k}}
\frac{\partial_r\check{p}^r_{k',vu}}{\check{p}^r_{k',vu}} \biggr)\tilde{\mathcal{X}}^1_{k,\tilde{k}}
\\
&&\hspace*{20pt}\qquad {}+b_n^{-{1}/{2}}\sum_{k'}
\biggl(\sum_{k;|k'-k|\leq h}\sum_{\tilde
{k};|\tilde{k}-k|\leq h,\tilde{k}<k'}
\tilde{\mathcal{X}}^1_{k,\tilde
{k}} \biggr)\frac{\partial_r\check{p}^r_{k',vu}}{\check{p}^r_{k',vu}} \biggr|
\bigl(Y^{0,vu}_{\check{U}}\bigr) \Big|\Pi \biggr]+\mathrm{o}_p(1)\\
&=& \mathrm{o}_p(1),
\end{eqnarray*}
by a similar argument to the estimate of $\Phi_{n,2}$ in the proof of
Lemma~\ref{coeff-move-lemma1}.

Moveover, 4. of Lemma~\ref{tech-lemma} in the \hyperref[app]{Appendix}
and estimates in the proof of Lemma~\ref{coeff-move-lemma1} yield
\[
\Phi'_{2,n}\leq C\sup_vE \biggl[
\biggl|b_n^{-{1}/{2}}\sum_{l,k}\tilde {
\mathcal{X}}^1_{l,k}\sum_{k';|k'-l|\leq h}
\tilde{\mathcal {X}}^2_{l,k',k} \biggr| \Big|\Pi \biggr]+\mathrm{o}_p(1).
\]

Therefore, we can see $\Phi'_{2,n}\to^p 0$ by using the
Burkholder--Davis--Gundy inequality
and estimates in the proof of Lemma~\ref{coeff-move-lemma1}.
\end{pf}

\begin{pf*}{Proof of Theorem~\ref{main}}
By virtue of (\ref{ratio-equ1}), Theorem~\ref{ogi-yos-thm}, Lemmas \ref
{A1toA1} and \ref{coeff-move-lemma2}, it is sufficient to show that
\[
b_n^{-{1}/{2}}u\int^1_0\sum
_k\frac{\bar{\mathbb{P}}^{2,vu}_k(\tilde
{f}^{vu,(1)}_k)}{\bar{\mathbb{P}}^{2,vu}_k(1)}\,\mathrm{d}v(Y_{\Pi})-
\bigl(H_n(\sigma _u)-H_n(
\sigma_{\ast})\bigr)\circ(\Pi, Y_{\Pi})\to^p 0
\]
as $n\to\infty$ under $[A1'],[A2]$ and $[A3']$.

Lemmas \ref{integration-lemma} and \ref{induction-lemma} yield
\begin{eqnarray*}
&& b_n^{-{1}/{2}}u\int^1_0\sum
_k\frac{\bar{\mathbb{P}}^{2,vu}_k(\tilde
{f}^{vu,(1)}_k)}{\bar{\mathbb{P}}^{2,vu}_k(1)}\,\mathrm{d}v(Y_{\Pi})
\\
&& \quad =b_n^{-{1}/{2}}u\int^1_0
\sum_k\frac{\int\tilde{f}^{vu,(1)}_k\exp
(\sum_{k'}\tilde{f}^{vu}_{k'})\,\mathrm{d}\hat{z}}{\int\exp(\sum_{k'}\tilde
{f}^{vu}_{k'})\,\mathrm{d}\hat{z}}\,\mathrm{d}v(Y_{\Pi})+\mathrm{o}_p(1)
\\
&& \quad =\int^1_0\partial_{v} \biggl(-
\frac{1}{2}\tilde{Z}^{\star}S_{1,L_0(\Pi
)}^{-1}
\tilde{Z}-\frac{1}{2}\log\det S_{1,L_0(\Pi)} \biggr)\,\mathrm{d}v+\mathrm{o}_p(1),
\end{eqnarray*}
where $\tilde{Z}=((Y^1_{S^{n,i}}-Y^1_{S^{n,i-1}})_i^{\star
},(Y^2_{T^{n,j}}-Y^2_{T^{n,j-1}})_j^{\star})^{\star}$.

Let $\tilde{D}=\operatorname{diag}((|I^i|)_i,(|J^j|)_j)$, then the difference
between $\tilde{D}^{-1/2}S_{1,L_0(\Pi)}\tilde{D}^{-1/2}$ and $S(\sigma
^n_{vu})$ in (\ref{S-def})
is only the substituted values of $b^1,b^2$.
Then we can see the right-hand side of the above equation is equal to
\[
\int^1_0\partial_vH_n
\bigl(\sigma^n_{vu}\bigr)\,\mathrm{d}v\circ(\Pi, Y_{\Pi
})+\mathrm{o}_p(1)=
\bigl(H_n\bigl(\sigma^n_u
\bigr)-H_n(\sigma_{\ast})\bigr)\circ(\Pi,
Y_{\Pi})+\mathrm{o}_p(1),
\]
by $[A2]$, Lemma~\ref{doubleThetaSum-est} and a similar argument to the
proof of Lemma~13 in Ogihara and Yoshida \cite{ogi-yos}. We omit the
details.
\end{pf*}

\begin{appendix}

\section*{Appendix}\label{app}

\setcounter{lemma}{0}
\begin{lemma}\label{tech-lemma}
Let $g,L\in\mathbb{N}$, $(\tilde{\Omega},\tilde{\mathcal{F}},\tilde
{P})$ be a probability space,
$\tilde{\mathbf{F}}=\{\tilde{\mathcal{F}}_k\}_{k=0}^L$ be a filtration.
Denote by $\tilde{E}$ the integral with respect to $\tilde{P}$.
\begin{enumerate}[4.]
\item[1.] Let $\{\mathcal{X}_{k,k'}\}_{1\leq k,k'\leq L}$ be random variables.\vspace*{1pt}
Suppose $\{\sum_{1\leq k'\leq l}\mathcal{X}_{k,k'}\}_{l=0}^L$
is $\tilde{\mathbf{F}}$-martingale for $1\leq k\leq L$.
Moreover, assume that there exists a sequence $\{\mathcal{C}_{k,k'}\}
_{1\leq k'\leq L}$ of positive numbers such that
$\tilde{E}[|\mathcal{X}_{k,k'}|^2]^{1/2}\leq\mathcal{C}_{k,k'}$
for $1\leq k,k'\leq L$. Then\vspace*{-2pt}
\[
\tilde{E} \Biggl[ \Biggl|\sum_{k=1}^L\sum
_{k'=1}^L\mathcal{X}_{k,k'} \Biggr|
\Biggr]\leq \Biggl(\sum_{l_1,l_2=1}^L\sum
_{k=1}^L\mathcal {C}_{l_1,k}
\mathcal{C}_{l_2,k} \Biggr)^{{1}/{2}}.
\]

\item[2.] Let $\{\mathcal{X}^1_k\}_{k=1}^L,\{\mathcal{X}^2_{k,k'}\}_{1\leq
k,k'\leq L}$ be random variables.
Suppose $\{\sum_{1\leq k\leq l}\mathcal{X}^1_k\}_{l=0}^L$ is
$\tilde{\mathbf{F}}$-martingale,
$\{\sum_{1\leq k'\leq l}\mathcal{X}^2_{k,k'}\}_{l=0}^L$ is
$\tilde{\mathbf{F}}$-martingale for $1\leq k\leq L$.
Moreover, assume that there exist a positive constant $\mathcal{C}^1$
and a sequence $\{\mathcal{C}^2_k\}_{1\leq k\leq L}$ of positive
numbers such that
$\tilde{E}[|\mathcal{X}^1_k|^4]^{1/4}\leq\mathcal{C}^1$ and $\tilde
{E}[|\mathcal{X}^2_{k,k'}|^4]^{1/4}\leq\mathcal{C}^2_{k'}$
for $1\leq k,k'\leq L$. Then\vspace*{-2pt}
\[
\tilde{E} \Biggl[ \Biggl|\sum_{k=1}^L
\mathcal{X}^1_k\sum_{k';0<|k'-k|\leq
g}
\mathcal{X}^2_{k,k'} \Biggr| \Biggr]\leq2(2g+1)
\mathcal{C}^1 \biggl\{\sum_{l_1,l_2;|l_1-l_2|\leq2g}
\mathcal{C}^2_{l_1}\mathcal {C}^2_{l_2}
\biggr\}^{{1}/{2}}. 
\]

\item[3.] Let $\{\mathcal{X}^1_k\}_{k=1}^L,\{\mathcal{X}^2_{k,k',\tilde
{k}}\}_{1\leq k,k',\tilde{k}\leq L}$ be random variables.
Suppose $\{\sum_{1\leq k\leq l}\mathcal{X}^1_k\}_{l=0}^L$ is
$\tilde{\mathbf{F}}$-martingale,
$\{\sum_{1\leq\tilde{k}\leq l}\mathcal{X}^2_{k,k',\tilde
{k}}\}_{l=0}^L$ is $\tilde{\mathbf{F}}$-martingale for $1\leq k,k'\leq L$.
Moreover, assume that there exist a positive constant $\mathcal{C}^1$
and sequences $\{\mathcal{C}^2_k\}_{1\leq k\leq L}$\break and $\{\mathcal
{C}^3_{k',\tilde{k}}\}_{1\leq k',\tilde{k}\leq L}$ of positive numbers
such that
$\tilde{E}[|\mathcal{X}^1_k|^4]^{1/4}\leq\mathcal{C}^1$,\break $\tilde
{E}[|\sum_{\tilde{k};\tilde{k}<l}\mathcal{X}^2_{k,k',\tilde
{k}}|^4]^{1/4}\leq\mathcal{C}^2_{k'}$
and $\tilde{E}[|\mathcal{X}^2_{k,k',\tilde{k}}|^4]^{1/4}\leq\mathcal
{C}^3_{k',\tilde{k}}$ for $1\leq l,k,k',\tilde{k}\leq L$.
Then\vspace*{-2pt}
\begin{eqnarray*}
&& \tilde{E} \Biggl[ \Biggl|\sum_{k=1}^L
\mathcal{X}^1_k\sum_{k';|k'-k|\leq
g}\sum
_{\tilde{k};\tilde{k}\neq k}\mathcal{X}^2_{k,k',\tilde{k}} \Biggr|
\Biggr] \\[-2pt]
&& \quad \leq\sqrt{2}(2g+1)\mathcal{C}^1 \biggl\{\sum
_{l_1,l_2;|l_1-l_2|\leq
2g}\mathcal{C}^2_{l_1}
\mathcal{C}^2_{l_2} +\sum_{l_1,l_2}
\sum_{k}\mathcal{C}^3_{l_1,k}
\mathcal{C}^3_{l_2,k} \biggr\} ^{{1}/{2}}.
\end{eqnarray*}

\item[4.] Let $\{\mathcal{X}^1_{l,k}\}_{1\leq k,l\leq L},\{\mathcal
{X}^2_{l,k',\tilde{k}}\}_{1\leq l,k',\tilde{k}\leq L}$ be random variables.
Suppose $\{\sum_{1\leq k\leq p}\mathcal{X}^1_{l,k}\}_{p=0}^L$
and $\{\sum_{1\leq\tilde{k}\leq p}\mathcal{X}^2_{l,k',\tilde
{k}}\}_{p=0}^L$ is $\tilde{\mathbf{F}}$-martingale for $1\leq l,k'\leq L$.
Moreover, assume that there  exist sequences $\{\mathcal{C}^1_{l,k}\}
_{1\leq l,k\leq L}$, $\{\mathcal{C}^2_{k'}\}_{1\leq k'\leq L}$ and $\{
\mathcal{C}^3_{k',\tilde{k}}\}_{1\leq k',\tilde{k}\leq L}$ of positive
numbers such that
$\tilde{E}[|\mathcal{X}^1_{l,k}|^4]^{1/4}\leq\mathcal{C}^1_{l,k}$,
$\tilde{E}[|\sum_{\tilde{k};\tilde{k}<k}\mathcal{X}^2_{l,k',\tilde
{k}}|^4]^{1/4}\leq\mathcal{C}^2_{k'}$
and $\tilde{E}[|\mathcal{X}^2_{l,k',\tilde{k}}|^4]^{1/4}\leq\mathcal
{C}^3_{k',\tilde{k}}$ for $1\leq l,k,k',\tilde{k}\leq L$.\vspace*{-2pt}
Then
\begin{eqnarray*}
&&\tilde{E} \biggl[ \biggl|\sum_{k,l}
\mathcal{X}^1_{l,k}\sum_{k';|k'-l|\leq
g}
\sum_{\tilde{k};\tilde{k}\neq k}\mathcal{X}^2_{l,k',\tilde{k}} \biggr|
\biggr]
\\[-2pt]
&& \quad \leq\sqrt{2} \biggl\{\sum_{l_1} \biggl(\sum
_l\mathcal{C}^1_{l,l_1}\sum
_{k';|k'-l|\leq g}\mathcal{C}^2_{k'}
\biggr)^2 +\sum_{l_1} \biggl(\sum
_l\sum_{l_2<l_1}
\mathcal{C}^1_{l,l_2}\sum_{k';|k'-l|\leq g}
\mathcal{C}^3_{k',l_1} \biggr)^2 \biggr
\}^{{1}/{2}}.
\end{eqnarray*}
\end{enumerate}
\end{lemma}

\begin{pf}
We first prove 4.
By using the Cauchy--Schwarz inequality and Lemma~9 in Ogihara and
Yoshida \cite{ogi-yos} repeatedly, we obtain
\begin{eqnarray*}
&&\tilde{E} \biggl[ \biggl|\sum_{k,l}
\mathcal{X}^1_{l,k}\sum_{k';|k'-l|\leq
g}
\sum_{\tilde{k};\tilde{k}\neq k}\mathcal{X}^2_{l,k',\tilde{k}} \biggr|
\biggr]^2
\\
&&\quad \leq\tilde{E} \biggl[ \biggl|\sum_{l_1,l}\sum
_{l_2<l_1} \biggl(\mathcal {X}^1_{l,l_1}\sum
_{k';|k'-l|\leq g}\mathcal{X}^2_{l,k',l_2}+
\mathcal {X}^1_{l,l_2}\sum_{k';|k'-l|\leq g}
\mathcal{X}^2_{l,k',l_1} \biggr) \biggr|^2 \biggr]
\\
&& \quad \leq\sum_{l_1} \biggl\{2\tilde{E} \biggl[ \biggl(
\sum_{l}\mathcal {X}^1_{l,l_1}
\sum_{l_2<l_1}\sum_{k';|k'-l|\leq g}
\mathcal {X}^2_{l,k',l_2} \biggr)^2 \biggr]\\
&&\quad \hspace*{28pt} {}+2
\tilde{E} \biggl[ \biggl(\sum_l\sum
_{l_2<l_1}\mathcal{X}^1_{l,l_2}\sum
_{k';|k'-l|\leq g}\mathcal{X}^2_{l,k',l_1}
\biggr)^2 \biggr] \biggr\}
\\
&& \quad \leq 2\sum_{l_1} \biggl(\sum
_l\mathcal{C}^1_{l,l_1}\tilde{E} \biggl[
\biggl(\sum_{l_2<l_1}\sum_{k';|k'-l|\leq g}
\mathcal{X}^2_{l,k',l_2} \biggr)^4
\biggr]^{{1}/{4}} \biggr)^2\\
&&\qquad {}+2\sum_{l_1}
\biggl(\sum_{l}\sum_{l_2<l_1}
\tilde{E} \biggl[ \biggl(\mathcal {X}^1_{l,l_2}\sum
_{k';|k'-l|\leq g}\mathcal{X}^2_{l,k',l_1}
\biggr)^2 \biggr]^{{1}/{2}} \biggr)^2
\\
&& \quad \leq 2\sum_{l_1} \biggl(\sum
_l\mathcal{C}^1_{l,l_1}\sum
_{k';|k'-l|\leq
g}\mathcal{C}^2_{k'}
\biggr)^2 +2\sum_{l_1} \biggl(\sum
_l\sum_{l_2<l_1}
\mathcal{C}^1_{l,l_2}\sum_{k';|k'-l|\leq g}
\mathcal{C}^3_{k',l_1} \biggr)^2.
\end{eqnarray*}
Then we obtain 4.

We obtain 3. by setting $\mathcal{X}^1_{l,k}=\mathcal{X}^1_k1_{\{l=k\}
}$ in 4.
We can prove 2. by setting $\mathcal{X}^2_{k,k',\tilde{k}}=\mathcal
{X}^2_{k,k'}1_{\{k'=\tilde{k}\} }$ in 3.
Moreover, we can easily check 1.
\end{pf}

The following  lemma is  proved by elementary calculation. We
omit proofs.

\begin{lemma}\label{dist-comb-lemma}
Let $A=\{A_{ij}\}_{i,j=1}^2$ be a $2\times2$ symmetric matrix, $V_1$
be a $2\times2$ symmetric, positive definite matrix, $\alpha,\beta\in
\mathbb{R}$ and $v_2>0$. Then
\begin{eqnarray*}
&&\int_{\mathbb{R}} (x_1-y_1+
\alpha,x_2-y_2+\beta)A(x_1-y_1+
\alpha ,x_2-y_2+\beta)^{\star}\\
&&\hspace*{1pt}\quad {}\times \varphi
\bigl((x_1-y_1,x_2-y_2)^{\star};V_1
\bigr)\varphi (y_1-w;v_2)\,\mathrm{d}y_1
\\
&& \quad = \biggl\{(\mathcal{W},x_2-y_2+\beta)A(
\mathcal{W},x_2-y_2+\beta )^{\star}+
\frac{A_{11}}{(V_1^{-1})_{11}+v_2^{-1}} \biggr\} \\
&&\qquad {}\times \int_{\mathbb{R}}\varphi
\bigl((x_1-y_1,x_2-y_2)^{\star};V_1
\bigr)\varphi (y_1-w;v_2)\,\mathrm{d}y_1
\end{eqnarray*}
for $x_1,x_2,y_2,w\in\mathbb{R}$,
where $(V_1^{-1})_{ij}$ denotes the element of $V_1^{-1}$ and
$\mathcal
{W}=(v_2^{-1}(x_1-w)-(V_1^{-1})_{12}(x_2-y_2))/((V_1^{-1})_{11}+v_2^{-1})+\alpha
$.
\end{lemma}

\begin{pf*}{Proof of Lemma~\ref{B1toA2}}
Let $\delta\in(3/q, \delta_2\wedge\delta_3)$ and
\[
\mathsf{A}=\bigcap_{i=1}^2\bigcap_{k=1}^{[Tb_n]+1}
\bigl\{\mathbf{N}^i_{(b_n^{-1}k)\wedge
T}-\mathbf{N}^i_{b_n^{-1}(k-1)}
\leq b_n^{\delta}\bigr\}.
\]
Then for sufficiently large $n$, we obtain
\begin{eqnarray*}
P\bigl[\mathsf{A}^c\bigr]  &\leq &   \sum_{k=1}^{[Tb_n]+1}
\sum_{i=1}^2P\bigl[\mathbf{N}^i_{(b_n^{-1}k) \wedge T}-
\mathbf{N}^i_{b_n^{-1}(k-1)}>b_n^{\delta}\bigr]
\\
&\leq &  b_n^{-q\delta}\sum_{k=1}^{[Tb_n]+1}
\sum_{i=1}^2E\bigl[\bigl(
\mathbf{N}^i_{(b_n^{-1}k) \wedge T}-\mathbf{N}^i_{b_n^{-1}(k-1)}
\bigr)^q\bigr]\leq Cb_n^{1-q\delta}.
\end{eqnarray*}

On the other hand, for any $k\in\mathbb{Z}_+$,
\[
\bigl|S^{n,j_2}-S^{n,j_1}\bigr|\leq kb_n^{-1}
\quad \Rightarrow \quad |j_2-j_1|\leq (k+1)b_n^{\delta}
\]
on $\mathsf{A}$. Hence, we have
\[
|j_2-j_1|>(k+1)b_n^{\delta}
\quad \Rightarrow \quad  \bigl|S^{n,j_2}-S^{n,j_1}\bigr|>kb_n^{-1}\qquad  \mbox{on }  \mathsf{A}.
\]

For sufficiently large $n$, if $|j_2-j_1|\geq b_n^{\delta_2}$ and
$\omega\in\mathsf{A}$, there exists $k\in\mathbb{N}$ such that
$(k+1)b_n^{\delta}<|j_2-j_1|\leq(k+2)b_n^{\delta}$.
Then since $|S^{n,j_2}-S^{n,j_1}|>kb_n^{-1}$, we have
\[
\frac{|S^{n,j_2}-S^{n,j_1}|}{|j_2-j_1|}>\frac
{kb_n^{-1}}{(k+2)b_n^{\delta}}\geq\frac{1}{3}b_n^{-1-\delta}
\geq b_n^{-1-\delta_3}.
\]
Therefore, we obtain
\[
b_n^2\sup_{j_1,j_2\in\mathbb{N},|j_2-j_1|\geq b_n^{\delta_2}}P \biggl[\ell
_{1,n}\geq j_1\vee j_2 \mbox{ and }
\frac
{|S^{n,j_2}-S^{n,j_1}|}{|j_2-j_1|}\leq b_n^{-1-\delta_3} \biggr] \leq
b_n^2P\bigl[\mathsf{A}^c\bigr]\leq
Cb_n^{3-q\delta} \to0
\]
as $n\to\infty$.
Similar estimates for $\{T^{n,j}\}$ hold true.

In particular, under $[B1]$, Proposition~8 in Ogihara and Yoshida \cite{ogi-yos} yields\break  $\limsup_{n\to\infty}E[b_n^{q-1}r_n^q]<\infty$ for any $q>0$.
Then we have $[A2]$.
\end{pf*}

\end{appendix}





\printhistory

\begin{thebibliography}{35}


\bibitem{ada}
\begin{bbook}[mr]
\bauthor{\bsnm{Adams},~\bfnm{Robert~A.}\binits{R.A.}}
(\byear{1975}).
\btitle{Sobolev Spaces}.
\bseries{Pure and Applied Mathematics}
\bvolume{65}.
\blocation{New York}:
\bpublisher{Academic Press}.
\bid{mr={0450957}}
\end{bbook}
\bptok{imsref}%
\endbibitem

\bibitem{ada-fou}
\begin{bbook}[mr]
\bauthor{\bsnm{Adams},~\bfnm{Robert~A.}\binits{R.A.}} \AND
\bauthor{\bsnm{Fournier},~\bfnm{John~J.~F.}\binits{J.J.F.}}
(\byear{2003}).
\btitle{Sobolev Spaces},
\bedition{2nd} ed.
\bseries{Pure and Applied Mathematics (Amsterdam)}
\bvolume{140}.
\blocation{Amsterdam}:
\bpublisher{Elsevier/Academic Press}.
\bid{mr={2424078}}
\end{bbook}
\bptok{imsref}%
\endbibitem

\bibitem{ait-et-al}
\begin{barticle}[mr]
\bauthor{\bsnm{A{\"{\i}}t-Sahalia},~\bfnm{Yacine}\binits{Y.}},
\bauthor{\bsnm{Fan},~\bfnm{Jianqing}\binits{J.}} \AND
\bauthor{\bsnm{Xiu},~\bfnm{Dacheng}\binits{D.}}
(\byear{2010}).
\btitle{High-frequency covariance estimates with noisy and asynchronous financial data}.
\bjournal{J. Amer. Statist. Assoc.}
\bvolume{105}
\bpages{1504--1517}.
\bid{doi={10.1198/jasa.2010.tm10163}, issn={0162-1459}, mr={2796567}}
\end{barticle}
\bptok{imsref}%
\endbibitem

\bibitem{ald-eag}
\begin{barticle}[mr]
\bauthor{\bsnm{Aldous},~\bfnm{D.~J.}\binits{D.J.}} \AND
\bauthor{\bsnm{Eagleson},~\bfnm{G.~K.}\binits{G.K.}}
(\byear{1978}).
\btitle{On mixing and stability of limit theorems}.
\bjournal{Ann. Probab.}
\bvolume{6}
\bpages{325--331}.
\bid{mr={0517416}}
\end{barticle}
\bptok{imsref}%
\endbibitem

\bibitem{aro}
\begin{barticle}[mr]
\bauthor{\bsnm{Aronson},~\bfnm{D.~G.}\binits{D.G.}}
(\byear{1967}).
\btitle{Bounds for the fundamental solution of a parabolic equation}.
\bjournal{Bull. Amer. Math. Soc.}
\bvolume{73}
\bpages{890--896}.
\bid{issn={0002-9904}, mr={0217444}}
\end{barticle}
\bptok{imsref}%
\endbibitem


\bibitem{bar-et-al}
\begin{barticle}[mr]
\bauthor{\bsnm{Barndorff-Nielsen},~\bfnm{Ole~E.}\binits{O.E.}},
\bauthor{\bsnm{Hansen},~\bfnm{Peter~Reinhard}\binits{P.R.}},
\bauthor{\bsnm{Lunde},~\bfnm{Asger}\binits{A.}} \AND
\bauthor{\bsnm{Shephard},~\bfnm{Neil}\binits{N.}}
(\byear{2011}).
\btitle{Multivariate realised kernels: Consistent positive semi-definite estimators of the covariation of equity prices with noise and non-synchronous trading}.
\bjournal{J. Econometrics}
\bvolume{162}
\bpages{149--169}.
\bid{doi={10.1016/j.jeconom.2010.07.009}, issn={0304-4076}, mr={2795610}}
\end{barticle}
\bptok{imsref}%
\endbibitem



\bibitem{bib-et-al}
\begin{barticle}[mr]
\bauthor{\bsnm{Bibinger},~\bfnm{Markus}\binits{M.}},
\bauthor{\bsnm{Hautsch},~\bfnm{Nikolaus}\binits{N.}},
\bauthor{\bsnm{Malec},~\bfnm{Peter}\binits{P.}} \AND
\bauthor{\bsnm{Rei{\ss}},~\bfnm{Markus}\binits{M.}}
(\byear{2014}).
\btitle{Estimating the quadratic covariation matrix from noisy observations: {L}ocal method of moments and efficiency}.
\bjournal{Ann. Statist.}
\bvolume{42}
\bpages{1312--1346}.
\bid{doi={10.1214/14-AOS1224}, issn={0090-5364}, mr={3226158}}
\bptnote{check year}%
\end{barticle}
\bptok{imsref}%
\endbibitem

\bibitem{chr-et-al}
\begin{barticle}[mr]
\bauthor{\bsnm{Christensen},~\bfnm{Kim}\binits{K.}},
\bauthor{\bsnm{Kinnebrock},~\bfnm{Silja}\binits{S.}} \AND
\bauthor{\bsnm{Podolskij},~\bfnm{Mark}\binits{M.}}
(\byear{2010}).
\btitle{Pre-averaging estimators of the ex-post covariance matrix in noisy diffusion models with non-synchronous data}.
\bjournal{J. Econometrics}
\bvolume{159}
\bpages{116--133}.
\bid{doi={10.1016/j.jeconom.2010.05.001}, issn={0304-4076}, mr={2720847}}
\end{barticle}
\bptok{imsref}%
\endbibitem


\bibitem{doh}
\begin{barticle}[mr]
\bauthor{\bsnm{Dohnal},~\bfnm{Gejza}\binits{G.}}
(\byear{1987}).
\btitle{On estimating the diffusion coefficient}.
\bjournal{J. Appl. Probab.}
\bvolume{24}
\bpages{105--114}.
\bid{issn={0021-9002}, mr={0876173}}
\end{barticle}
\bptok{imsref}%
\endbibitem


\bibitem{gob}
\begin{barticle}[mr]
\bauthor{\bsnm{Gobet},~\bfnm{Emmanuel}\binits{E.}}
(\byear{2001}).
\btitle{Local asymptotic mixed normality property for elliptic diffusion: A {M}alliavin calculus approach}.
\bjournal{Bernoulli}
\bvolume{7}
\bpages{899--912}.
\bid{doi={10.2307/3318625}, issn={1350-7265}, mr={1873834}}
\end{barticle}
\bptok{imsref}%
\endbibitem

\bibitem{gob02}
\begin{barticle}[mr]
\bauthor{\bsnm{Gobet},~\bfnm{Emmanuel}\binits{E.}}
(\byear{2002}).
\btitle{L{AN} property for ergodic diffusions with discrete observations}.
\bjournal{Ann. Inst. Henri Poincar\'e Probab. Stat.}
\bvolume{38}
\bpages{711--737}.
\bid{doi={10.1016/S0246-0203(02)01107-X}, issn={0246-0203}, mr={1931584}}
\end{barticle}
\bptok{imsref}%
\endbibitem

\bibitem{hay-yos01}
\begin{barticle}[mr]
\bauthor{\bsnm{Hayashi},~\bfnm{Takaki}\binits{T.}} \AND
\bauthor{\bsnm{Yoshida},~\bfnm{Nakahiro}\binits{N.}}
(\byear{2005}).
\btitle{On covariance estimation of non-synchronously observed diffusion processes}.
\bjournal{Bernoulli}
\bvolume{11}
\bpages{359--379}.
\bid{doi={10.3150/bj/1116340299}, issn={1350-7265}, mr={2132731}}
\end{barticle}
\bptok{imsref}%
\endbibitem





\bibitem{ibr-has03}
\begin{bbook}[mr]
\bauthor{\bsnm{Ibragimov},~\bfnm{I.~A.}\binits{I.A.}} \AND
\bauthor{\bsnm{Has'minski{\u\i}},~\bfnm{R.~Z.}\binits{R.Z.}}
(\byear{1981}).
\btitle{Statistical Estimation: Asymptotic Theory}.
\bseries{Applications of Mathematics}
\bvolume{16}.
\blocation{New York}:
\bpublisher{Springer}.
\bid{mr={0620321}}
\end{bbook}
\bptok{imsref}%
\endbibitem

\bibitem{jac}
\begin{bincollection}[mr]
\bauthor{\bsnm{Jacod},~\bfnm{Jean}\binits{J.}}
(\byear{1997}).
\btitle{On continuous conditional {G}aussian martingales and stable convergence in law}.
In \bbooktitle{S\'eminaire de {P}robabilit\'es, {XXXI}}.
\bseries{Lecture Notes in Math.}
\bvolume{1655}
\bpages{232--246}.
\blocation{Berlin}:
\bpublisher{Springer}.
\bid{doi={10.1007/BFb0119308}, mr={1478732}}
\end{bincollection}
\bptok{imsref}%
\endbibitem


\bibitem{jeg}
\begin{barticle}[mr]
\bauthor{\bsnm{Jeganathan},~\bfnm{P.}\binits{P.}}
(\byear{1983}).
\btitle{Some asymptotic properties of risk functions when the limit of the experiment is mixed normal}.
\bjournal{Sankhy\=a Ser. A}
\bvolume{45}
\bpages{66--87}.
\bid{issn={0581-572X}, mr={0749355}}
\end{barticle}
\bptok{imsref}%
\endbibitem

\bibitem{mal-man}
\begin{barticle}[mr]
\bauthor{\bsnm{Malliavin},~\bfnm{Paul}\binits{P.}} \AND
\bauthor{\bsnm{Mancino},~\bfnm{Maria~Elvira}\binits{M.E.}}
(\byear{2002}).
\btitle{Fourier series method for measurement of multivariate volatilities}.
\bjournal{Finance Stoch.}
\bvolume{6}
\bpages{49--61}.
\bid{doi={10.1007/s780-002-8400-6}, issn={0949-2984}, mr={1885583}}
\end{barticle}
\bptok{imsref}%
\endbibitem

\bibitem{nua}
\begin{bbook}[mr]
\bauthor{\bsnm{Nualart},~\bfnm{David}\binits{D.}}
(\byear{2006}).
\btitle{The {M}alliavin Calculus and Related Topics},
\bedition{2nd} ed.
\bseries{Probability and Its Applications (New York)}.
\blocation{Berlin}:
\bpublisher{Springer}.
\bid{mr={2200233}}
\end{bbook}
\bptok{imsref}%
\endbibitem

\bibitem{nua-par}
\begin{barticle}[mr]
\bauthor{\bsnm{Nualart},~\bfnm{D.}\binits{D.}} \AND
\bauthor{\bsnm{Pardoux},~\bfnm{{\'E}.}\binits{{\'E}.}}
(\byear{1988}).
\btitle{Stochastic calculus with anticipating integrands}.
\bjournal{Probab. Theory Related Fields}
\bvolume{78}
\bpages{535--581}.
\bid{doi={10.1007/BF00353876}, issn={0178-8051}, mr={0950346}}
\end{barticle}
\bptok{imsref}%
\endbibitem

\bibitem{ogi-yos}
\begin{bmisc}[auto:STB|2014/08/04|07:23:14]
\bauthor{\bsnm{Ogihara},~\bfnm{T.}\binits{T.}} \AND
\bauthor{\bsnm{Yoshida},~\bfnm{N.}\binits{N.}}
(\byear{2012}).
\bhowpublished{Quasi-likelihood analysis for stochastic regression models with nonsynchronous observations. Available at \arxivurl{arXiv:1212.4911}}.
\end{bmisc}
\bptok{imsref}%
\endbibitem

\bibitem{ogi-yos02}
\begin{barticle}[mr]
\bauthor{\bsnm{Ogihara},~\bfnm{Teppei}\binits{T.}} \AND
\bauthor{\bsnm{Yoshida},~\bfnm{Nakahiro}\binits{N.}}
(\byear{2014}).
\btitle{Quasi-likelihood analysis for nonsynchronously observed diffusion processes}.
\bjournal{Stochastic Process. Appl.}
\bvolume{124}
\bpages{2954--3008}.
\bid{doi={10.1016/j.spa.2014.03.014}, issn={0304-4149}, mr={3217430}}
\end{barticle}
\bptok{imsref}%
\endbibitem

\bibitem{pro}
\begin{bbook}[mr]
\bauthor{\bsnm{Protter},~\bfnm{Philip}\binits{P.}}
(\byear{1990}).
\btitle{Stochastic Integration and Differential Equations: A New Approach}.
\bseries{Applications of Mathematics (New York)}
\bvolume{21}.
\blocation{Berlin}:
\bpublisher{Springer}.
\bid{doi={10.1007/978-3-662-02619-9}, mr={1037262}}
\end{bbook}
\bptok{imsref}%
\endbibitem


\bibitem{uch-yos02}
\begin{barticle}[mr]
\bauthor{\bsnm{Uchida},~\bfnm{Masayuki}\binits{M.}} \AND
\bauthor{\bsnm{Yoshida},~\bfnm{Nakahiro}\binits{N.}}
(\byear{2013}).
\btitle{Quasi likelihood analysis of volatility and nondegeneracy of statistical random field}.
\bjournal{Stochastic Process. Appl.}
\bvolume{123}
\bpages{2851--2876}.
\bid{doi={10.1016/j.spa.2013.04.008}, issn={0304-4149}, mr={3054548}}
\end{barticle}
\bptok{imsref}%
\endbibitem


\bibitem{yos05}
\begin{barticle}[mr]
\bauthor{\bsnm{Yoshida},~\bfnm{Nakahiro}\binits{N.}}
(\byear{2011}).
\btitle{Polynomial type large deviation inequalities and quasi-likelihood analysis for stochastic differential equations}.
\bjournal{Ann. Inst. Statist. Math.}
\bvolume{63}
\bpages{431--479}.
\bid{doi={10.1007/s10463-009-0263-z}, issn={0020-3157}, mr={2786943}}
\end{barticle}
\bptok{imsref}%
\endbibitem
\end{thebibliography}
\end{document}